\documentclass[10pt]{amsart}
\usepackage{amscd}
\usepackage{amssymb}

\usepackage{eepic}

\allowdisplaybreaks[1]

\newtheorem{thm}{Th\'eor\`eme}[section]
\newtheorem{rema}[thm]{Remarque}

\newtheorem{prop}[thm]{Proposition}

\newtheorem{defn}[thm]{D\'efinition}
\theoremstyle{definition}

\numberwithin{equation}{section}

\newcommand{\n}{\noindent}

\newcommand{\resumename}{R\'esum\'e}
\newenvironment{resume}{\narrower\footnotesize\bf
\noindent\resumename.\quad\footnotesize\rm}{\par\bigskip}

\font\hb=cmbx12
\newcommand\pt{\hbox{\hb .}}

\begin{document}

\title[alg\`ebres homotopiques]{Les $(a,b)$-alg\`ebres \`a homotopie pr\`es}
\author[Walid Aloulou ]{Walid ALOULOU}

\address{
D\'epartement de Math\'ematiques\\
Fa\-cult\'e des Sciences de Monastir\\
Avenue de l'environnement\\
5019 Monastir\\
Tuni\-sie\\}
\address{Institut de Math\'ematiques de Bourgogne\\
UMR CNRS 5584\\
Universit\'e de Bourgogne\\
U.F.R. Sciences et Techniques
B.P. 47870\\
F-21078 Dijon Cedex\\France} \email{Walid.Aloulou@ipeim.rnu.tn}

\thanks{
Ce travail a \'et\'e effectu\'e dans le cadre de l'accord CMCU 06 S 1502.
W. Aloulou remercie l'Universit\'e de Bourgogne pour l'accueil
dont il a b\'en\'efici\'e au cours de ses s\'ejours.}

\begin{abstract}
We study in this article the concepts of algebra up to homotopy for a structure defined by two operations $ \pt $ and $[~,~]$. Having determined the structure of $ G_\infty $ algebras and
$ P_\infty $ algebras, we generalize this construction and we define a
structure of $ (a, b)$-algebra up to homotopy. Given a structure of commutative and differential graded Lie algebra for two shifts degree given by $a$ and  $b$, we will give an explicit construction of the associate algebra up to homotopy.
\end{abstract}

\maketitle
\begin{resume}

On \'etudie dans cet article les notions d'alg\`ebre \`a homotopie pr\`es pour une structure d\'efinie par deux op\'erations $\pt$ et $[~,~]$. Ayant d\'etermin\'e la structure des $G_\infty$ alg\`ebres et
des $P_\infty$ alg\`ebres, on g\'en\'eralise cette construction et on d\'efinit la
stucture des $(a,b)$-alg\`ebres \`a homotopie pr\`es. Etant donn\'ee une structure d'alg\`ebre commutative et de Lie diff\'erentielle gradu\'ee pour deux d\'ecalages des degr\'es donn\'es par $a$ et $b$, on donnera une construction explicite de l'alg\`ebre \`a homotopie pr\`es associ\'ee.

\end{resume}


\

\section{Introduction}\label{sec1}

\

Consid\'erons une alg\`ebre $\mathcal A$ pour une op\'eration $m$ ($m$ est par exemple une multiplication associative ou un crochet de Lie) munie d'une diff\'erentielle $d$. On dira juste que $\mathcal A$ est une alg\`ebre de type $\mathcal P$. Dans beaucoup de cas, on sait d\'efinir la notion d'alg\`ebre de type $\mathcal P$ \`a homotopie pr\`es et construire l'alg\`ebre \`a homotopie pr\`es enveloppante de $\mathcal A$. Cette alg\`ebre donne naturellement les complexes d'homologie et de cohomologie associ\'es \`a ce type d'alg\`ebre pour $\mathcal A$ et ses modules.\vskip0.2cm

Lorsque $\mathcal A$ poss\`ede deux op\'erations avec des relations de compatibilit\'e, la construction de l'alg\`ebre \`a homotopie pr\`es enveloppante correpondante est plus difficile.
En particulier, le complexe de Poisson et celui de Gerstenhaber consiste \`a composer les structures de la cog\`ebre cocommutative colibre et de la cog\`ebre de Lie colibre associ\'ees \`a chaque alg\`ebre. On obtient une bicog\`ebre codiff\'erentielle colibre. Cette composition tient compte des degr\'es. Les structures de bicog\`ebre obtenues diff\`erent. Pour chaque type, on appelle la bicog\`ebre obtenue l'alg\`ebre \`a homotopie pr\`es, c'est la $P_\infty$ alg\`ebre pour une alg\`ebre de Poisson et la $G_\infty$ alg\`ebre pour une alg\`ebre de Gerstenhaber.\vskip0.2cm

Le pr\'esent travail consiste \`a unifier ces constructions d'alg\`ebre \`a homotopie pr\`es dans les cas des alg\`ebres de Poisson et des alg\`ebres de Gerstenhaber, ce qui nous permet de d\'efinir la structure d'une $(a,b)$-alg\`ebre \`a homotopie pr\`es.
Disons
qu'une $(a,b)$-alg\`ebre diff\'erentielle est un espace vectoriel gradu\'e
$\mathcal{A}$ muni d'un produit commutatif $\pt$ de degr\'e
$a\in\mathbb{Z}$ $(|\pt|=a)$, d'un crochet $[~,~]$ de degr\'e
$b\in\mathbb{Z}$ $(|[~,~]|=b)$ et d'une diff\'erentielle $d$ de degr\'e $1$ tel que
$\big(\mathcal{A}[-a],\pt,d\big)$ est une alg\`ebre commutative
associative diff\'erentielle gradu\'ee et $\big(\mathcal{A}[-b],[~,~],d\big)$ est une alg\`ebre de Lie diff\'erentielle
gradu\'ee. Le produit et le crochet v\'erifient une relation de
compatibilit\'e entre eux dite identit\'e de Leibniz donn\'ee
par:
$$
\forall\alpha,\beta,\gamma\in\mathcal{A}, \hskip0.8cm
[\alpha,\beta\pt\gamma]=[\alpha,\beta]\pt\gamma+(-1)^{(|\beta|+a)(|\alpha|+b)}\beta\pt[\alpha,\gamma].
$$

Dans le cas o\`u $a=0$ et $b=-1$, on retrouve les alg\`ebres de
Gerstenhaber et dans le cas o\`u $a=b=0$, on retrouve les
alg\`ebres de Poisson gradu\'ees.\vskip0.2cm

Dans cet article, on reviendra dans la deuxi\`eme section \`a la
construction des alg\`ebres de Gerstenhaber \`a homotopie pr\`es.
Dans la troisi\`eme section, on donnera la construction de la
structure des alg\`ebres de Poisson \`a homotopie pr\`es. Dans la
quatri\`eme section, on g\'en\'eralisera la construction des
alg\`ebres de Gerstenhaber et de Poisson \`a homotopie pr\`es et
on donnera la structure des $(a,b)$-alg\`ebres \`a homotopie
pr\`es. Dans la cinqui\`eme section, on donnera quelques exemples de $(a,b)$-alg\`ebres diff\'erentielles gradu\'ees en introduisant, en particulier, les super-crochets de Schouten et de Poisson pour un super-espace $\mathbb{R}^{p|q}$.


\section{Les alg\`ebres de Gerstenhaber \`a homotopie pr\`es}

\

\subsection{D\'efinitions}

\

Soit $V=\displaystyle\bigoplus_{i\in \mathbb Z}V_i$ un espace vectoriel gradu\'e. On note par $V[k]$ l'espace $V$ muni de la graduation $V[k]_i=V_{i+k}$.

Soit $W$ un autre espace gradu\'e. Une application lin\'eaire $f$ de degr\'e $k$ de $V$ dans $W$ sera not\'ee $f:V\longrightarrow W[k]$.\vskip0.25cm

Une alg\`ebre de Gerstenhaber est un espace vectoriel gradu\'e
$\mathcal G$ muni d'une multiplication commutative et associative
gradu\'ee $\pt: \mathcal G\otimes\mathcal
G\longrightarrow\mathcal G$ de degr\'e $0$ et d'un crochet
$[~~,~~]: \mathcal G\otimes\mathcal G\longrightarrow\mathcal G$ de
degr\'e $-1$ tel que $\left(\mathcal G, \pt\right)$ est une alg\`ebre commutative associative gradu\'ee et $\left(\mathcal G[1], [~~,~~]\right)$ est
une alg\`ebre de Lie gradu\'ee et que, pour tout $\alpha$
homog\`ene, l'application $[\alpha, ~.~]$ est une d\'erivation
gradu\'ee pour la multiplication $\pt$. En notant $|\alpha|$ le
degr\'e d'un \'el\'ement homog\`ene $\alpha$ de $\mathcal G$, on a
donc:
$$
\begin{aligned}
\alpha\pt\beta&=(-1)^{|\alpha||\beta|}\beta\pt\alpha,\\
\alpha\pt(\beta\pt\gamma)&=(\alpha\pt\beta)\pt\gamma,\\
[\alpha,\beta]&=-(-1)^{(|\alpha|-1)(|\beta|-1)}[\beta,\alpha],\\
(-1)^{(|\alpha|-1)(|\gamma|-1)}\big[[\alpha,\beta],\gamma\big]&
+(-1)^{(|\beta|-1)(|\alpha|-1)}\big[[\beta,\gamma],\alpha\big]
+(-1)^{(|\gamma|-1)(|\beta|-1)}\big[[\gamma,\alpha],\beta\big]=0,\\
[\alpha,\beta\pt\gamma]&=[\alpha,\beta]\pt\gamma+(-1)^{|\beta|(|\alpha|-1)}\beta\pt[\alpha,\gamma]
\end{aligned}
$$
et donc aussi:
$$
[\alpha\pt\beta,\gamma]=\alpha\pt[\beta,\gamma]+(-1)^{|\beta|(|\gamma|-1)}[\alpha,\gamma]\pt\beta.
$$
Si de plus, on a une diff\'erentielle
$d:\mathcal{G}\longrightarrow\mathcal{G}[1]$ (ou $d:\mathcal{G}[1]\longrightarrow\mathcal{G}[2]$) de degr\'e $1$ telle que
$d\circ d=0$,
$d(\alpha\pt\beta)=d\alpha\pt\beta+(-1)^{|\alpha|}\alpha\pt
d\beta $ \ et \
$d([\alpha,\beta])=[d\alpha,\beta]+(-1)^{|\alpha|-1}[\alpha,d\beta]
$, on dit que $\big(\mathcal{G},\pt, [~,~],d\big)$ est une
alg\`ebre de Gerstenhaber diff\'erentielle gradu\'ee. (Voir \cite{[G]})

\vskip0.15cm

Le degr\'e sur
$\mathcal{G}[1]$ sera not\'e $dg$, d\'efini par:
$dg(\alpha)=|\alpha|-1$. Alors, le degr\'e du
crochet devient $0$ et celui du produit vaut $1$ ($dg([~,~])=0$,
$dg(\pt)=1$).

Le produit $\pt$ \'etait commutatif et associatif pour $|~|$
mais il ne l'est plus pour $dg$.

Pour chaque $k$, l'espace $\bigwedge^k(\mathcal G[1])$ est lin\'eairement isomorphe \`a l'espace $S^k(\mathcal G)[k]$. Cet isomorphisme n'est pas canonique. Nous choisissons ici l'isomorphisme d\'efini dans \cite{[AAC1]}:
$$\alpha_1\wedge\dots\wedge\alpha_k\longmapsto (-1)^{\sum_{i=1}^k(k-i)dg(\alpha_i)}\alpha_1\dots\alpha_k.$$
On transforme donc le produit $\pt$ en l'application $\mu$ sur
$\mathcal{G}[1]$ d\'efini par:
$\mu(\alpha,\beta)=(-1)^{dg(\alpha)}\alpha\pt\beta$. Le produit
$\mu$ est anticommutatif et antiassociatif de degr\'e
$1$:
\begin{align*}\mu(\alpha,\beta)&=-(-1)^{dg(\alpha)dg(\beta)}\mu(\beta,\alpha)\\\mu\left(\mu(\alpha,\beta),\gamma\right)&=-(-1)^{dg(\alpha)}\mu\left(\alpha,\mu(\beta,\gamma)\right)\end{align*}
et
$$d(\mu(\alpha,\beta))=-\mu(d\alpha,\beta)+(-1)^{dg(\alpha)+1}\mu(\alpha,
d\beta).$$

Le crochet $[~,~]$ est antisym\'etrique de degr\'e $0$ dans
$\mathcal{G}[1]$ v\'erifiant Jacobi et Leibniz:
\begin{align*}
[\alpha,\beta]&=-(-1)^{dg(\alpha)dg(\beta)}[\beta,\alpha],\\
(-1)^{dg(\alpha)dg(\gamma)}\big[[\alpha,\beta],\gamma\big]+&(-1)^{dg(\beta)dg(\alpha)}\big[[\beta,\gamma],\alpha\big]
+(-1)^{dg(\gamma)dg(\beta)}\big[[\gamma,\alpha],\beta\big]=0,\\
[\alpha,\mu(\beta,\gamma)]&=(-1)^{dg(\alpha)}\mu([\alpha,\beta],\gamma)+(-1)^{dg(\alpha)(dg(\beta)+1)}\mu(\beta,[\alpha,\gamma])
\end{align*}
ou encore
$$[\mu(\alpha,\beta),\gamma)]=\mu(\alpha,[\beta,\gamma])+(-1)^{dg(\beta)dg(\gamma)}
\mu([\alpha,\gamma],\beta).$$

De plus, on a
$$d([\alpha,\beta])=[d\alpha,\beta]+(-1)^{dg(\alpha)}[\alpha,d\beta]
.$$

\subsection{Produit battement}

\

Soit $V$ un espace vectoriel gradu\'e. Notons ici le degr\'e d'un vecteur homog\`ene $\alpha_i$ de $V$ par la m\^eme lettre $\alpha_i$ et
$\varepsilon_\alpha(\sigma^{-1})=\varepsilon\left(\begin{smallmatrix}\alpha_1,\dots,\alpha_{p+q}\\\alpha_{\sigma^{-1}(1)},
\dots,\alpha_{\sigma^{-1}(p+q)}\end{smallmatrix}\right)$ la signature de la permutation $\sigma^{-1}$ en tenant compte des degr\'es, c'est \`a dire, la signature de la restriction de $\sigma^{-1}$  aux indices des vecteurs de degr\'e impair.

Rappelons qu'un $(p,q)$-battement $(p,q\geq1)$ est une permutation
$\sigma\in S_{p+q}$ telle que:
$$\sigma(1)<\dots<\sigma(p) \ \ \hbox{et} \ \
\sigma(p+1)<\dots<\sigma(p+q).$$

On appelle $Bat(p,q)=\{\sigma\in S_{p+q}/ \sigma \ \hbox{est un} \
(p,q)\hbox{-battement}\}$ l'ensemble de tous les
$(p,q)$-battements. Pour deux tenseurs
$\alpha=\alpha_1\otimes\dots\otimes\alpha_p$ et
$\beta=\alpha_{p+1}\otimes\dots\otimes\alpha_{p+q}$, on d\'efinit
le produit battement de $\alpha$ et $\beta$ par:
\begin{align*}
bat_{p,q}(\alpha,\beta)&=\sum_{\sigma\in
Bat(p,q)}\varepsilon_\alpha(\sigma^{-1})
\alpha_{\sigma^{-1}(1)}\otimes\dots\otimes
\alpha_{\sigma^{-1}(p+q)}\cr&=\sum_{\sigma\in
Bat(p,q)}\varepsilon\left(\begin{smallmatrix}\alpha_1,\dots,\alpha_{p+q}\\\alpha_{\sigma^{-1}(1)},
\dots,\alpha_{\sigma^{-1}(p+q)}\end{smallmatrix}\right)
\alpha_{\sigma^{-1}(1)}\otimes\dots\otimes
\alpha_{\sigma^{-1}(p+q)}.
\end{align*}
Ceci repr\'esente la somme sign\'ee de tous les tenseurs
$\alpha_{i_1}\otimes\dots\otimes \alpha_{i_n}$ dans lesquels les
vecteurs $\alpha_1,\dots,\alpha_p$ et $\alpha_{p+1},\dots,\alpha_{p+q}$
apparaissent rang\'es dans leur ordre naturel.
(Voir \cite{[GH]} et \cite{[Lo]})\\

On note
$\underline{\displaystyle\bigotimes}^n(\mathcal{G}[1])$ l'espace quotient $\bigotimes^n(\mathcal{G}[1])\diagup_{\displaystyle\sum_{p+q=n\atop
p,q\geq1}Im(bat_{p,q})}$ de
$\bigotimes^n\mathcal{G}[1]$ par la somme de toutes les images des
applications lin\'eaires $bat_{p,q}$ telles que $p+q=n$.\vskip0.1cm

On pose
$$\mathcal{H}=\underline{\displaystyle\bigotimes}^+(\mathcal{G}[1])=\displaystyle\bigoplus_{n\geq1}
\left(\underline{\bigotimes}^n(\mathcal{G}[1])\right).$$

L'espace $\mathcal H$ est engendr\'e par les `paquets' $X$ qui sont les classes  des tenseurs $\alpha_1\otimes\dots\otimes\alpha_n$ dans le quotient.
On notera ces classes par
$X=\alpha_{[1,n]}=\alpha_1\underline{\otimes}\dots\underline{\otimes}\alpha_n$.

Leur degr\'e est:
$$x=dg(X)=dg(\alpha_1)+\dots+dg(\alpha_n).$$

\begin{prop}

\

Le produit battement est associatif et commutatif gradu\'e de
degr\'e 0: Pour tout $\alpha\in\bigotimes^p\mathcal{G}[1]$,
$\beta\in \bigotimes^q\mathcal{G}[1]$, $\gamma\in\bigotimes^r
\mathcal{G}[1]$,
$$
\begin{aligned}
(i)&\quad bat_{p,q}(\alpha,\beta)=(-1)^{dg(\alpha)dg(\beta)}bat_{q,p}(\beta,\alpha),\\
(ii)&\quad
bat_{p+q,r}(bat_{p,q}(\alpha,\beta),\gamma)=bat_{p,q+r}(\alpha,bat_{q,r}(\beta,\gamma)).
\end{aligned}
$$
\end{prop}

\subsection{Complexe de Gerstenhaber}

\

La th\'eorie des op\'erades dit que le dual d'une structure
d'alg\`ebre commutative est une structure de cog\`ebre de Lie. On
consid\`ere, alors, la cog\`ebre de Lie
$\big(\mathcal{H}=\underline{\displaystyle\bigotimes}^+(\mathcal{G}[1]),\delta)$
o\`u $\delta$ est le cocrochet sur
$\underline{\displaystyle\bigotimes}^+(\mathcal{G}[1])$ d\'efini
pour $X=\alpha_1\underline{\otimes}\dots\underline{\otimes}\alpha_n\in\mathcal{H}$ par
$$
\begin{aligned}
\delta(X)&=\sum_{j=1}^{n-1}\alpha_1\underline{\otimes}\dots\underline{\otimes}
\alpha_j\bigotimes \alpha_{j+1}
\underline{\otimes}\dots\underline{\otimes} \alpha_n
-\varepsilon\left(\begin{smallmatrix}\alpha_1\dots
\alpha_{j}&\alpha_{j+1}\dots \alpha_n\\ \alpha_{j+1}\dots
\alpha_n&\alpha_1\dots
\alpha_j\end{smallmatrix}\right)\alpha_{j+1}\underline{\otimes}\dots\underline{\otimes}
\alpha_n\bigotimes
\alpha_1\underline{\otimes}\dots\underline{\otimes}
\alpha_j\\&=\sum_{U\underline{\otimes} V=X\atop U,V\neq\emptyset} U\bigotimes
V-(-1)^{v u}V\bigotimes U.
\end{aligned}
$$

Le cocrochet $\delta$ est coantisym\'etrique de degr\'e $0$ et
v\'erifie l'identit\'e de coJacobi:\vskip0.15cm

{\bf (i)} $\tau\circ\delta=-\delta$,\vskip0.15cm

{\bf (ii)}
$(id^{\otimes3}+\tau_{12}\circ\tau_{23}+\tau_{23}\circ\tau_{12})\circ(\delta\otimes
id)\circ\delta=0$,\vskip0.15cm

\noindent o\`u $\tau$ est la volte d\'efinie pour
$X,Y\in\mathcal{H}$ par:
$\tau(X\bigotimes Y)=(-1)^{x
y}Y\bigotimes X$,
$\tau_{12}=\tau\otimes id$ et $\tau_{23}=id\otimes \tau$.\vskip
0.25cm

Gr\^ace \`a notre d\'ecalage, $\mu$ et $d$ sont de m\^eme degr\'e $1$ dans $\mathcal{G}[1]$. On peut maintenant les r\'eunir en une seule op\'eration homog\`ene sur $\mathcal{H}$.

On prolonge le produit $\mu$ et la diff\'erentielle $d$ \`a
$\mathcal{H}=\underline{\bigotimes}^+\mathcal{G}[1]$ comme des
cod\'erivations $\mu_1$ et $d_1$ du cocrochet $\delta$  de degr\'e
$1$ en posant:
$$
d_1(\alpha_1\underline\otimes \dots \underline\otimes
\alpha_{n})=\sum_{1\leq k\leq
n}(-1)^{\sum_{i<k}dg(\alpha_i)}\alpha_1\underline\otimes \dots
\underline\otimes d(\alpha_k)\underline\otimes
\dots\underline\otimes \alpha_{n}
$$et
$$
\mu_1(\alpha_1\underline\otimes \dots \underline\otimes
\alpha_{n})=\sum_{1\leq k<
n}(-1)^{\sum_{i<k}dg(\alpha_i)}\alpha_1\underline\otimes \dots
\underline\otimes \mu(\alpha_k,\alpha_{k+1})\underline\otimes
\dots\underline\otimes \alpha_{n}.
$$
Alors, $$(\mu_1\otimes
id+id\otimes\mu_1)\circ\delta=\delta\circ\mu_1,
\mu_1\circ\mu_1=0, (d_1\otimes id+id\otimes
d_1)\circ\delta=\delta\circ d_1 \ \hbox{et} \ d_1\circ d_1=0.$$

Rappelons qu'une cod\'erivation $D$ de la cog\`ebre $(\mathcal H, \delta)$ est enti\`erement d\'etermin\'ee par une suite d'application $D_r:\underline\bigotimes^r \mathcal G[1]\longrightarrow\mathcal G[2]$. Plus pr\'ecis\'ement, on a:

$$D(\alpha_1\underline{\otimes}\dots\underline{\otimes}
\alpha_n)=\displaystyle\sum_{1\leq r\leq n\atop 0\leq j\leq
n-r}(-1)^{\sum_{i\leq
j}dg(\alpha_i)}\alpha_1\underline{\otimes}\dots\underline{\otimes}
\alpha_j\underline\otimes
D_r(\alpha_{j+1}\underline{\otimes}\dots\underline{\otimes}
\alpha_{j+r})\underline{\otimes}
\alpha_{j+r+1}\underline\otimes\dots\underline{\otimes}
\alpha_n.$$

Si on pose $D_1=d$, $D_2=\mu$, $D_k=0$, pour $k\geq3$,
alors, $D=d_1+\mu_1$ est l'unique
cod\'erivation de $\delta$ de degr\'e $1$ qui prolonge $d$ et
$\mu$ \`a $\mathcal{H}$. Elle v\'erifie
$$D\circ D=0 \ \hbox{et} \ (D\otimes
id+id\otimes D)\circ\delta=\delta\circ D.$$

 On obtient que
$\big(\mathcal{H}=
 \underline{\bigotimes}^+(\mathcal{G}[1]),\delta,D=d_1+\mu_1\big)$ est une cog\`ebre de Lie
codiff\'erentielle, c'est \`a dire, une alg\`ebre commutative \`a
homotopie pr\`es ou une $C_\infty$-alg\`ebre.\vskip0.15cm

Etudions maintenant le crochet $[~,~]$ d\'efini sur $\mathcal{G}[1]$. D'abord, on le prolonge  \`a
$\mathcal{H}$ comme l'unique crochet de degr\'e $0$ compatible
avec le cocrochet $\delta$ v\'erifiant:
$$\delta\circ[~,~]=([~,~]\otimes id)\circ\big(\tau_{23}\circ(\delta\otimes id)
+id\otimes\delta\big)+(id\otimes[~,~])\circ\big(\delta\otimes
id+\tau_{12}\circ(id\otimes \delta)\big).$$ Pour
$X=\alpha_1\underline{\otimes}\dots\underline{\otimes} \alpha_p$
et $Y=\alpha_{p+1}\underline{\otimes}\dots\underline{\otimes}
\alpha_{p+q}$, ce prolongement est
donn\'e  par:
$$ [X,Y]=\hskip-0,8cm\displaystyle\sum_{\sigma\in Bat(p,q)\atop
k,\sigma^{-1}(k)\leq p<\sigma^{-1}(k+1)}\hskip-0,8cm
\varepsilon_\alpha(\sigma^{-1})
\alpha_{\sigma^{-1}(1)}\underline{\otimes}\dots\underline{\otimes}[\alpha_{\sigma^{-1}(k)},
\alpha_{\sigma^{-1}(k+1)}]\underline{\otimes}
\dots\underline{\otimes}\alpha_{\sigma^{-1}(p+q)}.$$ Ce crochet
est bien d\'efini sur $\mathcal{H}$. On a

\begin{prop}

\

L'espace $\mathcal H$, muni du crochet $[~,~]$ et de l'op\'erateur
$D$ est une alg\`ebre de Lie diff\'erentielle gradu\'ee: Pour tout
$X$, $Y$ et $Z$ de $\mathcal H$, on a:

\begin{itemize}
\item[(i)] \quad$[X,Y]=-(-1)^{xy}[Y,X]$,

\item[(ii)] \quad$(-1)^{xz}\left[[X,Y],Z\right]+(-1)^{yx}\left[[Y,Z],X\right]+(-1)^{zy}\left[[Z,X],Y\right]=0$,

\item[(iii)] \quad$D\left([X,Y]\right)=\left[D(X),Y\right]+(-1)^x\left[X,D(Y)\right]$.
\end{itemize}
\end{prop}

L'espace $\mathcal{H}$ est maintenant muni d'une op\'eration $D$ \`a un argument de degr\'e $1$ et d'un crochet $[~,~]$ \`a deux arguments de degr\'e $0$.

La th\'eorie des op\'erades dit aussi que le dual d'une structure
d'alg\`ebre de Lie est une structure de cog\`ebre cocommutative
coassociative. ($\mathcal H, [~,~],D)$ \'etant une alg\`ebre de
Lie diff\'erentielle gradu\'ee, on consid\`ere l'espace
$\mathcal{H}[1]$ et le degr\'e $dg'(X)=dg(X)-1:=x'$ et on construit la
cog\`ebre cocommutative coassociative
$(S^+(\mathcal{H}[1]),\Delta)$, o\`u
$S^+(\mathcal{H}[1])=\bigoplus_{n\geq1}S^n(\mathcal{H}[1])$ et
$\Delta$ est le coproduit d\'efini par: $\forall X_1\dots X_n\in
S^n(\mathcal{H}[1]) $,
\begin{align*}\Delta(X_1\dots X_n)=\sum_{I\cup J=\{1,\dots, n\}\atop \#I, \#J>0}
\varepsilon\left(\begin{smallmatrix}x_1'\dots x_n'\\
x_I' x_J'\end{smallmatrix}\right)X_I\bigotimes X_J.\end{align*}

On a not\'e par $X_I$ le produit $X_{i_1}\dots X_{i_r}$ si
$I=\{i_1<\dots<i_r\}$ et par $X_J$ le produit $X_{j_1}\dots
X_{j_{n-r}}$ si $J=\{j_1<\dots<j_{n-r}\}$.\vskip0.15cm

Le coproduit $\Delta$ est cocommutatif et coassociatif: $$\tau'\circ\Delta=\Delta \ \hbox{et} \ (\Delta\otimes
id)\circ\Delta=(id\otimes\Delta)\circ\Delta,$$ o\`u $\tau'$ est la volte dans $S^+(\mathcal H [1])$.
\vskip0.15cm

Le crochet $[~,~]$ \'etait antisym\'etrique sur $\mathcal{H}$ de
degr\'e $0$, il devient sur
$\mathcal{H}[1]$ un crochet $\ell_2(X,Y)=(-1)^{x}[X,Y]$ sym\'etrique, de degr\'e $1$ et
v\'erifie l'identit\'e de Jacobi gradu\'ee. De plus, $D$ est une
diff\'erentielle gradu\'ee de $\ell_2$ qui est aussi de degr\'e $1$. On a donc:

\begin{itemize}
\item[(i)] \quad$\ell_2(X,Y)=(-1)^{x'y'}\ell_2(Y,X)$,\vskip0.15cm

\item[(ii)] \quad$(-1)^{x'z'}\ell_2\left(\ell_2(X,Y),Z\right)+(-1)^{y'x'}\ell_2\left(\ell_2(Y,Z),X\right)
+(-1)^{z'y'}\ell_2\left(\ell_2(Z,X),Y\right)=0$,\vskip0.15cm

\item[(iii)] \quad$D\left(\ell_2(X,Y)\right)=-\ell_2\left(D(X),Y\right)+(-1)^{x'+1}\ell_2\left(X,D(Y)\right)$.
\end{itemize}\vskip0.12cm

On peut maintenant r\'eunir $\ell_2$ et $D$ en une seule op\'eration $Q$ sur $S^+(\mathcal{H}[1])$.\vskip0.12cm

On prolonge $\ell_2$ \`a $S^+(\mathcal{H}[1])$ comme une
cod\'erivation $\ell$ de $\Delta$ de degr\'e $1$:
$$\ell(X_1\dots
X_n)=\displaystyle\sum_{i<j}\varepsilon\left(\begin{smallmatrix}x_1'\dots x_n'\\
x_i' x_j'x_1'\dots \widehat{ij}\dots x_n'
\end{smallmatrix}\right)\ell_{2}
(X_i,X_j).X_1\dots\widehat{ij}\dots X_n.$$

On prolonge aussi $D$ \`a $S^+(\mathcal{H}[1])$ comme une
cod\'erivation $m$ de $\Delta$ toujours de degr\'e $1$:
$$m(X_1\dots
X_n)=\displaystyle\sum_{i=1}^{n}\varepsilon\left(\begin{smallmatrix}x_1'\dots x_n'\\
x_i' x_1'\dots \widehat{i}\dots x_n'
\end{smallmatrix}\right)D(X_i).X_1\dots\widehat{i}\dots
X_n.$$

En posant $Q_1=D$, $Q_2=\ell_2$, $Q_k=0$, ~si
$k\geq3$ et
$$Q(X_1\dots X_n)=\displaystyle\sum_{I\cup
J=\{1,\ldots,n\}\atop I\neq\emptyset}\varepsilon\left(\begin{smallmatrix}x_1'\dots x_n'\\
x_I' x_J'\end{smallmatrix}\right)Q_{\#I} (X_{I}).X_J.$$

Alors,
$Q=m+\ell$ est l'unique cod\'erivation de $\Delta$, de degr\'e $1$, prolongeant $D$ et $\ell_2$ \`a $S^+(\mathcal{H}[1])$. Elle v\'erifie
$$Q^2=0 \ \hbox{et} \ \left(Q\otimes id+id\otimes
Q\right)\circ\Delta=\Delta\circ Q.$$

Donc, $\left(S^+(\mathcal{H}[1]),\Delta,Q=\ell+m\right)$ est une
cog\`ebre cocommuative coassociative et codiff\'erentielle, c'est
\`a dire, c'est une $L_\infty$ alg\`ebre ou alg\`ebre de Lie \`a
homotopie pr\`es. (Voir \cite{[AMM]})\vskip0.24cm

D'autre part, $\big(\mathcal{H},\delta,D\big)$ est une cog\`ebre
de Lie codiff\'erentielle, on d\'efinit un cocrochet $\kappa$ sur
$\mathcal{H}[1]$ par:
$$\kappa(X)=\sum_{U\underline\otimes V=X\atop U,V\neq\emptyset}(-1)^{u'}\left( U\bigotimes
V+(-1)^{v' u'}V\bigotimes U\right),\hskip0.5cm \forall
X\in\mathcal{H}[1].$$
Puis $\kappa$ se prolonge \`a $S^+(\mathcal{H}[1])$ par:
$$\kappa(X_1\dots X_n)=\hskip-0.8cm\displaystyle\sum_{1\leq s\leq n\atop I\cup
J=\{1,\dots,n\}\setminus\{s\}}\hskip-0.8cm\varepsilon\left(\begin{smallmatrix}x_1'\dots x_n'\\
x_I'x_s' x_J'\end{smallmatrix}\right)\hskip-0.2cm\sum_{U_s\underline\otimes V_s=X_s\atop
U_s,V_s\neq\emptyset }(-1)^{x_I'+u_s'}\left(X_I. U_s\bigotimes
V_s . X_J+(-1)^{v_s'u_s'}X_I. V_s\bigotimes U_s. X_J\right).$$
Alors, $\kappa$ est un cocrochet cosym\'etrique sur $S^+(\mathcal{H}[1])$ de
degr\'e $1$ v\'erifiant l'identit\'e de coJacobi:
$$\tau'\circ\kappa=\kappa \ \hbox{et} \ \Big(id^{\otimes3}+\tau_{12}'\circ\tau_{23}'+\tau_{23}'\circ\tau_{12}'\Big)\circ(\kappa\otimes
id)\circ\kappa=0.
$$
De plus, $Q=\ell+m$ est une cod\'erivation de
$S^+(\mathcal{H}[1])$ pour $\kappa$ de degr\'e $1$. \vskip0.2cm

Alors,
$\left(S^+(\mathcal{H}[1]),\kappa,Q\right)$ est une
cog\`ebre de Lie codiff\'erentielle gradu\'ee, donc,
c'est encore une $C_\infty$ alg\`ebre.\vskip0.15cm

Le cocrochet $\kappa$ et le coproduit $\Delta$ v\'erifient, enfin,
l'identit\'e de
coLeibniz:$$(id\otimes\Delta)\circ\kappa=(\kappa\otimes
id)\circ\Delta+\tau_{12}'\circ(id\otimes\kappa)\circ\Delta.
$$ On
dit que $\left(S^+(\mathcal{H}[1]),\Delta,\kappa,Q\right)$
est une cog\`ebre de Gerstenhaber codiff\'erentielle gradu\'ee.

\begin{defn}

\

 Une $G_\infty$ alg\`ebre $\left(V,\Delta,\kappa,Q\right)$est une bicog\`ebre
 codiff\'erentielle colibre
 telle que $Q$ est une cod\'erivation pour le coproduit cocommutatif et coassociatif $\Delta$ et pour le cocrochet de Lie $\kappa$, de degr\'e $1$ et de carr\'e nul.

La $G_\infty$ alg\`ebre $\left(S^+
 (\underline{\displaystyle\bigotimes}^+(\mathcal{G}[1])[1]),\Delta,\kappa,Q=\ell+m\right)$ est appel\'ee la
$G_\infty$ alg\`ebre enveloppante de l'alg\`ebre de Gerstenhaber $\mathcal G$.
\end{defn}

(Voir \cite{[AAC2]} et \cite{[BGHHW]})

\vskip0.5cm


\

\section{Les alg\`ebres de Poisson \`a homotopie pr\`es}\label{sec3}\

\

\subsection{D\'efinition}

\

Une alg\`ebre de Poisson gradu\'ee est un espace vectoriel gradu\'e
$\mathcal P$ muni d'un produit commutatif et associatif gradu\'e
$\pt: \mathcal P\otimes\mathcal P\longrightarrow\mathcal P$ de
degr\'e $0$ et d'un crochet $\{~~,~~\}: \mathcal P\otimes\mathcal
P\longrightarrow\mathcal P$ de degr\'e $0$ tel que $\left(\mathcal
P, \pt\right)$ est une alg\`ebre commutative gradu\'ee, $\left(\mathcal
P, \{~~,~~\}\right)$ est une alg\`ebre de Lie gradu\'ee et que,
pour tout $\alpha$ homog\`ene, l'application $\{\alpha, ~.~\}$
soit une d\'erivation gradu\'ee pour le produit $\pt$. En notant
$|\alpha|$ le degr\'e d'un \'el\'ement homog\`ene $\alpha$ de
$\mathcal P$, on a donc:
$$\begin{aligned}
\alpha\pt\beta&=(-1)^{|\alpha||\beta|}\beta\pt\alpha,\\
\alpha\pt(\beta\pt\gamma)&=(\alpha\pt\beta)\pt\gamma,\\
\{\alpha,\beta\}&=-(-1)^{|\alpha||\beta|}\{\beta,\alpha\},\\
(-1)^{|\alpha||\gamma|}\big\{\{\alpha,\beta\},\gamma\big\}&+(-1)^{|\beta||\alpha|}\big\{\{\beta,\gamma\},\alpha\big\}
+(-1)^{|\gamma||\beta|}\big\{\{\gamma,\alpha\},\beta\big\}=0,\\
\{\alpha,\beta\pt\gamma\}&=\{\alpha,\beta\}\pt\gamma+(-1)^{|\beta||\alpha|}\beta\pt\{\alpha,\gamma\}
\end{aligned}$$
et donc aussi:
$$
\hskip-3cm\{\alpha\pt\beta,\gamma\}=\alpha\pt\{\beta,\gamma\}+(-1)^{|\beta||\gamma|}\{\alpha,\gamma\}\pt\beta.
$$

Si de plus, on a une diff\'erentielle
$d:\mathcal{P}\longrightarrow\mathcal{P}[1]$ de degr\'e $1$ telle que
$$d\circ d=0, \
d(\alpha\pt\beta)=d\alpha\pt\beta+(-1)^{|\alpha|}\alpha\pt d\beta
 \ \hbox{et} \
d(\{\alpha,\beta\})=\{d\alpha,\beta\}+(-1)^{|\alpha|}\{\alpha,d\beta\}
,$$ on dit que $\big(\mathcal{P},\pt, \{~,~\},d\big)$ est une
alg\`ebre de Poisson diff\'erentielle gradu\'ee. (Voir \cite{[F]})\vskip0.3cm

Pour construire la structure de l'alg\`ebre de Poisson \`a
homotopie pr\`es, on va reprendre la construction pr\'ec\'edente dans le cas des alg\`ebres de Poisson diff\'erentielles gradu\'ees.

On consid\`ere, donc, l'espace $\mathcal{P}[1]$ et la graduation
$dg(\alpha)=|\alpha|-1$. On transforme le produit $\pt$ en $\mu$ en posant
$\mu(\alpha,\beta)=(-1)^{dg(\alpha)}\alpha\pt\beta$ et le crochet $\{~,~\}$ en $[~,~]$ en posant $[\alpha,\beta]=(-1)^{dg(\alpha)}\{\alpha,\beta\}$.

Le produit $\mu$ devient anticommutatif,
antiassociatif de degr\'e $1$ et le crochet $[~,~]$ devient
sym\'etrique  de degr\'e $1$, v\'erifiant l'identit\'e de Jacobi
gradu\'ee et celle de Leibniz avec $\mu$:
\begin{align*}
dg(\mu)&=dg([~,~])=dg(d)=1,\\
\mu(\alpha,\beta)&=-(-1)^{dg(\alpha)dg(\beta)}\mu(\beta,\alpha),\\
\mu\left(\mu(\alpha,\beta),\gamma\right)&=-(-1)^{dg(\alpha)}\mu\left(\alpha,\mu(\beta,\gamma)\right),\\
[\alpha,\beta]&=-(-1)^{dg(\alpha)dg(\beta)}[\beta,\alpha],\\
(-1)^{dg(\alpha)dg(\gamma)}\big[[\alpha,\beta],\gamma\big]&+(-1)^{dg(\beta)dg(\alpha)}\big[[\beta,\gamma],\alpha\big]
+(-1)^{dg(\gamma)dg(\beta)}\big[[\gamma,\alpha],\beta\big]=0,\\
[\alpha,\mu(\beta,\gamma)]&=(-1)^{dg(\alpha)+1}\mu([\alpha,\beta],\gamma)+(-1)^{(dg(\alpha)+1)(dg(\beta)+1)}
\mu(\beta,[\alpha,\gamma]),
\end{align*}

ou encore
$$[\mu(\alpha,\beta),\gamma)]=(-1)^{dg(\alpha)+1}\mu(\alpha,[\beta,\gamma])+(-1)^{dg(\beta)dg(\gamma)+1}
\mu([\alpha,\gamma],\beta).$$

De plus, on a
$$d(\mu(\alpha,\beta))=-\mu(d\alpha,\beta)+(-1)^{dg(\alpha)+1}\mu(\alpha,d\beta)
\ \hbox{et} \
d([\alpha,\beta])=-[d\alpha,\beta]+(-1)^{dg(\alpha)+1}[\alpha,d\beta].
$$

\subsection{Complexe de Poisson}

\

Comme dans le cas d'une alg\`ebre de Gerstenhaber, on r\'eunit d'abord $\mu$ et $d$:

On consid\`ere l'espace
$\mathcal{H}=\underline{\displaystyle\bigotimes}^+(\mathcal{P}[1])=\displaystyle\bigoplus_{n\geq1}
\left(\underline{\bigotimes}^n(\mathcal{P}[1])\right)$ avec le degr\'e
$$dg(X)=x=dg(\alpha_1)+\dots+dg(\alpha_n), \ \hbox{pour} \ X=\alpha_1\underline{\otimes}\dots\underline{\otimes}\alpha_n\in\mathcal{H}.$$
On munit cet espace du cocrochet $\delta$ d\'efini par
$$\delta(X)=\displaystyle\sum_{U\underline\otimes V=X\atop U,V\neq\emptyset}
U\bigotimes V-(-1)^{v u}V\bigotimes U.$$

Alors, le complexe $(\mathcal H, \delta)$ est une cog\`ebre de
Lie.\vskip0.15cm

On prolonge $d$ et $\mu$ \`a $\mathcal{H}$ comme des
cod\'erivations $d_1$ et $\mu_1$ du cocrochet $\delta$  de degr\'e
$1$ en posant:
$$
d_1(\alpha_1\underline\otimes \dots \underline\otimes
\alpha_{n})=\sum_{1\leq k\leq
n}(-1)^{\sum_{i<k}dg(\alpha_i)}\alpha_1\underline\otimes \dots
\underline\otimes d(\alpha_k)\underline\otimes
\dots\underline\otimes \alpha_{n}
$$et
$$
\mu_1(\alpha_1\underline\otimes \dots \underline\otimes
\alpha_{n})=\sum_{1\leq k<
n}(-1)^{\sum_{i<k}dg(\alpha_i)}\alpha_1\underline\otimes \dots
\underline\otimes \mu(\alpha_k,\alpha_{k+1})\underline\otimes
\dots\underline\otimes \alpha_{n}.
$$
Alors, $$(\mu_1\otimes
id+id\otimes\mu_1)\circ\delta=\delta\circ\mu_1,\mu_1\circ\mu_1=0, (d_1\otimes id+id\otimes
d_1)\circ\delta=\delta\circ d_1 \ \hbox{et} \ d_1\circ d_1=0.$$

Or, on sait qu'une cod\'erivation $D$ de la cog\`ebre $(\mathcal H, \delta)$ est enti\`erement d\'etermin\'ee par une suite d'application $D_r:\underline\bigotimes^r \mathcal P[1]\longrightarrow\mathcal P[2]$. Plus pr\'ecis\'ement, on a:

$$D(\alpha_1\underline{\otimes}\dots\underline{\otimes}
\alpha_n)=\displaystyle\sum_{1\leq r\leq n\atop 0\leq j\leq
n-r}(-1)^{\sum_{i\leq
j}dg(\alpha_i)}\alpha_1\underline{\otimes}\dots\underline{\otimes}
\alpha_j\underline\otimes
D_r(\alpha_{j+1}\underline{\otimes}\dots\underline{\otimes}
\alpha_{j+r})\underline{\otimes}
\alpha_{j+r+1}\underline\otimes\dots\underline{\otimes}
\alpha_n.$$

Si on pose $D_1=d$, $D_2=\mu$, $D_k=0$, pour $k\geq3$,
alors, $D=d_1+\mu_1$ est l'unique
cod\'erivation de $\delta$ de degr\'e $1$ qui prolonge $d$ et
$\mu$ \`a $\mathcal{H}$. Elle v\'erifie
$$D\circ D=0 \ \hbox{et} \ (D\otimes
id+id\otimes D)\circ\delta=\delta\circ D.$$
Alors, $\big(\mathcal{H}=
 \underline{\bigotimes}^+(\mathcal{P}[1]),\delta,D=d_1+\mu_1\big)$ est une cog\`ebre de Lie
codiff\'erentielle, donc, c'est une
$C_\infty$-alg\`ebre.\vskip0.15cm

Bien qu'il ne soit pas antisym\'etrique, on prolonge le crochet $[~,~]$ d\'efini sur $\mathcal{P}[1]$ \`a
$\mathcal{H}$ comme l'unique "crochet" de degr\'e $1$ compatible
avec le cocrochet $\delta$, c'est \`a dire, v\'erifiant:
$$\delta\circ[~,~]=([~,~]\otimes id)\circ\big(\tau_{23}\circ(\delta\otimes id)
+id\otimes\delta\big)+(id\otimes[~,~])\circ\big(\delta\otimes
id+\tau_{12}\circ(id\otimes \delta)\big).$$ Ce prolongement est
donn\'e pour
$X=\alpha_1\underline{\otimes}\dots\underline{\otimes} \alpha_p$
et $Y=\alpha_{p+1}\underline{\otimes}\dots\underline{\otimes}
\alpha_{p+q}$ par:
$$ [X,Y]=\hskip-0,8cm\displaystyle\sum_{\sigma\in Bat(p,q)\atop
k,\sigma^{-1}(k)\leq p<\sigma^{-1}(k+1)}\hskip-0,8cm
\varepsilon_\alpha(\sigma^{-1})(-1)^{\sum_{i<k}dg(\alpha_{\sigma^{-1}(i)})}
\alpha_{\sigma^{-1}(1)}\underline{\otimes}\dots\underline{\otimes}[\alpha_{\sigma^{-1}(k)},
\alpha_{\sigma^{-1}(k+1)}]\underline{\otimes}
\dots\underline{\otimes}\alpha_{\sigma^{-1}(p+q)}.$$
Puisque $[~,~]$ n'\'etait pas antisym\'etrique sur $\mathcal{P}[1]$, son prolongement ne l'est pas non plus sur $\mathcal{H}$. On se ram\`ene \`a un vrai crochet de degr\'e $0$ en d\'ecalant $\mathcal{H}$ par $-1$.

On consid\`ere l'espace $\mathcal{H}[-1]$ muni de la graduation
$dg'(X)=x'=dg(X)+1$, pour $X\in\mathcal{H}[-1]$. On construit sur $\mathcal H[-1]$ le crochet
$\{X,Y\}=(-1)^{-x'}[X,Y]$. Alors, on a

\begin{prop}

\

L'espace $\mathcal H[-1]$, muni du crochet $\{~,~\}$ et de
l'op\'erateur $D$ est une alg\`ebre de Lie diff\'erentielle
gradu\'ee: Pour tout $X$, $Y$ et $Z$ de $\mathcal H[-1]$, on a:

\begin{itemize}
\item[(i)] \quad$\{X,Y\}=-(-1)^{x'y'}\{Y,X\}$,\vskip0.15cm

\item[(ii)] \quad$(-1)^{x'z'}\{\{X,Y\},Z\}+(-1)^{y'x'}\{\{Y,Z\},X\}+(-1)^{z'y'}\{\{Z,X\},Y\}=0$,\vskip0.15cm

\item[(iii)] \quad$D\left(\{X,Y\}\right)=\{D(X),Y\}+(-1)^{x'}\{X,D(Y)\}$.
\end{itemize}
\end{prop}

Maintenant, $\big(\mathcal{H}[-1], \{~,~\}, D\big)$ est une alg\`ebre de Lie diff\'erentielle gradu\'ee. On peut donc lui associer une cog\`ebre cocommutative coassociative $(S^+(\mathcal{H}[-1][1]),\Delta)$ avec une cod\'erivation $Q$ construite \`a partir de $\{~,~\}$ et $D$. Mais ici $\mathcal{H}[-1][1]=\mathcal{H}$ et on pr\'ef\`ere r\'eutiliser le crochet $[~,~]$ sur $\mathcal{H}$ d\'efini plus haut plut\^ot que $(-1)^x\{X,Y\}=-[X,Y]$.\vskip0.2cm

On prolonge donc $D$ et $[~,~]$ \`a
$S^+(\mathcal{H})$ comme des cod\'erivations $m$ et $\ell$ de
$\Delta$
 de degr\'e $1$ et de carr\'e nul en posant:
$$m(X_1\dots
X_n)=\displaystyle\sum_{i=1}^{n}\varepsilon\left(\begin{smallmatrix}x_1\dots x_n\\
x_i x_1\dots \widehat{i}\dots x_n
\end{smallmatrix}\right)D(X_i).X_1\dots\widehat{i}\dots
X_n$$
et
$$\ell(X_1\dots
X_n)=\displaystyle\sum_{i<j}\varepsilon\left(\begin{smallmatrix}x_1\dots x_n\\
x_i x_jx_1\dots \widehat{ij}\dots x_n
\end{smallmatrix}\right)[X_i,X_j].X_1\dots\widehat{ij}\dots X_n.$$

Puis, on pose $Q_1=D$, $Q_2=\ell_2$, $Q_k=0$, ~si $k\geq3$ et
$$Q(X_1\dots X_n)=\displaystyle\sum_{I\cup
J=\{1,\ldots,n\}\atop I\neq\emptyset}\varepsilon\left(\begin{smallmatrix}x_1\dots x_n\\
x_I x_J\end{smallmatrix}\right)Q_{ \#I} (X_{I}).X_J.$$

Alors,
$Q=m+\ell$ est l'unique cod\'erivation de $\Delta$, de degr\'e $1$, prolongeant $D$ et $[~,~]$ \`a $S^+(\mathcal{H})$.
Elle v\'erifie
$$Q^2=0 \ \hbox{et} \ \left(Q\otimes id+id\otimes
Q\right)\circ\Delta=\Delta\circ Q.$$

Donc, $\left(S^+(\mathcal{H}),\Delta,Q\right)$ est une
$L_\infty$ alg\`ebre.\vskip0.24cm

D'autre part, $(\mathcal{H},\delta,D)$ \'etant une $C_\infty$
alg\`ebre, on prolonge le cocrochet $\delta$ \`a
$S^+(\mathcal{H})$ par:
$$\delta(X_1\dots X_n)=\hskip-0.6cm\displaystyle\sum_{1\leq s\leq n\atop I\cup
J=\{1,\dots,n\}\setminus\{s\}}\hskip-0.6cm\varepsilon\left(\begin{smallmatrix}x_1\dots x_n\\
x_Ix_s x_J\end{smallmatrix}\right)\sum_{U_s\underline\otimes V_s=X_s\atop
U_s,V_s\neq\emptyset }\left(X_I. U_s\bigotimes V_s . X_J-(-1)^{v_s
u_s}X_I. V_s\bigotimes U_s. X_J\right).$$

Alors, $\delta$ est un cocrochet coantisym\'etrique sur
$S^+(\mathcal{H})$ de degr\'e $0$ v\'erifiant l'identit\'e de
coJacobi gradu\'ee. De plus, $Q=\ell+m$ est une cod\'erivation de
$S^+(\mathcal{H})$ pour $\delta$. Donc,
$\left(S^+(\mathcal{H}),\delta,Q\right)$ est une cog\`ebre
de Lie codiff\'erentielle gradu\'ee, c'est \`a dire, c'est encore une
$C_\infty$ alg\`ebre.\vskip0.24cm

De plus, le cocrochet $\delta$ et le coproduit $\Delta$ v\'erifient
l'identit\'e de
coLeibniz:$$(id\otimes\Delta)\circ\delta=(\delta\otimes
id)\circ\Delta+\tau_{12}\circ(id\otimes\delta)\circ\Delta.
$$ On
dit que $\left(S^+(\mathcal{H}),\Delta,\delta,Q\right)$ est
une cog\`ebre de Poisson codiff\'erentielle gradu\'ee.

\begin{defn}

\

 Une $P_\infty$ alg\`ebre $(V,\Delta,\delta,Q)$ est une bicog\`ebre
 codiff\'erentielle colibre
 telle que $Q$ est une cod\'erivation pour le coproduit cocommutatif et coassociatif $\Delta$ et pour le cocrochet de Lie $\delta$, de degr\'e $1$ et de carr\'e nul.

La $P_\infty$ alg\`ebre $\left(S^+
 (\underline{\displaystyle\bigotimes}^+\mathcal{P}[1]),\Delta,\delta,Q=\ell+m\right)$ est appel\'ee la
$P_\infty$ alg\`ebre enveloppante de l'alg\`ebre de Poisson $\mathcal P$.

\end{defn}

\vskip0.4cm


\

\section{Les $(a,b)$-alg\`ebres \`a homotopie pr\`es}\label{sec4}\

\

\subsection{D\'efinition}

\

Consid\'erons un espace vectoriel $\mathcal{A}$ gradu\'e. Le
degr\'e d'un \'el\'ement homog\`ene $\alpha$ de $\mathcal{A}$ est
not\'e $|\alpha|$. Soient $a,b\in \mathbb{Z}$, l'espace
$\mathcal{A}$ est muni d'un produit $\pt$ de degr\'e $a$ ($|\pt|=a$)
et d'un crochet $[~,~]$ de degr\'e $b$ ($|[~,~]|=b$)
tel que $\Big(\mathcal{A}[-a],\pt\Big)$ est une alg\`ebre
commutative et associative gradu\'ee  et $\Big(\mathcal{A}[-b],[~,~]\Big)$
est une alg\`ebre de Lie gradu\'ee. De plus, l'application lin\'eaire
$ad:\mathcal{A}[-b]\longrightarrow\mathcal{D}er\Big(\mathcal{A}[-a],
\pt\Big)$; $\alpha\longmapsto ad_{\alpha}$
 est telle que $ad_{\alpha}$ soit une d\'erivation gradu\'ee pour le produit
 $\pt$.

 On dit que $\Big(\mathcal{A},\pt,[~,~]\Big)$ est une
 $(a,b)$-alg\`ebre gradu\'ee. Pour tout $\alpha,\beta,\gamma\in\mathcal{A}$, on a les propri\'et\'es suivantes:\vskip0.1cm

 {\bf(i)} $
 \alpha\pt\beta=(-1)^{(|\alpha|+a)(|\beta|+a)}\beta\pt\alpha$,\vskip0.15cm

{\bf(ii)}
$\alpha\pt(\beta\pt\gamma)=(\alpha\pt\beta)\pt\gamma$,\vskip0.15cm

{\bf(iii)}
$[\alpha,\beta]=-(-1)^{(|\alpha|+b)(|\beta|+b)}[\beta,\alpha]$,\vskip0.15cm

{\bf(iv)}
$(-1)^{(|\alpha|+b)(|\gamma|+b)}[[\alpha,\beta],\gamma]+(-1)^{(|\beta|+b)(|\alpha|+b)}\big[[\beta,\gamma]
,\alpha\big]+(-1)^{(|\gamma|+b)(|\beta|+b)}\big[[\gamma,\alpha],\beta\big]=0$,\vskip0.15cm

{\bf (v)} $[\alpha,\beta\pt\gamma]=[\alpha,\beta]\pt
\gamma+(-1)^{(|\beta|+a)(|\alpha|+b)}\beta\pt[\alpha,\gamma] $

qui s'\'ecrit encore $[\alpha\pt\beta,\gamma]=\alpha\pt[\beta,
\gamma]+(-1)^{(|\beta|+a)(|\gamma|+b)}[\alpha,\gamma]\pt\beta $.
\vskip0,3cm

\n De plus, si on a une
diff\'erentielle $d:\mathcal{A}[-a]\longrightarrow\mathcal{A}[-a+1]$ \Big(ou $d:\mathcal{A}[-b]\longrightarrow\mathcal{A}[-b+1]$\Big) de degr\'e $1$ v\'erifiant
$$d\circ d=0, \ d(\alpha\pt
\beta)=d\alpha\pt\beta+(-1)^{|\alpha|+a}\alpha\pt d\beta \ \hbox{et} \ d([\alpha,\beta])=[d\alpha,\beta]+(-1)^{|\alpha|+b}[\alpha,d\beta]
,$$ on dira que $\Big(\mathcal{A},\pt,[~,~],d\Big)$ est une
 $(a,b)$-alg\`ebre diff\'erentielle gradu\'ee.\vskip0,3cm

Comme ci-dessus, on utilise un d\'ecalage pour homog\`en\'eiser le produit et la diff\'erentielle.
 On consid\`ere l'espace $\mathcal{A}[-a+1]$ muni de la graduation
 $dg(\alpha)=|\alpha|+a-1$ que l'on note simplement par $\alpha$. Sur $\mathcal{A}[-a+1]$, le produit
 $\pt$ n'est plus commutatif et le crochet $[~,~]$ n'est plus
 antisym\'etrique. On construit, donc, un nouveau produit $\mu$ sur $\mathcal{A}[-a+1]=\mathcal{A}[-a][1]$ de degr\'e $1$ d\'efini
 par
 $$\mu(\alpha,\beta)=(-1)^{1.\alpha}\alpha\pt\beta$$
et un nouveau crochet $\ell$ sur
$\mathcal{A}[-a+1]=\mathcal{A}[-b][b-a+1]$ de degr\'e $b-a+1$ d\'efini par
$$\ell(\alpha,\beta)=(-1)^{(b-a+1).\alpha}[\alpha,\beta].$$ Et on a

{\bf(i)} $\mu(\alpha,\beta)=-(-1)^{\alpha
\beta}\mu(\beta,\alpha)$,\vskip0.15cm

{\bf(ii)}
$\mu\left(\mu(\alpha,\beta),\gamma\right)=-(-1)^{\alpha}\mu\left(\alpha,\mu(\beta,\gamma)\right)$,\vskip0.15cm

{\bf(iii)} $\ell(\alpha,\beta)=-(-1)^{b-a+1}(-1)^{\alpha
\beta}\ell(\beta,\alpha)$,\vskip0.15cm

{\bf(iv)} $(-1)^{\alpha
\gamma}\ell(\ell(\alpha,\beta),\gamma\big)+(-1)^{\beta
\alpha}\ell(\ell(\beta,\gamma),\alpha) +(-1)^{\gamma
\beta}\ell(\ell(\gamma,\alpha),\beta)=0$,\vskip0.15cm

{\bf(v)}
$\ell(\alpha,\mu(\beta,\gamma))=(-1)^{\alpha+b-a+1}\mu(\ell(\alpha,\beta),\gamma)+
(-1)^{(\alpha+b-a+1)(\beta+1)} \mu(\beta,\ell(\alpha,\gamma))
$,\vskip0.15cm

ou encore
$\ell(\alpha,\mu(\beta,\gamma))=(-1)^{(b-a+1)(\alpha+1)}\mu(\alpha,\ell(\beta,
\gamma))+(-1)^{b-a+1+\beta \gamma}\mu(\ell(\alpha,\gamma),\beta)
$. \vskip0,2cm

\noindent De plus, $d$ reste encore une d\'erivation pour $\mu$ et $\ell$, elle v\'erifie:\vskip0.2cm

\noindent$ d(\mu(\alpha,
\beta))=-\mu(d\alpha,\beta)+(-1)^{\alpha+1}\mu(\alpha, d\beta) $
et
$d(\ell(\alpha,\beta))=(-1)^{b-a+1}\ell(d\alpha,\beta)+(-1)^{\alpha+b-a+1}\ell(\alpha,d\beta)
$.

\subsection{Extension de la multiplication et du crochet \`a la cog\`ebre de Lie codiff\'erentielle}

\

On consid\`ere comme pr\'ec\'edemment l'espace
 $\mathcal{H}=\displaystyle\bigoplus_{n\geq1}\big(\underline{\bigotimes}^n\mathcal{A}[-a+1]\big)=
 \underline{\bigotimes}^+(\mathcal{A}[-a+1])$ et pour
$X=\alpha_1\underline{\otimes}\dots\underline{\otimes}\alpha_n\in\mathcal{H}$,
le degr\'e $dg(X)=\alpha_1+\dots+\alpha_n$ not\'e simplement par
$x$. Sur cet espace, on d\'efinit un cocrochet $\delta$ de degr\'e
$0$  par:
$$
\begin{aligned}
\delta(X)&=\sum_{j=1}^{n-1}\alpha_1\underline{\otimes}\dots\underline{\otimes}
\alpha_j\bigotimes \alpha_{j+1}
\underline{\otimes}\dots\underline{\otimes} \alpha_n
-\varepsilon\left(\begin{smallmatrix}\alpha_1\dots
\alpha_{j}&\alpha_{j+1}\dots \alpha_n\\ \alpha_{j+1}\dots
\alpha_n&\alpha_1\dots
\alpha_j\end{smallmatrix}\right)\alpha_{j+1}\underline{\otimes}\dots\underline{\otimes}
\alpha_n\bigotimes
\alpha_1\underline{\otimes}\dots\underline{\otimes}
\alpha_j\\&=\sum_{U\underline\otimes V=X\atop U,V\neq\emptyset} U\bigotimes
V-(-1)^{v u}V\bigotimes U.
\end{aligned}
$$

On prolonge $\mu$ et $d$ \`a $\mathcal{H}$ comme des
cod\'erivations $\mu_1$ et $d_1$ de $\delta$ de degr\'e $1$ en
posant:
$$
d_1(\alpha_1\underline\otimes \dots \underline\otimes
\alpha_{n})=\sum_{1\leq k\leq
n}(-1)^{\sum_{i<k}\alpha_i}\alpha_1\underline\otimes \dots
\underline\otimes d(\alpha_k)\underline\otimes
\dots\underline\otimes \alpha_{n}
$$et
$$
\mu_1(\alpha_1\underline\otimes \dots \underline\otimes
\alpha_{n})=\sum_{1\leq
k<n}(-1)^{\sum_{i<k}\alpha_i}\alpha_1\underline\otimes \dots
\underline\otimes \mu(\alpha_k,\alpha_{k+1})\underline\otimes
\dots\underline\otimes \alpha_{n}.
$$
Alors, $$(\mu_1\otimes
id+id\otimes\mu_1)\circ\delta=\delta\circ\mu_1,
\mu_1\circ\mu_1=0, (d_1\otimes id+id\otimes
d_1)\circ\delta=\delta\circ d_1 \ \hbox{et} \ d_1\circ d_1=0.$$
(Voir \cite{[AAC2]})
\vskip0.2cm

En posant $D_1=d$, $D_2=\mu$, $D_k=0$, si $k\geq3$ et
$$D(\alpha_1\underline{\otimes}\dots\underline{\otimes}
\alpha_n)=\displaystyle\sum_{1\leq r\leq n\atop 0\leq j\leq
n-r}(-1)^{\sum_{i\leq
j}\alpha_i}\alpha_1\underline{\otimes}\dots\underline{\otimes}
\alpha_j\underline\otimes
D_r(\alpha_{j+1}\underline{\otimes}\dots\underline{\otimes}
\alpha_{j+r})\underline{\otimes}
\alpha_{j+r+1}\underline\otimes\dots\underline{\otimes}
\alpha_n.$$

Alors, $D=d_1+\mu_1$ est l'unique cod\'erivation de
$\delta$ de degr\'e $1$ qui prolonge $d$ et $\mu$ \`a
$\mathcal{H}$.

Elle v\'erifie $$D\circ D=0 \ \hbox{et} \ (D\otimes
id+id\otimes D)\circ\delta=\delta\circ D.$$
On obtient que
$\big(\mathcal{H},\delta,D\big)$ est une cog\`ebre de Lie
codiff\'erentielle, donc, c'est une
$C_\infty$-alg\`ebre.\vskip0.15cm

On prolonge, ensuite, le crochet $\ell$ \`a $\mathcal{H}$.

\begin{prop}

\

Sur $\mathcal{H}$, il existe un unique
"crochet" $\ell_2$, de degr\'e $b-a+1$, v\'erifiant:
$$\delta\circ\ell_2=(\ell_2\otimes id)\circ\big(\tau_{23}\circ(\delta\otimes id)
+id\otimes\delta\big)+(id\otimes\ell_2)\circ\big(\delta\otimes
id+\tau_{12}\circ(id\otimes \delta)\big): \ (\ast).$$
Ce crochet est
d\'efini pour $X=\alpha_1\underline{\otimes}\dots\underline{\otimes} \alpha_p$
et $Y=\alpha_{p+1}\underline{\otimes}\dots\underline{\otimes}
\alpha_{p+q}$ par:
$$ \ell_2(X,Y)=\hskip-0,95cm\displaystyle\sum_{\sigma\in Bat(p,q)\atop
k,\sigma^{-1}(k)\leq p<\sigma^{-1}(k+1)}\hskip-0,95cm
\varepsilon_\alpha(\sigma^{-1})(-1)^{(b-a+1)\sum_{s<k}\alpha_{\sigma^{-1}(s)}}
\alpha_{\sigma^{-1}(1)}\underline{\otimes}\dots\underline{\otimes}\ell(\alpha_{\sigma^{-1}(k)},
\alpha_{\sigma^{-1}(k+1)})\underline{\otimes}
\dots\underline{\otimes}\alpha_{\sigma^{-1}(p+q)}.$$

\end{prop}

{\bf Preuve:}

Le prolongement $\ell_2$ d\'efini ci-dessus v\'erifie bien $(\ast)$, en effet:

\begin{align*}& \delta\circ\ell_2(X,Y)=\cr&\delta\Big(\sum_{\sigma\in Bat(p,q)\atop
k,\sigma^{-1}(k)\leq p<\sigma^{-1}(k+1)}\hskip-1,03cm
\varepsilon_\alpha(\sigma^{-1})(-1)^{(b-a+1)\sum_{s<k}\alpha_{\sigma^{-1}(s)}}
\alpha_{\sigma^{-1}(1)}\underline{\otimes}\dots\underline{\otimes}\ell(\alpha_{\sigma^{-1}(k)},
\alpha_{\sigma^{-1}(k+1)})\underline{\otimes}
\dots\underline{\otimes}\alpha_{\sigma^{-1}(p+q)}\Big)\cr&=\sum_{\begin{smallmatrix}
\sigma\in Bat(p,q)\cr 1\leq j<k\leq p+q-1\cr\sigma^{-1}(k)\leq
p<\sigma^{-1}(k+1)\end{smallmatrix}}
\varepsilon_\alpha(\sigma^{-1})(-1)^{(b-a+1)\sum_{s<k}\alpha_{\sigma^{-1}(s)}}~\times\\&\hskip2.1cm\times\Big[\alpha_{\sigma^{-1}(1)}\underline{\otimes}\dots\underline{\otimes}\alpha_{\sigma^{-1}(j)}\bigotimes
\alpha_{\sigma^{-1}(j+1)}\underline{\otimes}\dots\underline{\otimes}\ell(\alpha_{\sigma^{-1}(k)},
\alpha_{\sigma^{-1}(k+1)})\underline{\otimes}\dots\underline{\otimes}\alpha_{\sigma^{-1}(p+q)}\cr&\hskip2.1cm-
\varepsilon\left(\begin{smallmatrix}\alpha_{[\sigma^{-1}(1),\sigma^{-1}(j)]}&\alpha_{\sigma^{-1}(j+1)}
\dots\ell(\alpha_{\sigma^{-1}(k)},
\alpha_{\sigma^{-1}(k+1)})\dots\alpha_{\sigma^{-1}(p+q)}\\
\alpha_{\sigma^{-1}(j+1)} \dots\ell(\alpha_{\sigma^{-1}(k)},
\alpha_{\sigma^{-1}(k+1)})\dots\alpha_{\sigma^{-1}(p+q)}&\alpha_{[\sigma^{-1}(1),\sigma^{-1}(j)]}\end{smallmatrix}\right)\times
\cr&\hskip2.5cm\times\alpha_{\sigma^{-1}(j+1)}\underline{\otimes}\dots\underline{\otimes}\ell(\alpha_{\sigma^{-1}(k)},
\alpha_{\sigma^{-1}(k+1)})\underline{\otimes}\dots\underline{\otimes}\alpha_{\sigma^{-1}(p+q)}
\bigotimes\alpha_{\sigma^{-1}(1)}\underline{\otimes}\dots\underline{\otimes}\alpha_{\sigma^{-1}(j)}\Big]\cr&
+\sum_{\begin{smallmatrix} \sigma\in Bat(p,q)\cr 1\leq
k<j\leq p+q-1\cr\sigma^{-1}(k)\leq
p<\sigma^{-1}(k+1)\end{smallmatrix}}
\varepsilon_\alpha(\sigma^{-1})(-1)^{(b-a+1)\sum_{s<k}\alpha_{\sigma^{-1}(s)}}~
\times\\&\hskip2.5cm\times\Big[\alpha_{\sigma^{-1}(1)}\underline{\otimes}\dots\underline{\otimes}\ell(\alpha_{\sigma^{-1}(k)},
\alpha_{\sigma^{-1}(k+1)})\underline{\otimes}\dots\underline{\otimes}\alpha_{\sigma^{-1}(j)}\bigotimes
\alpha_{\sigma^{-1}(j+1)}\underline{\otimes}\dots\underline{\otimes}\alpha_{\sigma^{-1}(p+q)}\cr&\hskip2.5cm-
\varepsilon\left(\begin{smallmatrix}\alpha_{\sigma^{-1}(1)}\dots\ell(\alpha_{\sigma^{-1}(k)},
\alpha_{\sigma^{-1}(k+1)})\dots\alpha_{\sigma^{-1}(j)}&\alpha_{[\sigma^{-1}(j+1),\sigma^{-1}(p+q)]}\\
\alpha_{[\sigma^{-1}(j+1),\sigma^{-1}(p+q)]}&\alpha_{\sigma^{-1}(1)}\dots\ell(\alpha_{\sigma^{-1}(k)},
\alpha_{\sigma^{-1}(k+1)})\dots\alpha_{\sigma^{-1}(j)}\end{smallmatrix}\right)\times\cr&\hskip2.6cm\times\alpha_{\sigma^{-1}(j+1)}\underline{\otimes}\dots\underline{\otimes}
\alpha_{\sigma^{-1}(p+q)}\bigotimes\alpha_{\sigma^{-1}(1)}\underline{\otimes}\dots\underline{\otimes}
\ell(\alpha_{\sigma^{-1}(k)},
\alpha_{\sigma^{-1}(k+1)})\underline{\otimes}\dots\underline{\otimes}\alpha_{\sigma^{-1}(j)}\Big]
\cr&=\hbox{(I)+(II)+(III)+(IV)}.\end{align*}

Dans $\hbox{(I)}$, on fixe $j$ et $\sigma$ et on pose
$$I=\{\sigma^{-1}(1),\dots,\sigma^{-1}(j)\},
 \ I_1=I\cap\{1,\dots,p\}=\{i_1,\dots,i_r\} \ \hbox{et} \ I_2=I\cap\{p+1,\dots,p+q\}=\{i_{r+1},\dots,i_j\}.$$

Si $I\not\subset\{1,\dots,p\}$ ou $I\not\subset\{p+1,\dots,p+q\}$,
on fixe $k>j$ tel que $\sigma^{-1}(k)\leq p<\sigma^{-1}(k+1)$ et
on fixe
$\big(\sigma^{-1}(j+1),\dots,\sigma^{-1}(p+q)\big)=(i_{j+1},\dots,i_{p+q})$.

\noindent Alors, l'ensemble $$\left\{\sigma\in Bat(p,q)/
\sigma^{-1}(\{1,\dots,j\})=I \ \hbox{et} \
\big(\sigma^{-1}(j+1),\dots,\sigma^{-1}(p+q)\big)=(i_{j+1},\dots,i_{p+q})\right\}$$
est en bijection avec l'ensemble
$$\left\{\rho\in Bat(r,j-r) \
\hbox{d\'efini sur} \ I=I_1\cup I_2\right\}.$$

Et on a,
$\varepsilon_\alpha(\sigma^{-1})=\varepsilon_\alpha(\rho^{-1})$.
Dans ce cas, les termes de $\hbox{(I)}$ sont de la forme $$\pm
bat_{r,j-r}(\alpha_{i_1}\underline{\otimes}\dots\underline{\otimes}\alpha_{i_r},
\alpha_{i_{r+1}}\underline{\otimes}\dots\underline{\otimes}\alpha_{i_j}
)\bigotimes\alpha_{i_{j+1}}\underline{\otimes}\dots\underline{\otimes}\ell(\alpha_{i_k},
\alpha_{i_{k+1}})\underline{\otimes}\dots\underline{\otimes}\alpha_{i_{p+q}}.$$

Ils sont nuls dans $\underline{\bigotimes}^+\mathcal{A}[-a+1]$.
Donc, dans $\hbox{(I)}$, on peut supposer que
$I\subset\{1,\dots,p\}$, dans ce cas la somme correspondante est
no\'ee $\hbox{(I}_1)$, ou que $I\subset\{p+1,\dots,p+q\}$, dans ce
cas la somme correspondante est not\'ee $\hbox{(I}_2)$.\vskip0.15cm

Ce m\^eme raisonnement est vrai pour
$\hbox{(II)}=\hbox{(II}_1)+\hbox{(II}_2)$.\vskip0.15cm

Concernant $\hbox{(III)}$ et $\hbox{(IV)}$, on pose
$I=\{\sigma^{-1}(j+1),\dots,\sigma^{-1}(p+q)\}$ et on applique le
m\^eme raisonnement, pour la m\^eme raison, il ne nous reste que les termes
o\`u $I\subset\{1,\dots,p\}$ dans ce cas les sommes
correspondantes sont not\'ees $\hbox{(III}_1)$ et  $\hbox{(IV}_1)$
et les termes o\`u $I\subset\{p+1,\dots,p+q\}$ dans ce cas les sommes
correspondantes sont not\'ees $\hbox{(III}_2)$ et $\hbox{(IV}_2)$.

Alors,
\begin{align*}& \delta\circ\ell_2(X,Y)=\hbox{(I}_1)+\hbox{(II}_1)+\hbox{(I}_2)+\hbox{(II}_2)+\hbox{(III}_1)
+\hbox{(IV}_1)+\hbox{(III}_2)+\hbox{(IV}_2).\end{align*}

D'autre part,
\begin{align*}&
\left((\ell_2\otimes id)\circ\Big(\tau_{23}\circ(\delta\otimes
id)+id\otimes\delta\Big)+(id\otimes\ell_2)\circ\Big(\delta\otimes
id+\tau_{12}\circ(id\otimes\delta)\Big)\right)(X,Y)=
\cr&(\ell_2\otimes
id)\circ\tau_{23}\Big[\sum_{j=1}^{p-1}\alpha_1\underline{\otimes}\dots\underline{\otimes}
\alpha_j\bigotimes \alpha_{j+1}
\underline{\otimes}\dots\underline{\otimes}
\alpha_p\bigotimes\alpha_{p+1}
\underline{\otimes}\dots\underline{\otimes} \alpha_{p+q}
\cr&\hskip2.5cm-\varepsilon\left(\begin{smallmatrix}\alpha_1\dots
\alpha_{j}&\alpha_{j+1}\dots \alpha_p\\ \alpha_{j+1}\dots
\alpha_p&\alpha_1\dots
\alpha_j\end{smallmatrix}\right)\alpha_{j+1}\underline{\otimes}\dots\underline{\otimes}
\alpha_p\bigotimes
\alpha_1\underline{\otimes}\dots\underline{\otimes}
\alpha_j\bigotimes\alpha_{p+1}
\underline{\otimes}\dots\underline{\otimes}
\alpha_{p+q}\Big]\cr&+(\ell_2\otimes
id)\Big[\sum_{i=p+1}^{p+q-1}\alpha_1\underline{\otimes}\dots\underline{\otimes}
\alpha_p\bigotimes \alpha_{p+1}
\underline{\otimes}\dots\underline{\otimes}
\alpha_i\bigotimes\alpha_{i+1}
\underline{\otimes}\dots\underline{\otimes} \alpha_{p+q}
\cr&\hskip2.3cm-\varepsilon\left(\begin{smallmatrix}\alpha_{p+1}\dots
\alpha_{i}&\alpha_{i+1}\dots \alpha_{p+q}\\ \alpha_{i+1}\dots
\alpha_{p+q}&\alpha_{p+1}\dots
\alpha_i\end{smallmatrix}\right)\alpha_{1}\underline{\otimes}\dots\underline{\otimes}
\alpha_p\bigotimes
\alpha_{i+1}\underline{\otimes}\dots\underline{\otimes}
\alpha_{p+q}\bigotimes\alpha_{p+1}
\underline{\otimes}\dots\underline{\otimes}
\alpha_{i}\Big]\cr&+(id\otimes
\ell_2)\Big[\sum_{j=1}^{p-1}\alpha_1\underline{\otimes}\dots\underline{\otimes}
\alpha_j\bigotimes \alpha_{j+1}
\underline{\otimes}\dots\underline{\otimes}
\alpha_p\bigotimes\alpha_{p+1}
\underline{\otimes}\dots\underline{\otimes} \alpha_{p+q}
\cr&\hskip2.5cm-\varepsilon\left(\begin{smallmatrix}\alpha_1\dots
\alpha_{j}&\alpha_{j+1}\dots \alpha_p\\ \alpha_{j+1}\dots
\alpha_p&\alpha_1\dots
\alpha_j\end{smallmatrix}\right)\alpha_{j+1}\underline{\otimes}\dots\underline{\otimes}
\alpha_p\bigotimes
\alpha_1\underline{\otimes}\dots\underline{\otimes}
\alpha_j\bigotimes\alpha_{p+1}
\underline{\otimes}\dots\underline{\otimes}
\alpha_{p+q}\Big]\cr&+(id\otimes
\ell_2)\circ\tau_{12}\Big[\sum_{i=p+1}^{p+q-1}\alpha_1\underline{\otimes}\dots\underline{\otimes}
\alpha_p\bigotimes \alpha_{p+1}
\underline{\otimes}\dots\underline{\otimes}
\alpha_i\bigotimes\alpha_{i+1}
\underline{\otimes}\dots\underline{\otimes} \alpha_{p+q}
\cr&\hskip2.3cm-\varepsilon\left(\begin{smallmatrix}\alpha_{p+1}\dots
\alpha_{i}&\alpha_{i+1}\dots \alpha_{p+q}\\ \alpha_{i+1}\dots
\alpha_{p+q}&\alpha_{p+1}\dots
\alpha_i\end{smallmatrix}\right)\alpha_{1}\underline{\otimes}\dots\underline{\otimes}
\alpha_p\bigotimes
\alpha_{i+1}\underline{\otimes}\dots\underline{\otimes}
\alpha_{p+q}\bigotimes\alpha_{p+1}
\underline{\otimes}\dots\underline{\otimes} \alpha_{i}\Big]
\cr&=(\ell_2\otimes
id)\Big[\sum_{i=1}^{p-1}(-1)^{\alpha_{[i+1,p]}\alpha_{[p+1,p+q]}}\alpha_1\underline{\otimes}\dots\underline{\otimes}
\alpha_i\bigotimes \alpha_{p+1}
\underline{\otimes}\dots\underline{\otimes} \alpha_{p+q}
\bigotimes\alpha_{i+1} \underline{\otimes}\dots\underline{\otimes}
\alpha_p \cr&\hskip2cm-\sum_{j=1}^{p-1}
(-1)^{\alpha_{[j+1,p+q]}\alpha_{[1,j]}}\alpha_{j+1}\underline{\otimes}\dots\underline{\otimes}
\alpha_p\bigotimes\alpha_{p+1}
\underline{\otimes}\dots\underline{\otimes} \alpha_{p+q}
\bigotimes\alpha_1\underline{\otimes}\dots\underline{\otimes}
\alpha_j \Big]\cr&+(\ell_2\otimes
id)\Big[\sum_{j=p+1}^{p+q-1}\alpha_1\underline{\otimes}\dots\underline{\otimes}
\alpha_p\bigotimes \alpha_{p+1}
\underline{\otimes}\dots\underline{\otimes}
\alpha_j\bigotimes\alpha_{j+1}
\underline{\otimes}\dots\underline{\otimes} \alpha_{p+q}
\cr&\hskip2cm-\sum_{i=p+1}^{p+q-1}(-1)^{\alpha_{[p+1,i]}\alpha_{[i+1,p+q]}}\alpha_{1}\underline{\otimes}\dots\underline{\otimes}
\alpha_p\bigotimes
\alpha_{i+1}\underline{\otimes}\dots\underline{\otimes}
\alpha_{p+q}\bigotimes\alpha_{p+1}
\underline{\otimes}\dots\underline{\otimes}
\alpha_{i}\Big]\cr&+(id\otimes
\ell_2)\Big[\sum_{j=1}^{p-1}\alpha_1\underline{\otimes}\dots\underline{\otimes}
\alpha_j\bigotimes \alpha_{j+1}
\underline{\otimes}\dots\underline{\otimes}
\alpha_p\bigotimes\alpha_{p+1}
\underline{\otimes}\dots\underline{\otimes} \alpha_{p+q}
\cr&\hskip2cm-\sum_{i=1}^{p-1}(-1)^{\alpha_{[i+1,p]}\alpha_{[1,i]}}
\alpha_{i+1}\underline{\otimes}\dots\underline{\otimes}
\alpha_p\bigotimes
\alpha_1\underline{\otimes}\dots\underline{\otimes}
\alpha_i\bigotimes\alpha_{p+1}
\underline{\otimes}\dots\underline{\otimes}
\alpha_{p+q}\Big]\cr&+(id\otimes
\ell_2)\Big[\sum_{i=p+1}^{p+q-1}(-1)^{\alpha_{[p+1,i]}\alpha_{[1,p]}}\alpha_{p+1}
\underline{\otimes}\dots\underline{\otimes} \alpha_i
\bigotimes\alpha_1\underline{\otimes}\dots\underline{\otimes}
\alpha_p \bigotimes\alpha_{i+1}
\underline{\otimes}\dots\underline{\otimes} \alpha_{p+q}
\cr&\hskip2cm-\sum_{j=p+1}^{p+q-1}(-1)^{\alpha_{[j+1,p+q]}\alpha_{[1,j]}}\alpha_{j+1}\underline{\otimes}\dots\underline{\otimes}
\alpha_{p+q}\bigotimes\alpha_{1}\underline{\otimes}\dots\underline{\otimes}
\alpha_p \bigotimes\alpha_{p+1}
\underline{\otimes}\dots\underline{\otimes} \alpha_{j}\Big]
\cr&=(1)+(2)+(3)+(4)+(5)+(6)+(7)+(8).\end{align*}

Dans $\hbox{(I}_1)$, il apparaissait des termes de la forme $$\alpha_{1}
\underline{\otimes}\dots\underline{\otimes}\alpha_{j}\bigotimes
\alpha_{\sigma^{-1}(j+1)}\underline{\otimes}\dots\underline{\otimes}\ell(\alpha_{\sigma^{-1}(k)},
\alpha_{\sigma^{-1}(k+1)})\underline{\otimes}\dots\underline{\otimes}\alpha_{\sigma^{-1}(p+q)}$$
avec le signe $\varepsilon_\alpha(\sigma^{-1})(-1)^{(b-a+1)\sum_{s<k}\alpha_{\sigma^{-1}(s)}}$.\vskip0.25cm

On construit une permutation $\rho$ sur $\{j+1,\dots,p+q\}$ telle
que $\rho^{-1}(s)=\sigma^{-1}(s)$, $\forall j+1\leq s\leq p+q$.
Alors, $\rho\in Bat(p-j,q)$ v\'erifiant $j+1\leq\rho^{-1}(k)\leq
p<\rho^{-1}(k+1)$ et
$\varepsilon_\alpha(\rho^{-1})=\varepsilon_\alpha(\sigma^{-1})$.\vskip0.2cm

De plus, on a
$(-1)^{(b-a+1)\sum_{s<k}\alpha_{\sigma^{-1}(s)}}=(-1)^{(b-a+1)\sum_{j<s<k}
\alpha_{\rho^{-1}(s)}}(-1)^{(b-a+1)\alpha_{[1,j]}}$.

\noindent Alors, le
terme pr\'ec\'edent s'\'ecrit:
\begin{align*}&(-1)^{(b-a+1)\alpha_{[1,j]}}\varepsilon_\alpha(\rho^{-1})(-1)^{(b-a+1)\sum_{j<s<k}\alpha_{\rho^{-1}(s)}}\times\\&\times\alpha_{[1,j]}
\bigotimes\alpha_{\rho^{-1}(j+1)}\underline{\otimes}\dots\underline{\otimes}\ell(\alpha_{\rho^{-1}(k)},
\alpha_{\rho^{-1}(k+1)})\underline{\otimes}\dots\underline{\otimes}\alpha_{\rho^{-1}(p+q)}.\end{align*}
Ce m\^eme terme appara\^it une seule fois dans le second membre et plus pr\'ecis\'ement dans
$(5)$ accompagn\'e du m\^eme signe provenant de
$(-1)^{(b-a+1)\alpha_{[1,j]}}\alpha_{[1,j]}\bigotimes\ell_2\Big(\alpha_{[j+1,p]},\alpha_{[p+1,p+q]}\Big)$.

\noindent On obtient, donc, $\hbox{(I}_1)=(5)$.\vskip0.35cm

\noindent Pour $\hbox{(I}_2)$, on a
$\{\sigma^{-1}(j+1),\dots,\sigma^{-1}(p+q)\}\subset\{p+1,\dots,p+q\}$,
on pose
$$\{\sigma^{-1}(j+1),\dots,\sigma^{-1}(p+q)\}=\{p+1,\dots,i\}$$ avec
$p+q-j=i-p$, alors, lorsque $j$ varie de $1$ \`a $q-1$, on trouve que $i$
varie de $p+1$ \`a $p+q-1$. Et on obtient
$\hbox{(I}_2)=(7)$.\vskip0.35cm

\noindent Pour $\hbox{(III}_1)$, on a
$\{\sigma^{-1}(j+1),\dots,\sigma^{-1}(p+q)\}\subset\{1,\dots,p\}$,
on pose
$$\{\sigma^{-1}(j+1),\dots,\sigma^{-1}(p+q)\}=\{i+1,\dots,p\}$$ avec
$p+q-j=p-i$, alors, lorsque $j$ varie de $q+1$ \`a $p+q-1$, on trouve que
$i$ varie de $1$ \`a $p-1$. Et on obtient $\hbox{(III}_1)=(1)$.\vskip0.15cm

Pour les autres termes, on d\'emontre que $\hbox{(II}_1)=(2)$,
$\hbox{(II}_2)=(4)$, $\hbox{(III}_2)=(3)$, $\hbox{(IV}_1)=(6)$ et
$\hbox{(IV}_2)=(8)$. \vskip0.15cm

L'unicit\'e de $\ell_2$ est une cons\'equence du fait que $(\mathcal{H},\delta)$ est une cog\`ebre de Lie colibre. (Voir \cite{[BGHHW]})
\hfill$\square$

\

\subsection{Alg\`ebre de Lie diff\'erentielle gradu\'ee associ\'ee \`a une $(a,b)$-alg\`ebre diff\'erentielle}

\

On consid\`ere, maintenant, l'espace $\mathcal{H}[a-b-1]$ muni de
la graduation $dg'(X)=dg(X)-a+b+1$ not\'e simplement par $x'$ pour
$X\in\mathcal{H}[a-b-1]$. On pose
$\ell'_2(X,Y)=(-1)^{(a-b-1)dg'(X)}\ell_2(X,Y)$. Alors, le crochet
$\ell'_2$ est de degr\'e $0$ dans $\mathcal{H}[a-b-1]$ et la
diff\'erentielle $D$ reste de degr\'e $1$. Et on a

\begin{prop}

\

L'espace $\mathcal H[a-b-1]$, muni du crochet $\ell'_2$ et de la
diff\'erentielle
 $D$ est une alg\`ebre de Lie diff\'erentielle
gradu\'ee: Pour tout $X$, $Y$ et $Z$ de $\mathcal H[a-b-1]$, on a:

\begin{itemize}
\item[{\bf(i)}] \quad$\ell'_2(X,Y)=-(-1)^{x'y'}\ell'_2(Y,X)$,

\item[{\bf(ii)}] \quad$(-1)^{x'z'}\ell'_2\left(\ell'_2(X,Y),Z\right)+(-1)^{y'x'}\ell'_2\left(\ell'_2(Y,Z),X\right)
+(-1)^{z'y'}\ell'_2\left(\ell'_2(Z,X),Y\right)=0$,

\item[{\bf(iii)}] \quad$D\left(\ell'_2(X,Y)\right)=\ell'_2\left(D(X),Y\right)+(-1)^{x'}\ell'_2\left(X,D(Y)\right)$.
\end{itemize}
\end{prop}

\noindent {\bf Preuve}

\noindent {\bf(i)} Soient
$X=\alpha_1\underline\otimes\dots\underline\otimes\alpha_p$ et
$Y=\alpha_{p+1}\underline\otimes\dots\underline\otimes\alpha_{p+q}$.
On sait que
$$
\begin{aligned}
&\ell'_2(X,Y)=(-1)^{(a-b-1)x'}\ell_2(X,Y)=\\&\hskip-0.25cm(-1)^{(a-b-1)x'}\hskip-1.32cm\sum_{\begin{smallmatrix}\sigma\in Bat(p,q)\\
k;\sigma^{-1}(k)\leq p<\sigma^{-1}(k+1)\end{smallmatrix}}
\hskip-1.2cm
\varepsilon_\alpha(\sigma^{-1})(-1)^{(b-a+1)\sum_{s<k}\alpha_{\sigma^{-1}(s)}}\alpha_{\sigma^{-1}(1)}\underline\otimes\dots\underline\otimes\ell(
\alpha_{\sigma^{-1}(k)},\alpha_{\sigma^{-1}(k+1)})\underline
\otimes\dots\underline\otimes\alpha_{\sigma^{-1}(p+q)}.
\end{aligned}
$$
Fixons un couple $(\sigma,k)$ dans cette somme tel que
$\sigma^{-1}(k)\leq p<\sigma^{-1}(k+1)$. On d\'efinit deux
permutations $\tau$ et $\rho$  de $S_{p+q}$ par:
$$
\tau(j)=\left\{
              \begin{array}{ll}
               j+p, & \hbox{si $1\leq j\leq q$} \\
                j-q, & \hbox{si $q< j\leq q+p$.}
              \end{array}
            \right.
\hbox{et}\qquad\rho=\sigma\circ\tau.
$$

On v\'erifie que $\rho\in Bat(q,p)$ tel que $\rho^{-1}(k+1)\leq
q<\rho^{-1}(k)$. En posant $\beta_j=\alpha_{\tau(j)}, (
1\leq j\leq p+q)$, on aura
$\beta_{\rho^{-1}(j)}=\alpha_{\sigma^{-1}(j)}$ et
$\varepsilon_{\beta}(\rho^{-1})=(-1)^{xy}\varepsilon_{\alpha}(\sigma^{-1})$.\vskip
0.2cm

On construit ensuite une nouvelle permutation $\nu$ de $S_{p+q}$
d\'efinie par:
$$\nu^{-1}(j)=\rho^{-1}(j),\hskip 0.35cm \forall
j\notin\{k,k+1\}, \ \nu^{-1}(k)=\rho^{-1}(k+1)\ \hbox{et} \
\nu^{-1}(k+1)=\rho^{-1}(k).$$

\vskip0.2cm

On v\'erifie que $\nu\in Bat(q,p)$ tel que \begin{align*}&\hskip3cm\nu^{-1}(k)\leq
q<\nu^{-1}(k+1) \ \hbox{et}\\& \varepsilon_\beta(\nu^{-1})=(-1)^{\beta_{\rho^{-1}(k)}\beta_{\rho^{-1}(k+1)}}\varepsilon_\beta(\rho^{-1})=(-1)^{\alpha
_{\sigma^{-1}(k)}\alpha_{\sigma^{-1}(k+1)}}(-1)^{xy}\varepsilon_\alpha(\sigma^{-1}).\end{align*}
\vskip 0.2cm

De plus l'application $(\sigma,k)\mapsto(\nu,k)$ est une bijection
sur les ensembles correspondants.

Alors,
\begin{align*}&\ell'_2(X,Y)=(-1)^{(a-b-1)x'}(-1)^{xy}\times\\&\times\sum_{\nu\in
Bat(q,p)\atop k ; \ \nu^{-1}(k)\leq q<\nu^{-1}(k+1)}\hskip -0.9cm
\varepsilon_{\beta}(\nu^{-1})(-1)^{\beta_{\nu^{-1}(k)}\beta_{\nu^{-1}(k+1)}}\beta_{\nu^{-1}(1)}\underline{\otimes}\dots\underline{\otimes}
\ell(\beta_{\nu^{-1}(k+1)},\beta_{\nu^{-1}(k)})
\underline{\otimes}\dots\underline{\otimes}\beta_{\nu^{-1}(p+q)}\\&=(-1)^{(a-b-1)x'}(-1)^{xy+b-a}\hskip-
0.9cm\sum_{\nu\in Bat(q,p)\atop k ; \ \nu^{-1}(k)\leq
q<\nu^{-1}(k+1)}\hskip -0.9cm
\varepsilon_{\beta}(\nu^{-1})\beta_{\nu^{-1}(1)}\underline{\otimes}\dots\underline{\otimes}
\ell(\beta_{\nu^{-1}(k)},\beta_{\nu^{-1}(k+1)})
\underline{\otimes}\dots\underline{\otimes}\beta_{\nu^{-1}(p+q)}\\&=-(-1)^{x'y'}\ell'_2(Y,X).\end{align*}

\

\noindent {\bf(ii)} Soient
$X=\alpha_1\underline\otimes\dots\underline\otimes\alpha_p$,
$Y=\beta_1\underline\otimes\dots\underline\otimes\beta_q$ et
$Z=\gamma_1\underline\otimes\dots\underline\otimes\gamma_r$.
Posons
$$
\xi_i=\left\{\begin{aligned}\alpha_i&\quad\text{si }~1\leq i\leq p\\
\beta_{i-p}&\quad\text{si }~p+1\leq i\leq p+q\\
\gamma_{i-p-q}&\quad\text{si }~p+q+1\leq i\leq p+q+r.
\end{aligned}
\right.
$$
On revient \`a $\ell_2$ en \'ecrivant l'identit\'e de Jacobi sous la forme:
\begin{align*}&(-1)^{x'z'}\ell'_2\left(\ell'_2(X,Y),Z\right)+(-1)^{y'x'}\ell'_2\left(\ell'_2(Y,Z),X\right)
+(-1)^{z'y'}\ell'_2\left(\ell'_2(Z,X),Y\right)=\\&(-1)^{(a-b-1)(x+y+z+1)}\left\{(-1)^{xz}\ell_2\left(\ell_2(X,Y),Z\right)
+(-1)^{yx}\ell_2\left(\ell_2(Y,Z),X\right)
+(-1)^{zy}\ell_2\left(\ell_2(Z,X),Y\right)\right\}.\end{align*}

En \'ecrivant $(-1)^{xz}\ell_2\left(\ell_2(X,Y),Z\right)$, on
trouve des termes de la forme:
$$
\begin{aligned}
&(I_1):\
(-1)^{xz}\varepsilon_{_\xi}(\rho_{_{I_1}}^{-1})\xi_{i_1}\underline\otimes\dots
\underline\otimes\xi_{i_{m_1}}\underline\otimes\ell(\alpha_i,\beta_j)
\underline\otimes
\xi_{i_{m_1+3}}\underline\otimes\dots\underline\otimes\xi_{i_{n_1}}\underline\otimes
\ell(\beta_k,\gamma_l)\underline\otimes
\xi_{i_{n_1+3}}\underline\otimes
\dots\underline\otimes\xi_{i_{p+q+r}}\\
&(I_2):\
(-1)^{xz}\varepsilon_{_\xi}(\rho_{_{I_2}}^{-1})\xi_{j_1}\underline\otimes\dots
\underline\otimes
\xi_{j_{m_2}}\underline\otimes\ell(\beta_j,\gamma_l)\underline\otimes
\xi_{j_{m_2+3}}\underline\otimes\dots\underline\otimes\xi_{j_{n_2}}\underline\otimes
\ell(\alpha_i,\beta_k)\underline\otimes
\xi_{j_{n_2+3}}\underline\otimes
\dots\underline\otimes \xi_{j_{p+q+r}}\\
&(I_3):\
(-1)^{xz}\varepsilon_{_\xi}(\rho_{_{I_3}}^{-1})\xi_{k_1}\underline\otimes\dots\underline\otimes
\xi_{k_{m_3}}\underline\otimes\ell(\alpha_i,\beta_j)\underline\otimes
\xi_{k_{m_3+3}}\underline\otimes\dots\underline\otimes\xi_{k_{n_3}}\underline\otimes\ell(\alpha_k,\gamma_l)
\underline\otimes \xi_{k_{n_3+3}}\underline\otimes
\dots\underline\otimes \xi_{k_{p+q+r}}\cr &(I_4):\
(-1)^{xz}\varepsilon_{_\xi}(\rho_{_{I_4}}^{-1})\xi_{l_1}\underline\otimes\dots\underline\otimes
\xi_{l_{m_4}}\underline\otimes\ell(\alpha_i,\gamma_l)\underline\otimes
\xi_{l_{m_4+3}}\underline\otimes\dots\underline\otimes\xi_{l_{n_4}}\underline\otimes\ell(\alpha_k,\beta_j)\underline\otimes
\xi_{l_{n_4+3}}\underline\otimes \dots\underline\otimes
\xi_{l_{p+q+r}}\cr &(I_5):\
(-1)^{xz}\varepsilon_{_\xi}(\rho_{_{I_5}}^{-1})\xi_{s_1}\underline\otimes\dots\underline\otimes
\xi_{s_{m_5}}\underline\otimes\ell(\ell(\alpha_i,\beta_j),\gamma_k)\underline\otimes
\xi_{s_{m_5+4}}\underline\otimes \dots\underline\otimes
\xi_{s_{p+q+r}}.
\end{aligned}
$$
En \'ecrivant $(-1)^{yx}\ell_2\left(\ell_2(Y,Z),X\right)$, on
trouve des termes de la forme:
$$
\begin{aligned}
&(II_1):\
(-1)^{yx}\varepsilon_{_\xi}(\rho_{_{II_1}}^{-1})\xi_{i_1}\underline\otimes\dots\underline\otimes
\xi_{i_{m_1}}
\underline\otimes\ell(\beta_j,\alpha_i)\underline\otimes
\xi_{i_{m_1+3}}\underline\otimes\dots\underline\otimes\xi_{i_{n_1}}\underline\otimes
\ell(\beta_k,\gamma_l)
\underline\otimes \xi_{i_{n_1+3}}\underline\otimes \dots\underline\otimes \xi_{i_{p+q+r}}\\
&(II_2):\
(-1)^{yx}\varepsilon_{_\xi}(\rho_{_{II_2}}^{-1})\xi_{j_1}\underline\otimes\dots\underline\otimes
\xi_{j_{m_2}}
\underline\otimes\ell(\beta_j,\gamma_l)\underline\otimes
\xi_{j_{m_2+3}}\underline\otimes\dots\underline\otimes\xi_{j_{n_2}}\underline\otimes\ell(\beta_k,\alpha_i)
\underline\otimes \xi_{j_{n_2+3}}\underline\otimes \dots\underline\otimes \xi_{j_{p+q+r}}\\
&(II_3):\
(-1)^{yx}\varepsilon_{_\xi}(\rho_{_{II_3}}^{-1})\xi_{t_1}\underline\otimes\dots\underline\otimes
\xi_{t_{m_6}}
\underline\otimes\ell(\beta_j,\gamma_k)\underline\otimes
\xi_{t_{m_6+3}}\underline\otimes\dots\underline\otimes\xi_{t_{n_6}}\underline\otimes
\ell(\gamma_l,\alpha_i)
\underline\otimes \xi_{t_{n_6+3}}\underline\otimes \dots\underline\otimes \xi_{t_{p+q+r}}\\
&(II_4):\
(-1)^{yx}\varepsilon_{_\xi}(\rho_{_{II_4}}^{-1})\xi_{r_1}\underline\otimes\dots\underline\otimes
\xi_{r_{m_7}}
\underline\otimes\ell(\gamma_k,\alpha_i)\underline\otimes
\xi_{r_{m_7+3}}\underline\otimes\dots\underline\otimes\xi_{r_{n_7}}\underline\otimes\ell(\beta_j,\gamma_l)
\underline\otimes \xi_{r_{n_7+3}}\underline\otimes \dots\underline\otimes \xi_{r_{p+q+r}}\\
&(II_5):\
(-1)^{yx}\varepsilon_{_\xi}(\rho_{_{II_5}}^{-1})\xi_{s_1}\underline\otimes\dots\underline\otimes
\xi_{s_{m_5}}
\underline\otimes\ell(\ell(\beta_j,\gamma_k),\alpha_i)\underline\otimes
\xi_{s_{m_5+4}}\underline\otimes \dots\underline\otimes
\xi_{s_{p+q+r}}.
\end{aligned}
$$
En \'ecrivant $(-1)^{zy}\ell_2\left(\ell_2(Z,X),Y\right)$, on
trouve des termes de la forme:
$$
\begin{aligned}
&(III_1)\hskip-0.1cm:
(-1)^{zy}\varepsilon_{_\xi}(\rho_{_{III_1}}^{-1})\xi_{k_1}\underline\otimes\dots\underline\otimes
\xi_{k_{m_3}}
\underline\otimes\ell(\alpha_i,\beta_j)\underline\otimes
\xi_{k_{m_3+3}}\underline\otimes\dots\underline\otimes\xi_{k_{n_3}}\underline\otimes
\ell(\gamma_l,\alpha_k)
\underline\otimes \xi_{k_{n_3+3}}\underline\otimes \dots\underline\otimes \xi_{k_{p+q+r}}\\
&(III_2)\hskip-0.1cm:
(-1)^{zy}\varepsilon_{_\xi}(\rho_{_{III_2}}^{-1})\xi_{l_1}\underline\otimes\dots\underline\otimes
\xi_{l_{m_4}}
\underline\otimes\ell(\gamma_l,\alpha_i)\underline\otimes
\xi_{l_{m_4+3}}\underline\otimes\dots\underline\otimes\xi_{l_{n_4}}\underline\otimes\ell(\alpha_k,\beta_j)
\underline\otimes \xi_{l_{n_4+3}}\underline\otimes \dots\underline\otimes \xi_{l_{p+q+r}}\\
&(III_3)\hskip-0.1cm:
(-1)^{zy}\varepsilon_{_\xi}(\rho_{_{III_3}}^{-1})\xi_{r_1}\underline\otimes\dots\underline\otimes
\xi_{r_{m_7}}
\underline\otimes\ell(\gamma_k,\alpha_i)\underline\otimes
\xi_{r_{m_7+3}}\underline\otimes\dots\underline\otimes\xi_{r_{n_7}}\underline\otimes
\ell(\gamma_l,\beta_j)
\underline\otimes \xi_{r_{n_7+3}}\underline\otimes \dots\underline\otimes \xi_{r_{p+q+r}}\\
&(III_4)\hskip-0.1cm:
(-1)^{zy}\varepsilon_{_\xi}(\rho_{_{III_4}}^{-1})\xi_{t_1}\underline\otimes\dots\underline\otimes
\xi_{t_{m_6}}
\underline\otimes\ell(\gamma_k,\beta_j)\underline\otimes
\xi_{t_{m_6+3}}\underline\otimes\dots\underline\otimes\xi_{t_{n_6}}\underline\otimes\ell(\gamma_l,\alpha_i)
\underline\otimes \xi_{t_{n_6+3}}\underline\otimes \dots\underline\otimes \xi_{t_{p+q+r}}\\
&(III_5)\hskip-0.1cm:
(-1)^{zy}\varepsilon_{_\xi}(\rho_{_{III_5}}^{-1})\xi_{s_1}\underline\otimes\dots\underline\otimes
\xi_{s_{m_5}}\underline\otimes\ell(\ell(\gamma_k,\alpha_i),\beta_j)
\underline\otimes \xi_{s_{m_5+4}}\underline\otimes
\dots\underline\otimes \xi_{s_{p+q+r}}.
\end{aligned}
$$
On v\'erifie gr\^ace \`a la $(b-a+1)$-antisym\'etrie de $\ell$ que: $$(I_1)+(II_1)=0, (I_2)+(II_2)=0,(I_3)+(III_1)=0, (I_4)+(III_2)=0, (II_3)+(III_4)=0, (II_4)+(III_3)=0.$$ Et gr\^ace \`a la relation de Jacobi pour $\ell$ que: $(I_5)+(II_5)+(III_5)=0$.\vskip 0.2cm

\

\noindent {\bf(iii)}
 On rappelle que $D=d_1+\mu_1$ o\`u
 $$d_1(\alpha_1\underline\otimes \dots \underline\otimes
\alpha_{n})=\displaystyle\sum_{1\leq i\leq
n}(-1)^{\sum_{s<i}\alpha_s}\alpha_1\underline\otimes \dots
\underline\otimes d(\alpha_i)\underline\otimes
\dots\underline\otimes \alpha_{n} $$
 et
$$
\mu_1(\alpha_1\underline\otimes \dots \underline\otimes
\alpha_{n})=\sum_{1\leq
i<n}(-1)^{\sum_{s<i}\alpha_s}\alpha_1\underline\otimes \dots
\underline\otimes \mu(\alpha_i,\alpha_{i+1})\underline\otimes
\dots\underline\otimes \alpha_{n}.
$$

V\'erifions d'abord que
$\mu_1\left(\ell'_2(X,Y)\right)=\ell'_2\left(\mu_1(X),Y\right)+(-1)^{x'}\ell'_2\left(X,\mu_1(Y)\right)$.\vskip0.15cm

Soient $X=\alpha_1\underline\otimes\dots\underline\otimes\alpha_p$
et
$Y=\alpha_{p+1}\underline\otimes\dots\underline\otimes\alpha_{p+q}$.

\vskip0.2cm

D'une part, on d\'eveloppe $\mu_1(\ell'_2(X,Y))$, on trouve les types de termes suivants:
$$(I): \alpha_{\sigma^{-1}(1)}\underline\otimes\dots\underline\otimes\mu(
\alpha_{\sigma^{-1}(i)},\alpha_{\sigma^{-1}(i+1)})\underline\otimes\dots \underline\otimes\ell(
\alpha_{\sigma^{-1}(k)},\alpha_{\sigma^{-1}(k+1)})\underline\otimes\dots\underline\otimes\alpha_{\sigma^{-1}(p+q)}$$
o\`u $\sigma\in Bat(p,q), i<k-1 \ \hbox{et} \ \sigma^{-1}(k)\leq p<\sigma^{-1}(k+1)$,
$$(II): \alpha_{\sigma^{-1}(1)}\underline\otimes\dots\underline\otimes\ell(
\alpha_{\sigma^{-1}(k)},\alpha_{\sigma^{-1}(k+1)})\underline\otimes\dots \underline\otimes\mu(
\alpha_{\sigma^{-1}(i)},\alpha_{\sigma^{-1}(i+1)})\underline\otimes\dots\underline\otimes\alpha_{\sigma^{-1}(p+q)}$$
o\`u $\sigma\in Bat(p,q), i>k+1 \ \hbox{et} \ \sigma^{-1}(k)\leq p<\sigma^{-1}(k+1)$,
$$\hskip-2cm(III): \alpha_{\sigma^{-1}(1)}\underline\otimes\dots\underline\otimes\mu\big(\alpha_{\sigma^{-1}(k-1)},
\ell(\alpha_{\sigma^{-1}(k)},\alpha_{\sigma^{-1}(k+1)})\big)\underline
\otimes\dots\underline\otimes\alpha_{\sigma^{-1}(p+q)}$$
o\`u $\sigma\in Bat(p,q), 1\leq k\leq p+q-1 \ \hbox{et} \ \sigma^{-1}(k)\leq p<\sigma^{-1}(k+1)$,
$$\hskip-2cm(IV): \alpha_{\sigma^{-1}(1)}\underline\otimes\dots\underline\otimes
\mu\big(\ell(\alpha_{\sigma^{-1}(k)},\alpha_{\sigma^{-1}(k+1)}),\alpha_{\sigma^{-1}(k+2)}\big)\underline
\otimes\dots\underline\otimes\alpha_{\sigma^{-1}(p+q)}\Big)$$
o\`u $\sigma\in Bat(p,q), 1\leq k\leq p+q-1 \ \hbox{et} \ \sigma^{-1}(k)\leq p<\sigma^{-1}(k+1)$.\vskip0.2cm

Parmi les termes de $(I)$, on distingue quatre cas:\vskip0.12cm

$(I_1)$  lorsque
$\{\sigma^{-1}(i),\sigma^{-1}(i+1)\}\subset\{1,\dots,p\}$,\vskip0.15cm

$(I_2)$ lorsque
$\{\sigma^{-1}(i),\sigma^{-1}(i+1)\}\subset\{p+1,\dots,p+q\}$,\vskip0.15cm

$(I_3)$ lorsque
$\sigma^{-1}(i)\leq
p <\sigma^{-1}(i+1)$,\vskip0.15cm

$(I_4)$ lorsque
$\sigma^{-1}(i+1)\leq p<\sigma^{-1}(i)$.\vskip0.15cm

Et de m\^eme pour $(II)$. Il y a des simplifications entre les termes de $(I_3)$ et $(I_4)$  et aussi entre les termes de $(II_3)$ et $(II_4)$.
Il ne reste que $(I_1),(I_2),(II_1) \ \hbox{et} \ (II_2)$.\vskip0.15cm

On s\'epare aussi les termes de type $(III)$ et $(IV)$ en deux cat\'egories. On a donc:\vskip0.15cm

$(III_1)$ lorsque $\sigma^{-1}(k-1)<\sigma^{-1}(k)\leq
p<\sigma^{-1}(k+1)$, \vskip0.15cm

$(III_2)$  lorsque
$\sigma^{-1}(k)\leq p<\sigma^{-1}(k-1)<\sigma^{-1}(k+1)$, \vskip0.15cm

$(IV_1)$
lorsque $\sigma^{-1}(k)<\sigma^{-1}(k+2)\leq
p<\sigma^{-1}(k+1)$, \vskip0.15cm

$(IV_2)$ lorsque
$\sigma^{-1}(k)\leq p<\sigma^{-1}(k+1)<\sigma^{-1}(k+2)$.\vskip0.15cm

Alors,
\begin{align*}&\mu_1(\ell'_2(X,Y))=(I_1)+(I_2)+(II_1)+(II_2)+(III_1)+(III_2)+(IV_1)+(IV_2).\end{align*}

D'autre part, dans le second membre, on a

\begin{align*}\ell'_2(\mu_1(X),Y)&=(-1)^{(a-b-1)(x'+1)}\hskip-1.8cm
\sum_{\begin{smallmatrix}i<k-1\\\sigma\in Bat(p,q)\\1\leq i<p\\
k\not\in\{ i, i+1\};\sigma^{-1}(k)\leq
p<\sigma^{-1}(k+1)\end{smallmatrix}}
\hskip-1.8cm\varepsilon_\alpha(\sigma^{-1})(-1)^{\sum_{s<i}\alpha_{\sigma^{-1}(s)}}(-1)^{(b-a+1)\sum_{s<k}\alpha_{\sigma^{-1}(s)}}
(-1)^{(b-a+1)}\times
\\&\hskip0.7cm\times\alpha_{\sigma^{-1}(1)}\underline\otimes\dots\underline\otimes\mu(
\alpha_{\sigma^{-1}(i)},\alpha_{\sigma^{-1}(i+1)})\underline\otimes\dots \underline\otimes\ell(
\alpha_{\sigma^{-1}(k)},\alpha_{\sigma^{-1}(k+1)})\underline\otimes\dots\underline\otimes\alpha_{\sigma^{-1}(p+q)}.\end{align*}
On trouve les types de termes suivants:
$$(I'): \alpha_{\sigma^{-1}(1)}\underline\otimes\dots\underline\otimes\mu(
\alpha_{\sigma^{-1}(i)},\alpha_{\sigma^{-1}(i+1)})\underline\otimes\dots \underline\otimes\ell(
\alpha_{\sigma^{-1}(k)},\alpha_{\sigma^{-1}(k+1)})\underline\otimes\dots\underline\otimes\alpha_{\sigma^{-1}(p+q)}$$
o\`u $\sigma\in Bat(p,q), i<k-1, \sigma^{-1}(k)\leq p<\sigma^{-1}(k+1) \ \hbox{et} \ 1\leq \sigma^{-1}(i),\sigma^{-1}(i+1)\leq p$,
$$(II'): \alpha_{\sigma^{-1}(1)}\underline\otimes\dots\underline\otimes\ell(
\alpha_{\sigma^{-1}(k)},\alpha_{\sigma^{-1}(k+1)})\underline\otimes\dots \underline\otimes\mu(
\alpha_{\sigma^{-1}(i)},\alpha_{\sigma^{-1}(i+1)})\underline\otimes\dots\underline\otimes\alpha_{\sigma^{-1}(p+q)}$$
o\`u $\sigma\in Bat(p,q), i>k+1, \sigma^{-1}(k)\leq p<\sigma^{-1}(k+1) \ \hbox{et} \  1\leq\sigma^{-1}(i),\sigma^{-1}(i+1)\leq p$
$$ \hskip-2.5cm \hbox{et} \hskip0.5cm (III'): \alpha_{\sigma^{-1}(1)}\underline\otimes\dots\underline\otimes\ell\big(\mu(\alpha_{\sigma^{-1}(k)},\alpha_{\sigma^{-1}(k+1)}),\alpha_{\sigma^{-1}(k+2)}
\big)\underline\otimes\dots\underline\otimes\alpha_{\sigma^{-1}(p+q)}$$
o\`u $\sigma\in Bat(p,q), 1\leq k\leq p+q-2 \ \hbox{et} \   \sigma^{-1}(k)<\sigma^{-1}(k+1)\leq p<\sigma^{-1}(k+2)$.\vskip0.2cm

Or, d'apr\`es l'identit\'e de Leibniz, on a
\begin{align*}\ell\big(\mu(\alpha_{\sigma^{-1}(k)},\alpha_{\sigma^{-1}(k+1)}),\alpha_{\sigma^{-1}(k+2)}&
\big)=(-1)^{(b-a+1)(\alpha_{\sigma^{-1}(k)}+1)}\mu\big(\alpha_{\sigma^{-1}(k)},\ell(\alpha_{\sigma^{-1}(k+1)},
\alpha_{\sigma^{-1}(k+2)})\big)\cr&+(-1)^{b-a+1+\alpha_{\sigma^{-1}(k+1)}\alpha_{\sigma^{-1}(k+2)}}\mu
\big(\ell(\alpha_{\sigma^{-1}(k)},\alpha_{\sigma^{-1}(k+2)}),
\alpha_{\sigma^{-1}(k+1)}\big),\end{align*} on s\'epare, alors, $(III')$ en $(III'_1)+(III'_2)$.

\vskip0.2cm

Donc,
$\ell'_2(\mu_1(X),Y)=(I')+(II')+(III'_1)+(III'_2)$.\vskip0.3cm

De m\^eme, on trouve que $(-1)^{x'}\ell'_2(X,\mu_1(Y))=(I'')+(II'')+(III''_1)+(III''_2)$.\vskip0.3cm

On v\'erifie, enfin, que $(I_1)=(I')$,  $(II_1)=(II')$,
$(III_1)=(III'_1)$,  $(IV_1)=(III'_2)$, $(I_2)=(I'')$,  $(II_2)=(II'')$,
$(IV_2)=(III''_1)$ et $(III_2)=(III''_2)$.\vskip0.3cm

De m\^eme, on v\'erifie que $D_1\left(\ell'_2(X,Y)\right)=\ell'_2\left(D_1(X),Y\right)+(-1)^{x'}\ell'_2\left(X,D_1(Y)\right)$.
\hfill$\square$

\

\subsection{La $L_\infty$-alg\`ebre $S^+({\mathcal H}[a-b])$}

\

Dans le paragraphe pr\'ec\'edent, on a montr\'e que $\Big(\mathcal H[a-b-1], \ell'_2,D\Big)$
est une alg\`ebre de Lie diff\'erentielle gradu\'ee. On
consid\`ere l'espace $\mathcal{H}[a-b]$ muni de la graduation
$$dg''(X)=dg'(X)-1=dg(X)-a+b:=x''\ , \ \hbox{pour tout} \ X\in\mathcal{H}[a-b].$$
On voudrait construire la cog\`ebre cocommutative coassociative
$(S^+(\mathcal{H}[a-b]),\Delta)$, o\`u
$S^+(\mathcal{H}[a-b])=\bigoplus_{n\geq1}S^n(\mathcal{H}[a-b])$ et
$\Delta$ est son coproduit qui est de degr\'e $0$ et d\'efini par: \vskip0.15cm

$\forall
X_1\dots X_n\in S^n(\mathcal{H}[a-b]) $,
\begin{align*}\Delta(X_1\dots X_n)=\sum_{I\cup J=\{1,\dots n\}\atop \#I, \#J>0}
\varepsilon\left(\begin{smallmatrix}x_1''\dots x_n''\\
x_I'' x_J''\end{smallmatrix}\right)X_I\bigotimes X_J.\end{align*}
Le crochet $\ell'_2$ \'etait antisym\'etrique de degr\'e $0$ sur
$\mathcal{H}[a-b-1]$. Comme l'on veut une cod\'erivation de
degr\'e $1$ pour $\Delta$, on pose
$\ell''_2(X,Y)=(-1)^{x''}\ell'_2(X,Y)$ qui est une application
sym\'etrique sur $\mathcal{H}[a-b]$ de degr\'e $1$. On a

\begin{prop}

\

Pour tout $X$, $Y$, $Z\in\mathcal H[a-b]$, on a:
\begin{align*}&{\bf (i)} \
 \ell''_2(X,Y)=(-1)^{x''y''}\ell''_2(Y,X), \cr&{\bf (ii)} \ (-1)^{x''z''}\ell''_2(\ell_2''(X,Y),Z)+(-1)^{y''x''}\ell''_2
 (\ell''_2(Y,Z),X)
 +(-1)^{z''y''}\ell''_2(\ell''_2(Z,X),Y)=0,\cr&{\bf (iii)} \  D(\ell''_2(X,Y))=-\ell''_2(D(X),Y)
+(-1)^{1+x''}\ell''_2(X,D(Y)).\end{align*}

\end{prop}

{\bf Preuve:}

{\bf (i)} On a
\begin{align*}\ell''_2(X,Y)&=(-1)^{x''}\ell'_2(X,Y)=(-1)^{x''}(-1)^{x'y'+1}\ell'_2(Y,X)
\\&=(-1)^{x''+(x''+1)(y''+1)+1}(-1)^{y''}\ell''_2(Y,X)\\&=(-1)^{x''y''}\ell''_2(Y,X).\end{align*}\vskip0.25cm

{\bf (ii)} On a
\begin{align*}&(-1)^{x''z''}\ell''_2(\ell_2''(X,Y),Z)+(-1)^{y''x''}\ell''_2
 (\ell''_2(Y,Z),X)
 +(-1)^{z''y''}\ell''_2(\ell''_2(Z,X),Y)=\\&(-1)^{x'+y'+z'+1}\left((-1)^{x'z'}\ell'_2\left(\ell'_2(X,Y),Z\right)
 +(-1)^{y'x'}\ell'_2(\left(\ell'_2(Y,Z),X\right)
+(-1)^{z'y'}\ell'_2\left(\ell'_2(Z,X),Y\right)\right)\\&=0
.\end{align*}

{\bf (iii)} On a
\begin{align*}D(\ell''_2(X,Y))&=(-1)^{x''}D(\ell'_2(X,Y))=(-1)^{x''}\left(\ell'_2\left(D(X),Y\right)+(-1)^{x'}
\ell'_2\left(X,D(Y)\right)\right)\\&=(-1)^{x''}\left((-1)^{x''+1}\ell''_2\left(D(X),Y\right)+(-1)^{x'+x''}
\ell''_2\left(X,D(Y)\right)\right)\\&=-\ell''_2(D(X),Y)
+(-1)^{1+x''}\ell''_2(X,D(Y)).\end{align*}\hfill$\square$

\

On prolonge $\ell''_2$ \`a $S^+(\mathcal{H}[a-b])$ de fa\c{c}on
unique comme une cod\'erivation $\ell''$ de $\Delta$ de degr\'e
$1$ en posant:
$$\ell''(X_1\dots
X_n)=\displaystyle\sum_{i<j}\varepsilon\left(\begin{smallmatrix}x_1''\dots x_n''\\
x_i'' x_j''x_1''\dots \widehat{ij}\dots x_n''
\end{smallmatrix}\right)\ell''_{2}
(X_i,X_j).X_1\dots\widehat{ij}\dots X_n.$$ En utilisant l'identit\'e de
Jacobi, on peut v\'erifier que $\ell''\circ\ell''=0$.\vskip0.2cm

On prolonge, aussi, la diff\'erentielle $D$ \`a
$S^+(\mathcal{H}[a-b])$ comme l'unique cod\'erivation $m$ de
$\Delta$ toujours de degr\'e $1$ en posant:
$$m(X_1\dots
X_n)=\displaystyle\sum_{i=1}^{n}\varepsilon\left(\begin{smallmatrix}x_1''\dots x_n''\\
x_i'' x_1''\dots \widehat{i}\dots x_n''
\end{smallmatrix}\right)D(X_i).X_1\dots\widehat{i}\dots
X_n.$$ Elle v\'erifie $m\circ m=0$.

Si on pose $Q_1=D$, $Q_2=\ell''_2$, $Q_k=0$, ~si
$ k\geq3$ et
$$Q(X_1\dots X_n)=\displaystyle\sum_{I\cup
J=\{1,\ldots,n\}\atop I\neq\emptyset}\varepsilon\left(\begin{smallmatrix}x_1''\dots x_n''\\
x_I'' x_J''\end{smallmatrix}\right)Q_{\#I} (X_{I}).X_J.$$ Alors,
$Q=m+\ell''$ et v\'erifie $Q^2=0$ et $\left(Q\otimes
id+id\otimes Q\right)\circ\Delta=\Delta\circ Q$.\vskip0.15cm

Donc, le complexe $\left(S^+(\mathcal{H}[a-b]),\Delta,Q\right)$ est
une cog\`ebre cocommuative coassociative et codif- f\'erentielle,
alors, c'est une $L_\infty$ alg\`ebre.\vskip0.4cm

\subsection{La $C_\infty$-alg\`ebre $S^+({\mathcal H}[a-b])$}

\

L'espace $\big(\mathcal{H},\delta,D\big)$ \'etant une cog\`ebre
de Lie codiff\'erentielle, on d\'efinit un cocrochet $\delta''$ de
degr\'e $a-b$ sur $\mathcal{H}[a-b]$ par:
$$\delta''(X)=\sum_{U\underline\otimes V=X\atop U,V\neq\emptyset}(-1)^{(a-b)u''}\left( U\bigotimes
V+(-1)^{u''v''+a-b+1}V\bigotimes U\right),\hskip0.5cm \forall
X\in\mathcal{H}[a-b].$$

On prolonge $\delta''$ \`a $S^+(\mathcal{H}[a-b])$ par:
\begin{align*}\delta''(X_1\dots X_n)=&\hskip-0.4cm\displaystyle\sum_{1\leq s\leq n\atop I\cup
J=\{1,\dots,n\}\setminus\{s\}}\hskip-0.6cm\varepsilon\left(\begin{smallmatrix}x_1''\dots x_n''\\
x_I''x_s'' x_J''\end{smallmatrix}\right)\sum_{U_s\underline\otimes
V_s=X_s\atop U_s,V_s\neq\emptyset
}(-1)^{(a-b)(x_I''+u_s'')}\times\\&\times\left(X_I.
U_s\bigotimes V_s . X_J+(-1)^{u_s''v_s''+a-b+1}X_I. V_s\bigotimes
U_s. X_J\right).\end{align*}
Alors, $\delta''$ est un cocrochet sur $S^+(\mathcal{H}[a-b])$
de degr\'e $a-b$. En notant $\tau''$ la volte dans
$S^+(\mathcal{H}[a-b])$, $\delta''$
v\'erifie:

\begin{prop}

\

\begin{align*}&{\bf (i)} \
 \tau''\circ\delta''=-(-1)^{a-b}\delta'': \hskip0.6cm\hbox{($\delta''$ est $(a-b)$-coantisym\'etrique)}, \cr&{\bf (ii)} \ \Big(id^{\otimes3}+\tau_{12}''\circ\tau_{23}''+\tau_{23}''\circ\tau_{12}''\Big)\circ(\delta''\otimes
id)\circ\delta''=0: \hskip0.6cm\hbox{(l'identit\'e de
coJacobi)},\cr&{\bf (iii)} \
(id\otimes\Delta)\circ\delta''=(\delta''\otimes
id)\circ\Delta+\tau_{12}''\circ(id\otimes\delta'')\circ\Delta:\hskip0.6cm\hbox{(l'identit\'e
de coLeibniz)}.\end{align*}

\end{prop}

{\bf Preuve:}

{\bf (i)} On a
\begin{align*}&\tau''\circ\delta''(X_1\dots X_n)=\tau''\Big[\hskip-0.6cm\sum_{1\leq s\leq n\atop I\cup
J=\{1,\dots,n\}\setminus\{s\}}\hskip-0.6cm\varepsilon\left(\begin{smallmatrix}x_1''\dots x_n''\\
x_J''x_s'' x_I''\end{smallmatrix}\right)\sum_{U_s\underline\otimes
V_s=X_s\atop U_s,V_s\neq\emptyset
}(-1)^{(a-b)(x_J''+u_s'')}\times\\&\hskip5cm\times\left(X_J. U_s\bigotimes V_s .
X_I+(-1)^{u_s''v_s''+a-b+1}X_J. V_s\bigotimes U_s.
X_I\right)\Big]\\&=\hskip-0.6cm\sum_{1\leq s\leq n\atop I\cup
J=\{1,\dots,n\}\setminus\{s\}}\hskip-0.6cm\varepsilon\left(\begin{smallmatrix}x_1''\dots x_n''\\
x_J''x_s'' x_I''\end{smallmatrix}\right)\sum_{U_s\underline\otimes
V_s=X_s\atop
U_s,V_s\neq\emptyset}(-1)^{(a-b)(x_J''+u_s'')}\times\\&\hskip0.6cm\times\left[(-1)^{(x''_I+v''_s)(x''_J+u''_s)}V_s
. X_I\bigotimes X_J.
U_s+(-1)^{(x''_J+v''_s)(x''_I+u''_s)}(-1)^{u_s''v_s''+a-b+1}U_s.
X_I\bigotimes X_J. V_s
\right]\\&=(-1)^{a-b+1}\Big[\hskip-0.6cm\sum_{1\leq s\leq n\atop
I\cup
J=\{1,\dots,n\}\setminus\{s\}}\hskip-0.6cm\varepsilon\left(\begin{smallmatrix}x_1''\dots x_n''\\
x_I''x_s'' x_J''\end{smallmatrix}\right)\sum_{U_s\underline\otimes
V_s=X_s\atop U_s,V_s\neq\emptyset
}(-1)^{(a-b)(x_I''+u_s'')}\times\\&\hskip3.6cm\times\left(X_I.
U_s\bigotimes V_s . X_J+(-1)^{u_s''v_s''+a-b+1}X_I. V_s\bigotimes
U_s. X_J\right)\Big]\\&=(-1)^{a-b+1}\delta''(X_1\dots
X_n).\end{align*}

\

{\bf (ii)} On calcule $(\delta''\otimes
id)\circ\delta''(X_1\dots X_n)$, on trouve des termes produit et produit tensoriel de facteurs $X_I,X_J,X_K,U_s,V_s,U_t,V_t$
pour $I\cup J\cup K\cup\{s,t\}=\{1,\dots,n\}$, $s\neq t$, $X_s=U_s\underline{\otimes}V_s$ et  $X_t=U_t\underline{\otimes}V_t$ et des termes produit et produit tensoriel de facteurs $X_I,X_J,X_K,U_s,V_s,W_s$
pour $I\cup J\cup K\cup\{s\}=\{1,\dots,n\}$ et $X_s=U_s\underline{\otimes}V_s\underline{\otimes}W_s$.

\noindent Pour retrouver les termes du premier type, on part de $X_I\bigotimes X_J\bigotimes X_K$ et on ajoute les $U_s,V_s,U_t,V_t$ en suivant le tableau:

$ U_t\bigotimes V_t. U_s\bigotimes V_s$, \ \ \
$ V_t\bigotimes U_t. U_s\bigotimes V_s$, \ \ \
$ U_t\bigotimes V_t. V_s\bigotimes U_s$, \ \ \
$ V_t\bigotimes U_t. V_s\bigotimes U_s$, \ \ \

$ U_s. U_t\bigotimes V_t\bigotimes V_s$, \ \ \
$ U_s. V_t\bigotimes U_t\bigotimes V_s$, \ \ \
$ V_s. U_t\bigotimes V_t\bigotimes U_s$, \ \ \
$ V_s. V_t\bigotimes U_t\bigotimes U_s$, \ \ \

$ U_s\bigotimes U_t\bigotimes V_s. V_t$, \ \ \
$ U_s\bigotimes V_t\bigotimes V_s. U_t$, \ \ \
$ V_s\bigotimes U_t\bigotimes U_s. V_t$, \ \ \
$ V_s\bigotimes V_t\bigotimes U_s. U_t.$\vskip0.12cm

\noindent Pour retrouver les termes du deuxi\`eme type, on part de $X_I\bigotimes X_J\bigotimes X_K$ et on ajoute les $U_s,V_s,W_s$ en suivant le tableau:
$$U_s\bigotimes V_s\bigotimes W_s, \ \ \  V_s\bigotimes U_s\bigotimes W_s, \ \ \ V_s\bigotimes W_s\bigotimes U_s, \ \ \   W_s\bigotimes V_s\bigotimes U_s.$$

Par exemple, le premier \'el\'ement $ U_t\bigotimes V_t. U_s\bigotimes V_s$ du tableau correspond au terme $$(1): X_I. U_t\bigotimes
V_t. X_J. U_s\bigotimes V_s. X_K.$$
Pour $s<t$, ce terme appara\^it deux fois
 dans $\Big(id^{\otimes3}+\tau_{12}''\circ\tau_{23}''+\tau_{23}''\circ\tau_{12}''\Big)\circ(\delta''\otimes id)\circ
 \delta''(X_1\dots X_n)$:

1) Une fois dans $id^{\otimes3}\circ(\delta''\otimes id)\circ\delta''(X_1\dots X_n)$ avec
 le signe $\varepsilon_1$ obtenu ainsi:

- On part de $X_1\dots X_n$, on le ram\`ene en $X_I.
X_t. X_J. X_s. X_K$ avec le signe
$\varepsilon\left(\begin{smallmatrix} x_1''\dots x_n''\\
x_I''~x_t''~x_J''~x_s''~x_K''\end{smallmatrix}\right)$.

\

- On applique $\delta''$ sur $X_I. X_t. X_J. X_s. X_K$ et
pr\'ecis\'ement lorsqu'on coupe $X_s$, il appara\^it le terme $X_I.
X_t. X_J. U_s\bigotimes V_s. X_K$ une seule fois avec le
signe
$$(-1)^{(a-b)(x_I''+x_t''+x_J''+u_s'')}\varepsilon\left(\begin{smallmatrix} x_1''
\dots x_n''\\
x_I''~x_t''~x_J''~x_s''~x_K''\end{smallmatrix}\right).$$

- Ensuite, on applique $(\delta''\otimes id)$ sur $X_I. X_t.
X_J. U_s\bigotimes V_s. X_K$ et pr\'ecis\'ement on coupe
$X_t$, on obtient une seule fois le terme $X_I. U_t\bigotimes
V_t. X_J. U_s\bigotimes V_s. X_K$ avec le signe
\begin{align*}&(-1)^{(a-b)(x_I''+x_t''+x_J''+u_s'')}\varepsilon\left(\begin{smallmatrix} x_1''
\dots x_n''\\
x_I''~x_t''~x_J''~x_s''~x_K''\end{smallmatrix}\right)(-1)^{(a-b)(x_I''+u_t'')}\\&=
(-1)^{(a-b)(x_J''+v_t''+u_s''+a-b)}
\varepsilon\left(\begin{smallmatrix} x_1''
\dots x_n''\\
x_I''~x_t''~x_J''~x_s''~x_K''\end{smallmatrix}\right)=\varepsilon_1.\end{align*}

2) Le terme $(1)$ appara\^it une autre fois dans
$\tau_{12}''\circ\tau_{23}''\circ(\delta''\otimes
id)\circ\delta''(X_1\dots X_n)$ avec le signe
$\varepsilon'_1$:

\

- On part de $X_1\dots X_n$, on le ram\`ene en $X_J.
X_s. X_K. X_t. X_I$ avec le signe
$\varepsilon\left(\begin{smallmatrix} x_1''\dots x_n''\\
x_J''~x_s''~x_K''~x_t''~x_I''\end{smallmatrix}\right)$.

\

- On applique $\delta''$ sur $X_J. X_s. X_K. X_t. X_I$ et plus
pr\'ecis\'ement, lorsqu'on coupe $X_t$, alors, le terme $X_J. X_s. X_K.
V_t\bigotimes U_t. X_I$ appara\^it une seule fois avec le signe
$$(-1)^{(a-b)(x_J''+x_s''+x_K''+u_t'')}(-1)^{u_t''v_t''+a-b+1}\varepsilon\left(\begin{smallmatrix} x_1''\dots x_n''
\\ x_J''~x_s''~x_K''~x_t''~x_I''\end{smallmatrix}
\right)$$ qui s'\'ecrit encore $V_t. X_J. X_s. X_K\bigotimes
X_I. U_t$ accompagn\'e du signe
$$(-1)^{(a-b)(x_J''+x_s''+x_K''+u_t'')+u_t''v_t''+a-b+1}(-1)^{u_t''x_I''+v_t''(x_K''+x_s''+x_J'')}
\varepsilon\left(\begin{smallmatrix} x_1''\dots x_n''
\\ x_J''~x_s''~x_K''~x_t''~x_I''\end{smallmatrix}\right).$$

- Ensuite, on applique $(\delta''\otimes id)$ sur $V_t. X_J.
X_s. X_K\bigotimes X_I. U_t$ et pr\'ecis\'ement on coupe
$X_s$, on obtient une seule fois le terme $V_t. X_J.
U_s\bigotimes V_s. X_K\bigotimes X_I. U_t$ avec le signe
\begin{align*}&(-1)^{(a-b)(x_J''+x_s''+x_K''+u_t'')+u_t''v_t''+a-b+1+u_t''x_I''+v_t''(x_K''+x_s''+x_J'')}
\varepsilon\left(\begin{smallmatrix} x_1''\dots x_n''
\\ x_J''~x_s''~x_K''~x_t''~x_I''\end{smallmatrix}\right)(-1)^{(a-b)(v_t''+x_J''+u_s'')}\\&=
(-1)^{(a-b)(x_t''+x_K''+v_s'')+u_t''v_t''+a-b+1}(-1)^{u_t''x_I''+v_t''(x_K''+x_s''+x_J'')}
\varepsilon\left(\begin{smallmatrix} x_1''\dots x_n''
\\ x_J''~x_s''~x_K''~x_t''~x_I''\end{smallmatrix}\right).\end{align*}

- Puis, on applique $\tau_{12}''\circ\tau_{23}''$, on obtient le
terme $(1)=X_I. U_t\bigotimes V_t. X_J. U_s\bigotimes V_s.
X_K$ avec le signe
\begin{align*}&(-1)^{(a-b)(x_t''+x_K''+v_s'')+u_t''v_t''+a-b+1}(-1)^{u_t''x_I''+v_t''(x_K''+x_s''+x_J'')}
\varepsilon\left(\begin{smallmatrix} x_1''\dots x_n''
\\ x_J''~x_s''~x_K''~x_t''~x_I''\end{smallmatrix}\right)\times\\&\times
(-1)^{(x_I''+u''_t)(v''_t+x''_J+u_s''+v_s''+x_K'')+x_I''u''_t}\cr&=(-1)^{(a-b)(x_J''+v_t''+u_s''+a-b)+1}
\varepsilon\left(\begin{smallmatrix} x_1''
\dots x_n''\\
x_I''~x_t''~x_J''~x_s''~x_K''\end{smallmatrix}\right)\\&=-\varepsilon_1.\end{align*}

Finalement, ces deux termes se simplifient.
\vskip0.25cm

Un calcul p\'enible montre, de m\^eme, que tous les termes du tableau se simplifient dans $\Big(id^{\otimes3}+\tau_{12}''\circ\tau_{23}''+\tau_{23}''\circ\tau_{12}''\Big)\circ(\delta''\otimes
id)\circ\delta''(X_1\dots X_n)$.

\

{\bf (iii)} D'une part, on a
\begin{align*}&(id\otimes\Delta)\circ\delta''(X_1\dots X_n)=\sum_{1\leq s\leq n\atop I\cup
J=\{1,\dots,n\}\setminus\{s\}}\varepsilon\left(\begin{smallmatrix}x_1''\dots x_n''\\
x_I''x_s'' x_J''\end{smallmatrix}\right)\sum_{U_s\underline\otimes
V_s=X_s\atop U_s,V_s\neq\emptyset
}(-1)^{(a-b)(x_I''+u_s'')}\times\\&\hskip5.1cm\times(id\otimes\Delta)\left[X_I.
U_s\bigotimes V_s . X_J+(-1)^{u_s''v_s''+a-b+1}X_I. V_s\bigotimes
U_s. X_J\right]\\&=\hskip-0.4cm\sum_{1\leq s\leq n\atop I\cup
J=\{1,\dots,n\}\setminus\{s\}}\hskip-0.8cm\varepsilon\left(\begin{smallmatrix}x_1''\dots x_n''\\
x_I''x_s'' x_J''\end{smallmatrix}\right)\sum_{U_s\underline\otimes
V_s=X_s\atop U_s,V_s\neq\emptyset
}(-1)^{(a-b)(x_I''+u_s'')}\times\\&\times\sum_{K\cup L=J}\Big[\varepsilon\left(\begin{smallmatrix}v_s'' x_J''\\
x_K''v_s'' x_L''\end{smallmatrix}\right)X_I.U_s\bigotimes X_K\bigotimes V_s.X_L+\varepsilon\left(\begin{smallmatrix}
v_s'' x_J''\\
v_s'' x_L''x_K''\end{smallmatrix}\right)X_I.U_s\bigotimes
V_s.X_L\bigotimes X_K \\&+(-1)^{u_s''v_s''+a-b+1}\hskip-0.12cm\left(\varepsilon\left(\begin{smallmatrix}u_s'' x_J''\\
x_K''u_s'' x_L''\end{smallmatrix}\right)X_I.V_s\bigotimes
X_K\bigotimes U_s.X_L+\varepsilon\left(\begin{smallmatrix}
u_s'' x_J''\\
u_s'' x_L''x_K''\end{smallmatrix}\right)X_I.V_s\bigotimes
U_s.X_L\bigotimes X_K \right) \Big]\\&=(I)+
(II)+(III)+(IV).\end{align*}

\

D'autre part, on a
\begin{align*}&(id\otimes\delta'')\circ\Delta(X_1\dots
X_n)= \sum_{K\cup J=\{1,\dots n\}\atop K, J\neq\emptyset}
\varepsilon\left(\begin{smallmatrix}x_1''\dots x_n''\\
x_K''
x_J''\end{smallmatrix}\right)(id\otimes\delta'')(X_K\bigotimes
X_J)\\&=\sum_{K\cup J=\{1,\dots n\}\atop K, J\neq\emptyset}
\varepsilon\left(\begin{smallmatrix}x_1''\dots x_n''\\
x_K'' x_J''\end{smallmatrix}\right)(-1)^{(a-b)x''_K}\sum_{ s\in
J\atop I\cup
L=J\setminus\{s\}}\hskip-0.4cm\varepsilon\left(\begin{smallmatrix}x_J''\\
x_I''x_s'' x_L''\end{smallmatrix}\right)\sum_{U_s\underline\otimes
V_s=X_s\atop U_s,V_s\neq\emptyset
}(-1)^{(a-b)(x_I''+u_s''+a-b)}~~\times\\&\hskip2cm\times\left[X_K\bigotimes
X_I. U_s\bigotimes V_s . X_L+(-1)^{u_s''v_s''+a-b+1}X_K\bigotimes
X_I. V_s\bigotimes U_s. X_L\right].\end{align*} Alors, en appliquant $\tau''_{12}$, on obtient

\begin{align*}&\tau''_{12}\circ(id\otimes\delta'')\circ\Delta(X_1\dots
X_n)=\\&\sum_{\begin{smallmatrix}1\leq s\leq n\\ I\cup K\cup L
=\{1,\dots,n\}\setminus\{s\}\\K\neq\emptyset\end{smallmatrix}}
\varepsilon\left(\begin{smallmatrix}x_1''\dots x_n''\\
x_K'' x_J''\end{smallmatrix}\right)(-1)^{(a-b)x''_K}\varepsilon\left(\begin{smallmatrix}x_K''x_J''\\
x_K''x_I''x_s'' x_L''\end{smallmatrix}\right)\sum_{U_s\underline\otimes
V_s=X_s\atop U_s,V_s\neq\emptyset
}(-1)^{(a-b)(x_I''+u_s'')}\times\\&\hskip-0.12cm\times\hskip-0.12cm\left[(-1)^{x''_K(x_I''+u_s'')}
X_I. U_s\bigotimes X_K\bigotimes V_s .
X_L+(-1)^{u_s''v_s''+a-b+1+x''_K(x_I''+v_s'')} X_I.
V_s\bigotimes X_K\bigotimes U_s.
X_L\right]\\&=(1)+(2).\end{align*}

De plus, on a
\begin{align*}&(\delta''\otimes id)\circ\Delta(X_1\dots
X_n)= \sum_{K\cup J=\{1,\dots n\}\atop K, J\neq\emptyset}
\varepsilon\left(\begin{smallmatrix}x_1''\dots x_n''\\
x_J'' x_K''\end{smallmatrix}\right)(\delta''\otimes
id)(X_J\bigotimes X_K)\\&=\sum_{K\cup J=\{1,\dots n\}\atop K,
J\neq\emptyset}
\varepsilon\left(\begin{smallmatrix}x_1''\dots x_n''\\
x_J'' x_K''\end{smallmatrix}\right)\sum_{ s\in J\atop I\cup
L=J\setminus\{s\}}\hskip-0.4cm\varepsilon\left(\begin{smallmatrix}x_J''\\
x_I''x_s'' x_L''\end{smallmatrix}\right)\sum_{U_s\underline\otimes
V_s=X_s\atop U_s,V_s\neq\emptyset
}(-1)^{(a-b)(x_I''+u_s'')}\times\\&\hskip2.1cm\times\left[
X_I. U_s\bigotimes V_s . X_L\bigotimes X_K
+(-1)^{u_s''v_s''+a-b+1} X_I. V_s\bigotimes U_s. X_L\bigotimes X_K
\right]\\&=\sum_{\begin{smallmatrix}1\leq s\leq n\\ I\cup K\cup L
=\{1,\dots,n\}\setminus\{s\}\\K\neq\emptyset\end{smallmatrix}}
\varepsilon\left(\begin{smallmatrix}x_1''\dots x_n''\\
x_J'' x_K''\end{smallmatrix}\right)\varepsilon\left(\begin{smallmatrix}x_J''x_K''\\
x_I''x_s'' x_L''x_K''\end{smallmatrix}\right)\sum_{U_s\underline\otimes
V_s=X_s\atop U_s,V_s\neq\emptyset
}(-1)^{(a-b)(x_I''+u_s'')}\times\\&\times\left[ X_I.
U_s\bigotimes V_s . X_L\bigotimes X_K+(-1)^{u_s''v_s''+a-b+1} X_I.
V_s\bigotimes U_s. X_L\bigotimes
X_K\right]\\&=(3)+(4).\end{align*}

En examinant les mon\^omes de type $(I),(II),(III)$ ou $(IV)$ du premier membre, on constate qu'ils apparaissent chaque fois dans un des types bien d\'etermin\'e du second membre et chacun appara\^it une fois et une seule. Il suffit, alors, de v\'erifier l'\'egalit\'e des signes. Par exemple:

\

Dans $(I)$, le terme $X_I. U_s\bigotimes X_K\bigotimes V_s .
X_L$ appara\^it avec le signe
\begin{align*}&(-1)^{(a-b)(x''_I+u''_s)}\varepsilon\left(\begin{smallmatrix} x''_1
\dots
x''_n\\x''_I~x''_s~x''_J\end{smallmatrix}\right)\varepsilon\left(\begin{smallmatrix}
v''_s  x''_J\\
x''_K~v''_s~x''_L\end{smallmatrix}\right)\\&=(-1)^{(a-b)(x''_I+u''_s)}\varepsilon\left(\begin{smallmatrix}
x''_1 \dots
x''_n\\x''_I~x''_s~x''_J\end{smallmatrix}\right)\varepsilon\left(\begin{smallmatrix}
x''_Iu''_sv''_s  x''_J\\
x''_Iu''_sx''_K~v''_s~x''_L\end{smallmatrix}\right).\end{align*} Or
$$\varepsilon\left(\begin{smallmatrix} x''_1 \dots x''_s\dots
x''_n\\x''_I~x''_s
~x''_J\end{smallmatrix}\right)=\varepsilon\left(\begin{smallmatrix}
x''_1 \dots u''_sv''_s\dots
x''_n\\x''_I~u''_sv''_s~x''_J\end{smallmatrix}\right)(-1)^{(a-b)\sum_{i<s\atop
i\in J }x''_i}(-1)^{(a-b)\sum_{i>s\atop i\in I }x''_i}.$$ Le signe
total de $(I)$ est:
\begin{align*}&(-1)^{(a-b)(x''_I+u''_s)}(-1)^{(a-b)\sum_{i<s\atop i\in K\cup L
}x''_i}(-1)^{(a-b)\sum_{i>s\atop i\in I
}x''_i}\varepsilon\left(\begin{smallmatrix} x''_1 \dots
u''_sv''_s\dots
x''_n\\x''_I~u''_sx''_Kv''_s~x''_L\end{smallmatrix}\right)\\&=(-1)^{(a-b)u''_s}(-1)^{(a-b)\sum_{i<s
}x''_i}\varepsilon\left(\begin{smallmatrix} x''_1 \dots
u''_sv''_s\dots
x''_n\\x''_I~u''_sx''_Kv''_s~x''_L\end{smallmatrix}\right).\end{align*}

\

Dans $(1)$, le terme $X_I. U_s\bigotimes X_K\bigotimes V_s .
X_L$ appara\^it avec le signe
\begin{align*}&(-1)^{(a-b)x''_K}(-1)^{(a-b)(x''_I+u''_s)}\varepsilon\left(\begin{smallmatrix}
x''_1 \dots
x''_n\\x''_K~x''_I~x''_s~x''_L\end{smallmatrix}\right)(-1)^{x''_K(x_I''+u_s'')}\\&\hskip-0.15cm=(-1)^{(a-b)(x''_I+x''_K+u''_s)+x''_K(x_I''+u_s'')}\varepsilon\left(\begin{smallmatrix}
x''_1\dots u''_sv''_s\dots  x''_n\\
x''_Kx''_Iu''_s~v''_s~x''_L\end{smallmatrix}\right)(-1)^{(a-b)\sum_{i<s\atop
i\in L }x''_i}(-1)^{(a-b)\sum_{i>s\atop i\in K\cup I
}x''_i}\\&\hskip-0.15cm=(-1)^{(a-b)u''_s}
(-1)^{x''_K(x_I''+u_s'')}\varepsilon\left(\begin{smallmatrix}
x''_1\dots u''_sv''_s\dots  x''_n\\
x''_Kx''_Iu''_s~v''_s~x''_L\end{smallmatrix}\right)(-1)^{(a-b)\sum_{i<s
}x''_i}\\&\hskip-0.15cm=(-1)^{(a-b)u''_s}\varepsilon\left(\begin{smallmatrix}
x''_1\dots u''_sv''_s\dots  x''_n\\
x''_Iu''_s~x''_Kv''_s~x''_L\end{smallmatrix}\right)(-1)^{(a-b)\sum_{i<s
}x''_i}.\end{align*} Alors, on a $(I)=(1)$.\vskip0.2cm

De m\^eme, on v\'erifie que
$(II)=(3)$, $(III)=(2)$ et
$(IV)=(4)$.\hfill$\square$

\vskip0.42cm

On a prouv\'e que $(S^+(\mathcal{H}[a-b]),\delta'')$ est une cog\`ebre de Lie. Montrons qu'avec $Q=m+\ell''$, elle est codiff\'erentielle.
Revenons, donc, aux op\'erateurs $m$ et $\ell''$ d\'efinis sur
$S^+(\mathcal{H}[a-b])$ par:$$m(X_1\dots
X_n)=\displaystyle\sum_{i=1}^{n}\varepsilon\left(\begin{smallmatrix}x_1''\dots x_n''\\
x_i'' x_1''\dots \widehat{i}\dots x_n''
\end{smallmatrix}\right)m(X_i).X_1\dots\widehat{i}\dots
X_n$$et
$$\ell''(X_1\dots
X_n)=\displaystyle\sum_{i<j}\varepsilon\left(\begin{smallmatrix}x_1''\dots x_n''\\
x_i'' x_j''x_1''\dots \widehat{ij}\dots x_n''
\end{smallmatrix}\right)\ell''_{2}
(X_i,X_j).X_1\dots \widehat{ij}\dots X_n.$$

\begin{prop}

\

 $m$ et $\ell''$ sont des cod\'erivations de $\delta''$ de
degr\'e $1$, ils v\'erifient:

{\bf(i)} $\left(m\otimes id+id\otimes
m\right)\circ\delta''=(-1)^{a-b}\delta''\circ m$.

{\bf(ii)} $\left(\ell''\otimes id+id\otimes
\ell''\right)\circ\delta''=(-1)^{a-b}\delta''\circ \ell''$.
\end{prop}

{\bf Preuve:}

{\bf(i)} D'une part, on a
\begin{align*}&\delta''\circ m(X_1\dots
X_n)=\delta''\Big[\sum_{s=1}^{n}(-1)^{\sum_{i<s}x_i''}X_1\dots
m(X_s)\dots
X_n\Big]\cr&=\hskip-0.8cm\sum_{\begin{smallmatrix}t<s\\U_t\underline\otimes V_t=X_t\\
I\cup J=\{1,\dots,n\}\setminus\{s,t\}
\end{smallmatrix}}\hskip-0.5cm(-1)^{\sum_{i<s}x_i''}\Big[(-1)^{(a-b)(x''_I+u_t'')}\varepsilon
\left(\begin{smallmatrix} x_1''\dots x_t''\dots m(x_s)''\dots x_n''\\
x_I''~x_t''~m(x_s)''~x_J''\end{smallmatrix}\right)\times\\&\hskip3.8cm\times\Big( X_I.
U_t\bigotimes V_t. m(X_s). X_J+
(-1)^{u''_tv_t''+a-b+1}X_I.
 V_t\bigotimes U_t. m(X_s).
X_J\Big)\\&\hskip3.8cm+(-1)^{(a-b)(x''_I+u_t''+x''_s+1)}\varepsilon
\left(\begin{smallmatrix} x_1''\dots x_t''\dots m(x_s)''\dots x_n''\\
x_I''~m(x_s)''~x_t''~x_J''\end{smallmatrix}\right)\times\cr&\hskip3.8cm\times\Big(X_I.
m(X_s). U_t \bigotimes V_t. X_J+(-1)^{u''_tv_t''+a-b+1}X_I.
m(X_s). V_t\bigotimes U_t.  X_J\Big)\Big]\cr&
+\hskip-0.8cm\sum_{\begin{smallmatrix}t>s\\U_t\underline\otimes V_t=X_t\\
I\cup J=\{1,\dots,n\}\setminus\{s,t\}
\end{smallmatrix}}\hskip-0.5cm(-1)^{\sum_{i<s}x_i''}\Big[(-1)^{(a-b)(x''_I+u_t'')}\varepsilon
\left(\begin{smallmatrix} x_1''\dots m(x_s)''\dots x_t''\dots x_n''\\
x_I''~x_t''~m(x_s)''~x_J''\end{smallmatrix}\right)\times\\&\hskip3.9cm\times\Big( X_I.
U_t\bigotimes V_t. m(X_s). X_J+
(-1)^{u''_tv_t''+a-b+1}X_I.
 V_t\bigotimes U_t. m(X_s).
X_J\Big)\\&\hskip3.6cm+(-1)^{(a-b)(x''_I+u_t''+x''_s+1)}\varepsilon
\left(\begin{smallmatrix} x_1''\dots  m(x_s)''\dots x_t''\dots x_n''\\
x_I''~m(x_s)''~x_t''~x_J''\end{smallmatrix}\right)\times\cr&\hskip3.9cm\times\Big(X_I.
m(X_s). U_t \bigotimes V_t. X_J+(-1)^{u''_tv_t''+a-b+1}X_I.
m(X_s). V_t\bigotimes U_t.  X_J\Big)\Big]\cr&
 +\sum_{\begin{smallmatrix}1\leq s\leq n\\U_s\underline\otimes V_s=X_s\\
I\cup J=\{1,\dots,n\}\setminus\{s\}
\end{smallmatrix}}(-1)^{\sum_{i<s}x_i''}\varepsilon
\left(\begin{smallmatrix} x_1''\dots m(x_s)'' \dots x_n''\\
x_I''~m(x_s)''~x_J''\end{smallmatrix}\right)\times\\&\hskip0.2cm\times\Big[(-1)^{(a-b)(x''_I+u_s''+1)}\Big(X_I.
m(U_s)\bigotimes
  V_s. X_J+ (-1)^{(u''_s+1)v_s''+a-b+1}X_I. V_s\bigotimes m(U_s) .
  X_J\Big)\\&\hskip0.7cm+(-1)^{(a-b)(x''_I+u_s'')}(-1)^{u_s}
  \Big(X_I.
U_s\bigotimes
  m(V_s). X_J+ (-1)^{(v''_s+1)u_s''+a-b+1}X_I. m(V_s)\bigotimes U_s . X_J\Big)\Big]\cr&=(1)+(2)+(3)+(4)+(5)+(6)+(7)+
 (8),\end{align*}

o\`u \begin{align*}
&(1)=\sum\pm\Big( X_I.
U_t\bigotimes V_t. m(X_s). X_J+
(-1)^{u''_tv_t''+a-b+1}X_I.
 V_t\bigotimes U_t. m(X_s).
X_J\Big),\\&\hbox{etc ...}\\&
(4)=\sum\pm\Big(X_I.
m(X_s). U_t \bigotimes V_t. X_J+(-1)^{u''_tv_t''+a-b+1}X_I.
m(X_s). V_t\bigotimes U_t.  X_J\Big),\\&
(5)=\sum\pm X_I.
m(U_s)\bigotimes
  V_s. X_J , \ \
  (6)=\sum\pm X_I. V_s\bigotimes m(U_s) .
  X_J ,\\&
  (7)=\sum\pm X_I.
U_s\bigotimes
  m(V_s). X_J , \ \
  (8)=\sum\pm X_I. m(V_s)\bigotimes U_s . X_J.
\end{align*}

D'autre part, on a

\begin{align*}&(m\otimes id)\circ \delta''(X_1\dots
X_n)=\\&\sum_{\begin{smallmatrix}t<s\\U_t\underline\otimes V_t=X_t\\
I\cup J=\{1,\dots,n\}\setminus\{s,t\}
\end{smallmatrix}}\hskip-0.5cm\varepsilon
\left(\begin{smallmatrix} x_1''\dots x_t''\dots x_s''\dots x_n''\\
x_I''~x_s''~x_t''~x_J''\end{smallmatrix}\right)(-1)^{(a-b)(x''_I+x''_s+u_t'')}(-1)^{x_I''}\times\\&\hskip2.2cm\times\Big(
X_I. m(X_s). U_t \bigotimes V_t.
X_J+(-1)^{u''_tv_t''+a-b+1}X_I. m(X_s). V_t\bigotimes U_t.
X_J\Big)\cr&+\sum_{\begin{smallmatrix}t>s\\U_t\underline\otimes V_t=X_t\\
I\cup J=\{1,\dots,n\}\setminus\{s,t\}
\end{smallmatrix}}\hskip-0.5cm\varepsilon
\left(\begin{smallmatrix} x_1''\dots x_s''\dots x_t''\dots x_n''\\
x_I''~x_s''~x_t''~x_J''\end{smallmatrix}\right)(-1)^{(a-b)(x''_I+x''_s+u_t'')}(-1)^{x_I''}\times\\&\hskip2.4cm\times\Big(X_I.
m(X_s). U_t \bigotimes V_t.
X_J+(-1)^{u''_tv_t''+a-b+1}X_I. m(X_s). V_t\bigotimes
U_t. X_J\Big)\cr&
 +\sum_{\begin{smallmatrix}1\leq s\leq n\\U_s\underline\otimes V_s=X_s\\
I\cup J=\{1,\dots,n\}\setminus\{s\}
\end{smallmatrix}}\varepsilon
\left(\begin{smallmatrix} x_1''\dots  x_s''\dots x_n''\\
x_I''~x_s''~x_J''\end{smallmatrix}\right)(-1)^{(a-b)(x''_I+u_s'')}(-1)^{x_I''}X_I.
m(U_s)\bigotimes
  V_s. X_J\cr&+ \sum_{\begin{smallmatrix}1\leq s\leq n\\U_s\underline\otimes V_s=X_s\\
I\cup J=\{1,\dots,n\}\setminus\{s\}
\end{smallmatrix}}\varepsilon
\left(\begin{smallmatrix} x_1''\dots  x_s''\dots x_n''\\
x_I''~x_s''~x_J''\end{smallmatrix}\right)(-1)^{(a-b)(x''_I+u_s'')}(-1)^{x_I''}(-1)^{u''_sv_s''+a-b+1}X_I. m(V_s)\bigotimes U_s.
  X_J\\&=(I_1)+(I_2)+(I_3)+(I_4).\end{align*}
et
\begin{align*}&(id\otimes m)\circ \delta''(X_1\dots
X_n)=\\&\sum_{\begin{smallmatrix}t<s\\U_t\underline\otimes V_t=X_t\\
I\cup J=\{1,\dots,n\}\setminus\{s,t\}
\end{smallmatrix}}\hskip-0.5cm\varepsilon
\left(\begin{smallmatrix} x_1''\dots x_t''\dots x_s''\dots x_n''\\
x_I''~x_t''~x_s''~x_J''\end{smallmatrix}\right)(-1)^{(a-b)(x''_I+u_t'')}(-1)^{x_I''+u''_t+v''_t}\times\\&\hskip2.2cm\times\Big(
X_I. U_t \bigotimes V_t. m(X_s).
X_J+(-1)^{u''_tv_t''+a-b+1}X_I.  V_t\bigotimes U_t. m(X_s).
X_J\Big)\cr&+\sum_{\begin{smallmatrix}t>s\\U_t\underline\otimes V_t=X_t\\
I\cup J=\{1,\dots,n\}\setminus\{s,t\}
\end{smallmatrix}}\hskip-0.5cm\varepsilon
\left(\begin{smallmatrix} x_1''\dots x_s''\dots x_t''\dots x_n''\\
x_I''~x_t''~x_s''~x_J''\end{smallmatrix}\right)(-1)^{(a-b)(x''_I+u_t'')}(-1)^{x_I''+u''_t+v''_t}\times\\&\hskip2.6cm\times\Big(X_I.
 U_t \bigotimes V_t. m(X_s).
X_J+(-1)^{u''_tv_t''+a-b+1}X_I.  V_t\bigotimes U_t.
m(X_s). X_J\Big)\cr&
 +\sum_{\begin{smallmatrix}1\leq s\leq n\\U_s\underline\otimes V_s=X_s\\
I\cup J=\{1,\dots,n\}\setminus\{s\}
\end{smallmatrix}}\varepsilon
\left(\begin{smallmatrix} x_1''\dots  x_s''\dots x_n''\\
x_I''~x_s''~x_J''\end{smallmatrix}\right)(-1)^{(a-b)(x''_I+u_s'')}(-1)^{x_I''+u_s''}X_I.
U_s\bigotimes
 m( V_s). X_J\cr&+\sum_{\begin{smallmatrix}1\leq s\leq n\\U_s\underline\otimes V_s=X_s\\
I\cup J=\{1,\dots,n\}\setminus\{s\}
\end{smallmatrix}}\varepsilon
\left(\begin{smallmatrix} x_1''\dots  x_s''\dots x_n''\\
x_I''~x_s''~x_J''\end{smallmatrix}\right)(-1)^{(a-b)(x''_I+u_s'')}(-1)^{x_I''+v_s''} (-1)^{u''_sv_s''+a-b+1}X_I. V_s\bigotimes m(U_s).
  X_J\\&=(II_1)+(II_2)+(II_3)+(II_4).\end{align*}\vskip0.2cm

On v\'erifie que les termes $(1),(2),(3),(4),(5),(6),(7)$ et
 $(8)$ de $\delta''\circ m(X_1\dots
X_n)$ apparaissent aussi dans $\left(m\otimes id+id\otimes
m\right)\circ\delta''(X_1\dots
X_n)$ \`a un signe $(-1)^{a-b}$ pr\`es. \vskip 0.15cm

\noindent Par exemple, le terme $X_I. U_t\bigotimes V_t. m(X_s). X_J+
(-1)^{u''_tv_t''+a-b+1}X_I.
 V_t\bigotimes U_t. m(X_s).
X_J$ de $(1)$ appara\^it dans $\delta''\circ m(X_1. \dots.
X_n)$ avec le signe
$$\varepsilon_1=(-1)^{\sum_{i<s}x_i''}(-1)^{(a-b)(x''_I+u_t'')}\varepsilon
\left(\begin{smallmatrix} x_1''\dots x_t''\dots m(x_s)''\dots x_n''\\
x_I''~x_t''~m(x_s)''~x_J''\end{smallmatrix}\right).$$ Or
$\varepsilon
\left(\begin{smallmatrix} x_1''\dots x_t''\dots m(x_s)''\dots x_n''\\
x_I''~x_t''~m(x_s)''~x_J''\end{smallmatrix}\right)=\varepsilon
\left(\begin{smallmatrix} x_1''\dots x_t''\dots x_s''\dots x_n''\\
x_I''~x_t''~x_s''~x_J''\end{smallmatrix}\right)(-1)^{\sum_{i<s\atop
i\in J }x_i''}(-1)^{\sum_{i>s\atop i\in I}x_i''}$, \vskip0.15cm

\noindent alors,
$$\varepsilon_1=(-1)^{(a-b)(x''_I+u_t'')}\varepsilon
\left(\begin{smallmatrix} x_1''\dots x_t''\dots x_s''\dots x_n''\\
x_I''~x_t''~x_s''~x_J''\end{smallmatrix}\right)(-1)^{x''_I+x''_t}.$$

Le m\^eme terme appara\^it dans $(id\otimes m)\circ\delta''(X_1
\dots X_n)$, et pr\'ecis\'ement dans $(II_1)$, avec le signe
\begin{align*}\varepsilon
\left(\begin{smallmatrix} x_1''\dots x_t''\dots x_s''\dots x_n''\\
x_I''~x_t''~x_s''~x_J''\end{smallmatrix}\right)(-1)^{(a-b)(x''_I+u_t'')}(-1)^{x_I''+u''_t+v''_t}=
(-1)^{a-b}\varepsilon_1.\end{align*}
Alors, on trouve que $(1)=(-1)^{(a-b)}(II_1)$.

\

De m\^eme, on d\'emontre que $(2)=(-1)^{(a-b)}(I_1)$,
$(3)=(-1)^{(a-b)}(II_2)$,
$(4)=(-1)^{(a-b)}(I_2)$,
$(5)=(-1)^{(a-b)}(I_3)$,
$(6)=(-1)^{(a-b)}(II_4)$,
$(7)=(-1)^{(a-b)}(II_3)$ et
$(8)=(-1)^{(a-b)}(I_4)$.

\

{\bf(ii)} On a
\begin{align*}&\delta''\circ\ell''(X_1\dots X_n)=\delta''\Big[\sum_{i<j\atop K=\{1,\dots,n\}\setminus\{i,j\}}\varepsilon\left(\begin
{smallmatrix} x_1''\dots x_n''\\
x_i''~x_j''x_J''\end{smallmatrix}\right)\ell''_2(X_i,X_j).
X_K\Big].
\end{align*}
Dans $\delta''\circ\ell''(X_1\dots X_n)$ appara\^it
des termes de la forme
\begin{align*}&(1):\ell''_2(X_i,X_j). X_{K}.U_s\bigotimes V_s. X_{L},\cr&
(2):\ell''_2(X_i,X_j). X_{K}.V_s\bigotimes U_s. X_{L},\cr&
(3):  X_{K}.U_s\bigotimes V_s.\ell''_2(X_i,X_j). X_{L},\cr&
(4): X_{K}.V_s\bigotimes U_s. \ell''_2(X_i,X_j). X_{L},\end{align*}
pour $K\cup L\cup\{i,j,s\}=\{1,\dots,n\}$, $i<j$ et $X_s=U_s\underline{\otimes}V_s$. \vskip0.15cm

\noindent Et d'autres termes dans le cas
o\`u on coupe $\ell''_2(X_i,X_j)$ par $\delta''$ de la forme
\begin{align*}&(5):X_K.U_{ij}\bigotimes V_{ij}. X_{L},\cr&(6):X_K.V_{ij}\bigotimes U_{ij}. X_{L},\end{align*}
pour $K\cup L\cup\{i,j\}=\{1,\dots,n\}$, $i<j$ et $\ell''_2(X_i,X_j)=U_{ij}\underline{\otimes}V_{ij}$.

\

De plus, en \'ecrivant le premier membre de {\bf (ii)}, on a d'une part,
\begin{align*}&(\ell''\otimes id)\circ\delta''(X_1
\dots X_n)=\\&(\ell''\otimes id)\Big[\sum_{i<j}\hskip-0.12cm\sum_{s\atop
K \cup L=\{1,\dots,n\}\setminus\{i,j,s\}}\hskip-0.12cm\sum_{U_s\underline\otimes
V_s=X_s\atop
U_s,V_s\neq\emptyset}\hskip-0.16cm\varepsilon\left(\begin{smallmatrix}
x_1''\dots x_n''\\
x_i''~x_j''~x_{K}''~x_s''~x_{L}''\end{smallmatrix}
\right)(-1)^{(a-b)(x_i''+x_j''+x_{K}''+u_s'')}\times\\&\hskip3.5cm\times\big(
X_i. X_j. X_{K}. U_s\bigotimes V_s.
X_{L}+(-1)^{u_s''v_s''+a-b+1} X_i. X_j. X_{K}.
V_s\bigotimes U_s. X_{L}\big)\Big]\cr&=\sum_{i<j}\sum_{s\atop
K \cup L=\{1,\dots,n\}\setminus\{i,j,s\}}\sum_{U_s\underline\otimes
V_s=X_s\atop
U_s,V_s\neq\emptyset}\varepsilon\left(\begin{smallmatrix}
x_1''\dots x_n''\\
x_i''~x_j''~x_{K}''~x_s''~x_{L}''\end{smallmatrix}
\right)(-1)^{(a-b)(x_i''+x_j''+x_{K}''+u_s'')}\times\\&\hskip2.4cm\times\big(
\ell''_2(X_i, X_j). X_{K}. U_s\bigotimes V_s.
X_{L}+(-1)^{u_s''v_s''+a-b+1} \ell''_2(X_i, X_j). X_{K}.
V_s\bigotimes U_s. X_{L}\big)\\&+\sum_{i<j\atop
K \cup L=\{1,\dots,n\}\setminus\{i,j\}}\sum_{U_j\underline\otimes
V_j=X_j\atop
U_j,V_j\neq\emptyset}\varepsilon\left(\begin{smallmatrix}
x_1''\dots x_n''\\
x_{K}''~x_i''~x_j''~x_{L}''\end{smallmatrix}
\right)(-1)^{(a-b)(x_{K}''+x_i''+u_j'')+x_K''}\times\\&\hskip3.4cm\times\big(
X_{K}.\ell''_2(X_i, U_j)\bigotimes V_j.
X_{L}+(-1)^{u_j''v_j''+a-b+1} X_{K}.\ell''_2(X_i, V_j)\bigotimes U_j. X_{L}\big)\\&+\sum_{i<j\atop
K \cup L=\{1,\dots,n\}\setminus\{i,j\}}\sum_{U_i\underline\otimes
V_i=X_i\atop
U_i,V_i\neq\emptyset}\varepsilon\left(\begin{smallmatrix}
x_1''\dots x_n''\\
x_{K}''~x_j''~x_i''~x_{L}''\end{smallmatrix}
\right)(-1)^{(a-b)(x_{K}''+x_j''+u_i'')+x_K''}\times\\&\hskip3.4cm\times\big(
X_{K}.\ell''_2(X_j, U_i)\bigotimes V_i.
X_{L}+(-1)^{u_i''v_i''+a-b+1} X_{K}.\ell''_2(X_j, V_i)\bigotimes U_i. X_{L}\big) \\&=(I_1)+(I_2)+(I_3)+(I_4)+(I_5)+(I_6).\end{align*}
Et d'autre part, on a

\begin{align*}&(id\otimes \ell'')\circ\delta''(X_1
\dots X_n)=\\&(id\otimes \ell'')\Big[\sum_{i<j}\sum_{s\atop
K \cup L=\{1,\dots,n\}\setminus\{i,j,s\}}\sum_{U_s\underline\otimes
V_s=X_s\atop
U_s,V_s\neq\emptyset}\varepsilon\left(\begin{smallmatrix} x_1''\dots x_n''\\
x_K''~x_s''~x_{i}''~x_j''~x_{L}''\end{smallmatrix}
\right)(-1)^{(a-b)(x_K''+u_s'')}\times\\&\hskip3.4cm\times\big(
X_{K}.U_s\bigotimes V_s.  X_i. X_j.
X_{L}+(-1)^{u_s''v_s''+a-b+1}X_{K}.V_s\bigotimes U_s.
X_i. X_j. X_{L}\big)\Big]\\&=\sum_{i<j}\hskip-0.12cm\sum_{s\atop
K \cup L=\{1,\dots,n\}\setminus\{i,j,s\}}\hskip-0.12cm\sum_{U_s\underline\otimes
V_s=X_s\atop
U_s,V_s\neq\emptyset}\hskip-0.13cm\varepsilon\left(\begin{smallmatrix} x_1''\dots x_n''\\
x_K''~x_s''~x_{i}''~x_j''~x_{L}''\end{smallmatrix}
\right)(-1)^{(a-b)(x_K''+u_s'')}(-1)^{x_{K}''+u_s''+v_s''}\times\\&\hskip2.4cm\times\big(
X_{K}. U_s\bigotimes V_s. \ell''_2(X_i, X_j).
X_{L}+(-1)^{u_s''v_s''+a-b+1}X_{K}.V_s\bigotimes U_s.
\ell''_2(X_i, X_j).
X_{L}\big)\\&+\sum_{i<j\atop
K \cup L=\{1,\dots,n\}\setminus\{i,j\}}\sum_{U_i\underline\otimes
V_i=X_i\atop
U_i,V_i\neq\emptyset}\varepsilon\left(\begin{smallmatrix}
x_1''\dots x_n''\\
x_{K}''~x_i''~x_j''~x_{L}''\end{smallmatrix}
\right)(-1)^{(a-b)(x_{K}''+u_i'')+x_K''}\times\\&\hskip2.1cm\times\big((-1)^{u_i''}
X_{K}.U_i\bigotimes \ell''_2( V_i,X_i).
X_{L}+(-1)^{u_i''v_i''+a-b+1+v_i''} X_{K}.V_i\bigotimes\ell''_2(U_i, X_j). X_{L}\big)\\&+\sum_{i<j\atop
K \cup L=\{1,\dots,n\}\setminus\{i,j\}}\sum_{U_j\underline\otimes
V_j=X_j\atop
U_j,V_j\neq\emptyset} \varepsilon\left(\begin{smallmatrix}
x_1''\dots x_n''\\
x_{K}''~x_j''~x_i''~x_{L}''\end{smallmatrix}
\right)(-1)^{(a-b)(x_{K}''+u_j'')+x_K''}\times\\&\hskip1.4cm\times\big(
(-1)^{u_j''}X_{K}.U_j\bigotimes\ell''_2(V_j, X_i)\bigotimes
X_{L}+(-1)^{u_j''v_j''+a-b+1+v_j''} X_{K}.V_j\bigotimes\ell''_2(U_j, X_i). X_{L}\big) \\&=(II_1)+(II_2)+(II_3)+(II_4)+(II_5)+(II_6).\end{align*}

Examinons les termes de $\delta''\circ\ell''(X_1\dots X_n)$ et cherchons les termes correspondants dans le premier membre.\vskip 0.2cm

- Dans $\delta''\circ\ell''(X_1\dots X_n)$, le terme
$(1):\ell_2''(X_i,X_j). X_{K}. U_s\bigotimes V_s.
X_{L}$ appara\^it avec le signe
$$\varepsilon_1=(-1)^{(a-b)(x_i''+x_j''+1+x''_{K}+u_s'')}\varepsilon\left(\begin{smallmatrix}
 x_1''\dots x_n''\\ x_i''~x_j''x_K''x_s''x_{L}''\end{smallmatrix}\right).$$

Il correspond au terme $(I_1)$ de $(\ell''\otimes
id)\circ\delta''(X_1. \dots. X_n)$ qui appara\^it avec le
signe $$\varepsilon\left(\begin{smallmatrix}
x_1''\dots x_n''\\
x_i''~x_j''~x_{K}''~x_s''~x_{L}''\end{smallmatrix}
\right)(-1)^{(a-b)(x_i''+x_j''+x_{K}''+u_s'')}=(-1)^{a-b}\varepsilon_1.$$

- Dans $\delta''\circ\ell''(X_1\dots X_n)$, le terme
$(2):\ell_2''(X_i,X_j). X_{K}. V_s\bigotimes U_s.
X_{L}$ appara\^it avec le signe
$$\varepsilon_2=(-1)^{(a-b)(x_i''+x_j''+1+x''_{K}+u_s'')}\varepsilon\left(\begin{smallmatrix}
 x_1''\dots x_n''\\ x_i''~x_j''x_K''x_s''x_{L}''\end{smallmatrix}\right)(-1)^{u_s''v_s''+a-b+1}.$$

Il correspond au terme $(I_2)$ de $(\ell''\otimes
id)\circ\delta''(X_1 \dots X_n)$ qui appara\^it avec le
signe $$\varepsilon\left(\begin{smallmatrix}
x_1''\dots x_n''\\
x_i''~x_j''~x_{K}''~x_s''~x_{L}''\end{smallmatrix}
\right)(-1)^{(a-b)(x_i''+x_j''+x_{K}''+u_s'')}(-1)^{u_s''v_s''+a-b+1}=(-1)^{a-b}\varepsilon_2.$$

- Dans $\delta''\circ\ell''(X_1\dots X_n)$, le terme
$(3): X_{K}.U_s\bigotimes V_s. \ell_2''(X_i,X_j).
X_{L}$ appara\^it avec le signe
$$\varepsilon_3=(-1)^{(a-b)(x_K''+u_s'')}(-1)^{x_s''+x_K''}\varepsilon\left(\begin{smallmatrix}
 x_1''\dots x_n''\\ x_K''~x_s''~x_{i}''~x_j''x_{L}''\end{smallmatrix}\right).$$

Il correspond au terme $(II_1)$ de $(id\otimes
\ell'')\circ\delta''(X_1 \dots X_n)$ qui appara\^it avec le
signe $$\varepsilon\left(\begin{smallmatrix} x_1''\dots x_n''\\
x_K''~x_s''~x_{i}''~x_j''~x_{L}''\end{smallmatrix}
\right)(-1)^{(a-b)(x_K''+u_s'')}(-1)^{x_{K}''+u_s''+v_s''}=(-1)^{a-b}\varepsilon_3.$$

- Dans $\delta''\circ\ell''(X_1\dots X_n)$, le terme
$(4):X_{K}.V_s\bigotimes U_s. \ell_2''(X_i,X_j).
X_{L}$ appara\^it avec le signe
$$\varepsilon_4=(-1)^{(a-b)(x_K''+u_s'')}(-1)^{x_s''+x_K''}\varepsilon\left(\begin{smallmatrix}
 x_1''\dots x_n''\\ x_K''~x_s''~x_{i}''~x_j''x_{L}''\end{smallmatrix}\right)(-1)^{u_s''v_s''+a-b+1}.$$

Il correspond au terme $(II_2)$ de $(id\otimes
\ell'')\circ\delta''(X_1 \dots X_n)$ qui appara\^it avec le
signe $$\varepsilon\left(\begin{smallmatrix} x_1''\dots x_n''\\
x_K''~x_s''~x_{i}''~x_j''~x_{L}''\end{smallmatrix}
\right)(-1)^{(a-b)(x_K''+u_s'')}(-1)^{x_{K}''+u_s''+v_s''}(-1)^{u_s''v_s''+a-b+1}=(-1)^{a-b}\varepsilon_4.$$
\vskip0.25cm

Le cas o\`u on coupe $\ell_2''(X_i,X_j)$ par $\delta''$ correspond aux termes $(5)$ et $(6)$ de $\delta''\circ\ell''(X_1\dots X_n)$ , $(I_3),(I_4),(I_5),(I_6)$ de $(\ell''\otimes id)\circ\delta''(X_1
\dots X_n)$ et $(II_3),(II_4),(II_5),(II_6)$ de $(id\otimes\ell'')\circ\delta''(X_1
\dots X_n)$.

\noindent On pourra voir ce cas comme si on a uniquement
deux paquets $X$ et $Y$. Il suffit de montrer que:
$$(-1)^{a-b}\delta''\circ \ell''_2(X,Y)=\left(\ell''_2\otimes id+id\otimes
\ell''_2\right)\circ\delta''(X,Y): \ (\ast\ast).$$

Mais, le crochet et le cocrochet $\ell''_2$ et $\delta''$ \'etaient $\ell_2$ et $\delta$ sur $\mathcal{H}$. Ils v\'erifient une relation de compatibilit\'e donn\'ee par:
$$\delta\circ\ell_2(X,Y)=\Big((\ell_2\otimes id)\circ\big(\tau_{23}\circ(\delta\otimes id)
+id\otimes\delta\big)+(id\otimes\ell_2)\circ\big(\delta\otimes
id+\tau_{12}\circ(id\otimes \delta)\big)\Big)(X,Y).$$

En d\'ecalant $\mathcal{H}$ par $a-b$ et en \'ecrivant $\ell_2''$ et $\delta''$ sur $S^2(\mathcal{H}[a-b])$, la sym\'etrisation de la relation pr\'ec\'edente donne $(\ast\ast)$.
Ceci ach\`eve la d\'emonstration.
\hfill$\square$

\vskip0.42cm

La proposition pr\'ec\'edente nous montre que $Q=\ell''+m$ est une cod\'erivation de
$S^+(\mathcal{H}[a-b])$ pour $\delta''$ de degr\'e $1$. Alors, le
complexe $\left(S^+(\mathcal{H}[a-b]),\delta'',Q\right)$
est une cog\`ebre de Lie codiff\'erentielle gradu\'ee, donc, c'est aussi
une $C_\infty$ alg\`ebre.\vskip0.15cm

Enfin, le cocrochet $\delta''$ et le coproduit $\Delta$
v\'erifient l'identit\'e de
coLeibniz:$$(id\otimes\Delta)\circ\delta''=(\delta''\otimes
id)\circ\Delta+\tau_{12}''\circ(id\otimes\delta'')\circ\Delta.
$$ Alors,
$\left(S^+\Big(\mathcal{H}[a-b]\Big),\Delta,\delta'',Q\right)$
est une bicog\`ebre codiff\'erentielle gradu\'ee.

\begin{defn}

\

Soit $\mathcal{A}$ est une $(a,b)$-alg\`ebre diff\'erentielle. Alors, la bicog\`ebre
colibre et codiff\'erentielle $$\left(S^+
\left(\Big(\underline{\displaystyle\bigotimes}^+\mathcal{A}[-a+1]\Big)[a-b]\right),\Delta,\delta'',Q=\ell''+m\right)$$ est la $(a,b)$-alg\`ebre \`a homotopie pr\`es enveloppante de $\mathcal{A}$.

\end{defn}

\vskip0.2cm

\begin{rema}

\

- Dans le cas o\`u $a=0$, $b=-1$ et $\mathcal{A}$ est une alg\`ebre de Gerstenhaber diff\'erentielle, on retrouve l'alg\`ebre de
Gerstenhaber \`a homotopie pr\`es enveloppante de $\mathcal{A}$: $$\left(S^+
\Big(\big(\underline{\displaystyle\bigotimes}^+\mathcal{A}[1]\big)[1]\Big),\Delta,\delta'',Q=\ell''+m\right).$$

- Dans le cas o\`u $a=b=0$ et $\mathcal{A}$ est une alg\`ebre de Poisson diff\'erentielle grdu\'ee, on retrouve le complexe de l'alg\`ebre
de Poisson \`a homotopie pr\`es enveloppante de $\mathcal{A}$: $$\left(S^+
\Big(\underline{\displaystyle\bigotimes}^+\mathcal{A}[1]\Big),\Delta,\delta'',Q=\ell''+m\right).$$

On conclut que notre construction g\'en\'eralise celle des
alg\`ebres de Gerstenhaber et de Poisson \`a homotopie pr\`es.

\end{rema}

\vskip0.4cm


\

\section{Exemples}\label{sec5}\

\

\subsection{Exemple 1}

\

On consid\`ere l'espace vectoriel $T_{poly}(\mathbb{R}^d)$ des multichamps de vecteurs totalement antisym\'etrique sur $\mathbb{R}^d$. Cet espace muni du crochet de Schouten $[~,~]_S$, du produit ext\'erieur $\wedge$ et de la graduation $\hbox{degr\'e}(\alpha)=k$, si $\alpha$ est un $k$-tenseur, est bien s\^ur une alg\`ebre de Gerstenhaber.

On appelle $E$ le sous espace  vectoriel de $T_{poly}(\mathbb{R}^d)$ engendr\'e par les tenseurs
$$\alpha=\sum_{i_1,\dots,i_k}\alpha^{i_1\dots i_k}(x) \partial x_{i_1}\wedge\dots\wedge\partial x_{i_k},$$o\`u
$\alpha^{i_1,\dots,i_k}(x)$ est un polyn\^ome homog\`ene de degr\'e $m$.\vskip 0.4cm

On peut munir $E$ de la graduation:  $|\alpha|=2m+k$, si $\alpha$ est un $k$-tenseur et $\alpha^{i_1\dots i_k}(x)$ est un polyn\^ome homog\`ene de degr\'e $m$.\vskip 0.4cm

En fait, l'espace $E$ est stable par le produit ext\'erieur $\wedge$ et par le crochet de Schouten $[~,~]_S$. Si on consid\`ere deux tenseurs $\alpha$ et $\beta$ de $E$ de degr\'e respectivement $|\alpha|=2m+k \ \ \hbox{et} \ \ |\beta|=2n+\ell.$ Alors, on a pour tout $\alpha,\beta,\gamma\in E$,
$$
\begin{aligned}
\alpha\wedge\beta=(-1)^{|\alpha||\beta|}\beta\wedge\alpha \ \ \hbox{et} \ \
\alpha\wedge(\beta\wedge\gamma)=(\alpha\wedge\beta)\wedge\gamma\end{aligned}.$$
Le produit ext\'erieur $\wedge$ est alors un produit commutatif et associatif gradu\'e de $E$. Il est de degr\'e $0$ puisque \ $|\alpha\wedge\beta|=2(n+m)+(k+\ell)=|\alpha|+|\beta|$.\vskip 0.2cm

On obtient que $(E,\wedge)$ est une alg\`ebre commutative et associative gradu\'ee. \vskip 0.2cm

D'autre part, le crochet de Schouten est de degr\'e $-3$ dans $E$ puisque \ $|[\alpha,\beta]_S|=2(n+m-1)+(k+\ell-1)=|\alpha|+|\beta|-3$.\vskip 0.2cm

On consid\`ere, alors, l'espace $E[3]$ et la graduation $deg(\alpha)=|\alpha|-3=2m+k-3$.\vskip 0.2cm

Sur $E[3]$, le crochet de Schouten sera de degr\'e $0$ et v\'erifie:
\begin{align*}
[\alpha,\beta]_S&=-(-1)^{deg(\alpha)deg(\beta)}[\beta,\alpha]_S\\
\hbox{et} \ (-1)^{deg(\alpha)deg(\gamma)}\big[[\alpha,\beta]_S,\gamma\big]_S+&(-1)^{deg(\beta)deg(\alpha)}\big[[\beta,\gamma]_S,\alpha\big]_S
+(-1)^{deg(\gamma)deg(\beta)}\big[[\gamma,\alpha]_S,\beta\big]_S=0.\end{align*}
Alors, $(E[3],[~,~]_S)$ est une alg\`ebre de Lie gradu\'ee. De plus, le crochet $[~,~]_S$ et le produit $\wedge$ v\'erifient l'identit\'e de Leibniz donn\'ee par:$$
[\alpha,\beta\wedge\gamma]_S=[\alpha,\beta]_S\wedge\gamma+(-1)^{|\beta|(|\alpha|-3)}\beta\wedge[\alpha,\gamma]_S.
$$
Donc, $(E,[~,~]_S,\wedge)$ est une $(0,-3)$-alg\`ebre. \vskip 0.2cm

Enfin, si on munit $E$ de la diff\'erentielle nulle $d=0$, on obtient que
$(E,[~,~]_S,\wedge,0)$ est une $(0,-3)$-alg\`ebre diff\'erentielle et on sait construire l'alg\`ebre \`a homotopie pr\`es associ\'ee \`a $E$.

\

\subsection{Exemple 2}

\

On consid\`ere $S(\mathbb{R}^d)$ l'alg\`ebre sym\'etrique de $\mathbb{R}^d$. Un \'el\'ement homog\`ene $f$ de $S(\mathbb{R}^d)$ est un polyn\^ome homog\`ene de degr\'e $k$. On munit $S(\mathbb{R}^d)$ de la graduation:
$$|f|=2k, \ \hbox{si \ $f$ \ est \ un \ polyn\^ome \ de \ degr\'e \ $k$}.$$
Sur $S(\mathbb{R}^d)$, on consid\`ere la multiplication usuelle  $\pt$ des polyn\^omes. Elle v\'erifie pour tout $f,g,h\in S(\mathbb{R}^d)$:
$$
\begin{aligned}
f\pt g&=(-1)^{|f||g|}g\pt f \ \ \hbox{et} \ \
f\pt(g\pt h)&=(f\pt g)\pt h\end{aligned}.$$
La multiplication $\pt$ est une op\'eration de degr\'e $0$ sur $S(\mathbb{R}^d)$ puisque \ $|f\pt g|=2(k+\ell)=|f|+|g|$.\vskip 0.2cm

D'autre part, on munit $S(\mathbb{R}^d)$ d'un crochet de Poisson d\'efini par:
$$
\{f,g\}=\sum_{i,j}\alpha^{i,j}(x)\ \partial_{x_i}(f).\partial_{x_j}(g),
$$
o\`u pour tout $i,j$, les coefficients $\alpha^{i,j}(x)$ est un polyn\^ome homog\`ene de degr\'e $m$.\vskip 0.2cm

Ce crochet est de degr\'e $2m-4$ puisque \ $|\{f,g\}|=2(m+k-1+\ell-1)=|f|+|g|+2m-4$.\vskip 0.2cm

On consid\`ere, alors, l'espace $S(\mathbb{R}^d)[-2m+4]$ et la graduation $deg(f)=|f|+2m-4$.\vskip 0.2cm

Sur $S(\mathbb{R}^d)[-2m+4]$, le crochet de Poisson sera de degr\'e $0$ et v\'erifie:
\begin{align*}
\{f,g\}&=-(-1)^{deg(f)deg(g)}\{g,f\}\\
\hbox{et} \ (-1)^{deg(f)deg(\gamma)}\{\{f,g\},\gamma\}+&(-1)^{deg(g)deg(f)}\{\{g,\gamma\},f\}
+(-1)^{deg(\gamma)deg(g)}\{\{\gamma,f\},g\}=0.\end{align*}

Alors, $\Big(S(\mathbb{R}^d)[-2m+4],\{~,~\}\Big)$ est une alg\`ebre de Lie gradu\'ee. De plus, le crochet $\{~,~\}$ et le produit $\pt$ v\'erifient l'identit\'e de Leibniz donn\'ee par:
$$
\{f,g\pt h\}=\{f,g\}\pt h+(-1)^{|g|(|f|+2m-4)}g\pt\{f,h\}.
$$
Donc, $\Big(S(\mathbb{R}^d),\{~,~\},\pt\Big)$ est une $(0,2m-4)$-alg\`ebre. \vskip 0.2cm

Enfin, si on munit $S(\mathbb{R}^d)$ de la diff\'erentielle nulle $d=0$, on obtient que
$\Big(S(\mathbb{R}^d),\{~,~\},\pt,0\Big)$ est une $(0,2m-4)$-alg\`ebre diff\'erentielle et on sait construire l'alg\`ebre \`a homotopie pr\`es associ\'ee \`a $S(\mathbb{R}^d)$.

\

\begin{rema}

\

Dans ces deux derniers exemples, pour tenir compte du degr\'e des fonctions polyn\^omes, on a muni les fonctions lin\'eaires $x\longmapsto x_i$ sur $\mathbb{R}^d$ du degr\'e $2$. En effet, puisqu'on veut que $\wedge$ restreint \`a $S(\mathbb{R}^d)\subset E$ ou $\pt$ sur $S(\mathbb{R}^d)$ soit la multiplication usuelle et puisque $x_ix_j=x_jx_i$, pour tout $i,j$, on doit donner \`a ces variables un degr\'e pair.

On veut g\'en\'eraliser, maintenant, ces exemples \`a un super-espace $\mathbb{R}^{p| q}$ pour lequel $S(\mathbb{R}^{p| q})\simeq S(\mathbb{R}^p)\bigotimes\bigwedge(\mathbb{R}^q)$ admet une base not\'ee \ $(x_{i_1}\dots x_{i_k}\xi_{j_1}\dots\xi_{j_r})$ \ avec $i_1\leq\dots\leq i_k$ et $j_1<\dots<j_r$.
\end{rema}

\

\subsection{Exemple 3}(Le crochet de Schouten sur le super-espace $\mathbb{R}^{p| q}$)

\

On reprend l'article \cite{[AMM]} dans le cas d'un super-espace $\mathbb{R}^{p| q}$ de coordonn\'ees $(x_1,\dots,x_p|\xi_1,\dots,\xi_q)$. Les variables $x_i$ commutent donc elles sont de degr\'e pair. On pose
$$
d^{\circ}x_i=2.
$$
Les variables $\xi_j$ anticommutent donc elles sont de degr\'e impair. On pose
$$
d^{\circ}\xi_j=1.
$$
Les champs de veceurs $\partial_{x_i}$ et $\partial_{\xi_j}$ associ\'es sont des d\'erivations gradu\'ees de l'alg\`ebre des polyn\^omes en $x_i$ et $\xi_j$:
$$
S(\mathbb{R}^{p| q})\simeq S(\mathbb{R}^p)\bigotimes\bigwedge(\mathbb{R}^q)
$$
de degr\'e
$$
d^\circ\partial_{x_i}=-2 \ \hbox{et} \ d^\circ\partial_{\xi_j}=-1.
$$
Comme dans \cite{[AMM]}, on modifie ces degr\'es pour obtenir la graduation naturelle de $T_{poly}(\mathbb{R}^{p| q})$. On va donc poser:
$$
|x_i|=2, \ |\xi_j|=1, \ |\partial_{x_i}|=-3 \ \hbox{et} \ |\partial_{\xi_j}|=-2.
$$
Le produit ext\'erieur usuel $\wedge_0$ \'etait d\'efini par:
$$
\partial_i\wedge_0\partial_j=\partial_i\otimes\partial_j-(-1)^{d^\circ(\partial_i)d^\circ(\partial_j)}\partial_j\otimes\partial_i,
$$
o\`u $\partial_i$ est soit $\partial_{x_i}$ soit $\partial_{\xi_i}$.\vskip0.15cm

Pour tenir compte du d\'ecalage, on pose maintenant
$$
\partial_i\wedge\partial_j=(-1)^{|\partial_i|}\partial_i\wedge_0\partial_j.
$$
On compl\`ete ce produit exr\'erieur avec le produit usuel, on aura
$$
x_i\wedge x_j=x_ix_j, \ x_i\wedge\xi_j=\xi_j\wedge x_i=\xi_jx_i, \ \xi_i\wedge\xi_j=-\xi_j\wedge\xi_i, \ x_i\wedge\partial_j=x_i\partial_j, \ \xi_i\wedge\partial_j=-\xi_i\partial_j.
$$
L'espace $T_{poly}(\mathbb{R}^{p| q})=\Big(S(\mathbb{R}^p)\bigotimes\bigwedge(\mathbb{R}^q)\Big)\bigotimes\Big(\bigwedge(\mathbb{R}^p)^*\bigotimes S(\mathbb{R}^q)^*\Big)$ est maintenant muni d'un produit $\wedge$ associatif et commutatif gradu\'e, de degr\'e $|\wedge|=0$.\vskip0.15cm

Un \'el\'ement $\alpha\in T_{poly}(\mathbb{R}^{p| q})$ s\'ecrit:
$$
\alpha=\sum_{1\leq i_1<\dots<i_n\leq p\atop 1\leq j_1\leq\dots\leq j_m\leq q}\sum_{1\leq s_1<\dots<s_r\leq p\atop 1\leq t_1\leq\dots\leq t_u\leq q}\alpha^{i_1\dots i_nj_1\dots j_m}_{s_1\dots s_rt_1\dots t_u}x_{s_1}\dots x_{s_r}\xi_{t_1}\dots \xi_{t_u}\partial_{x_{i_1}}\wedge\partial_{x_{i_n}}\wedge\partial_{\xi_{j_1}}\wedge\dots\partial_{\xi_{j_m}},
$$
ce qu'on notera simplement par:
$$
\alpha=\sum_{S,T,I,J}\alpha_{ST}^{IJ} \ x_S\xi_T\partial_{x_I}\wedge\partial_{\xi_J}.
$$
On a \ $|\alpha|=2\#S+\#T-3\#I-2\#J$.\vskip0.15cm

Sur $T_{poly}(\mathbb{R}^{p| q})$, la multiplication $\wedge$ est d\'efinie par:
$$
\alpha\wedge\beta=\sum_{S,T,I,J\atop A,B,C,D}(-1)^{\#B\#I}\alpha_{ST}^{IJ} \beta_{AB}^{CD}\ x_Sx_A\xi_T\xi_B\partial_{x_I}\wedge\partial_{x_C}\wedge\partial_{\xi_J}\wedge\partial_{\xi_D}.
$$
D'autre part, les champs de vecteurs $\mathfrak{X}(\mathbb{R}^{p| q})$ de la forme $X=\sum_i\tilde{X}^i(x,\xi)\partial_i$ sur le super-espace $\mathbb{R}^{p| q}$
ont \'et\'e \'etudi\'es par plusieurs auteurs \cite{[BB]}. Cet espace est muni d'un crochet d\'efini par:
$$[X,Y]=X\circ Y-(-1)^{(|X|-1)(|Y|-1)}Y\circ X.$$
On \'etend ce crochet au crochet de Schouten des multichamps de vecteurs en posant:
\begin{align*}&[X_1\wedge\dots\wedge X_n,Y_1\wedge\dots\wedge Y_m]_S=\\&\sum_{s=1}^n
\sum_{t=1}^m(-1)^{|X_s|(|X_{s+1}|+\cdots+|X_n|)+|Y_t|(|Y_1|+\dots+|Y_{t-1}|)}X_1\wedge\cdots\widehat{s}\dots\wedge X_n
\wedge[X_s,Y_t]\wedge Y_1\wedge\cdots\widehat{t}\dots\wedge Y_m
.\end{align*}(Voir par exemple \cite{[AAC1]}, \cite{[AIPP-B]}, \cite{[K]}, \cite{[Le]}...) \vskip0.2cm
Ceci revient \`a imposer une relation de Leibniz de la forme:
$$
[X,Y\wedge Z]=[X,Y]\wedge Z+(-1)^{|Y|(|X|-1)}Y\wedge[X,Z].
$$

En fait le degr\'e naturel sur $\mathfrak{X}(\mathbb{R}^{p| q})$ pour ce crochet est
$$deg(X)=|X|+1$$ puisque $|[~,~]|=1$. On trouve que $\mathfrak{X}(\mathbb{R}^{p| q})$ muni de $deg$ et $[~,~]$ est une alg\`ebre de Lie gradu\'ee.\vskip0.15cm

De m\^eme, si on pose
$$
deg(X_1\wedge\dots \wedge X_n)=|X_1\wedge\dots \wedge X_n|+1,
$$
alors, $\Big(T_{poly}(\mathbb{R}^{p| q})[-1], [~,~]_S\Big)$ est une alg\`ebre de Lie gradu\'ee (pour $deg$).
\vskip0.15cm

On peut, comme dans \cite{[AMM]}, r\'e\'ecrire ces op\'erations ainsi:\vskip0.15cm

\noindent Pour deux champs de vecteurs $X=\sum_i\tilde{X}^i(x,\xi)\partial_i$ et $Y=\sum_j\tilde{Y}^j(x,\xi)\partial_j$, on a:
$$
[X,Y]=\sum_{i,j}\tilde{X}^i(x,\xi)\partial_i(\tilde{Y}^j(x,\xi))\partial_j-(-1)^{deg(X)deg(Y)}\tilde{Y}^j(x,\xi)\partial_j(\tilde{X}^i(x,\xi))\partial_i.
$$
Pour deux tenseurs $$\alpha=\sum_{i_1,\dots,i_n}\tilde{\alpha}^{i_1\dots i_n}(x,\xi)\partial_{i_1}\wedge\dots\wedge\partial_{i_n}=\sum_I\tilde{\alpha}^I(x,\xi)\partial_I$$
et
$$\beta=\sum_{j_1,\dots,j_m}\tilde{\beta}^{j_1\dots j_m}(x,\xi)\partial_{j_1}\wedge\dots\wedge\partial_{j_m}=\sum_J\tilde{\beta}^J(x,\xi)\partial_J,$$on a:
\begin{align*}
&[\alpha,\beta]_S=\\&\sum_{I,J}\Big(\sum_{s=1}^n(-1)^{|\partial_{i_s}|(|\partial_{i_1}|+\dots+|\partial_{i_{s-1}}|)+(|\tilde{\beta}^J|+1)(|\partial_I|-|\partial_{i_s}|)}
\tilde{\alpha}^I(x,\xi)\partial_{i_s}(\tilde{\beta}^J(x,\xi))\partial_{i_1}\wedge\dots\widehat{s}\dots\wedge\partial_{i_n}\wedge\partial_J\\&-\sum_{t=1}^m
(-1)^{|\partial_{j_t}|(|\partial_{j_1}|+\dots+|\partial_{j_{t-1}}|)+(|\tilde{\beta}^J|+|\partial_{j_t}|+1)(|\alpha|+1)}
\tilde{\beta}^J(x,\xi)\partial_{j_t}(\tilde{\alpha}^I(x,\xi))\partial_I\wedge\partial_{j_1}\wedge\dots\widehat{t}\dots\wedge\partial_{j_m}\Big).
\end{align*}
Ce qu'on peut l'\'ecrire comme dans \cite{[AMM]}:
$$
[\alpha,\beta]_S=(-1)^{|\alpha|+1}\alpha\bullet\beta-(-1)^{|\alpha|(|\beta|+1)}\beta\bullet\alpha,
$$
avec
$$
\alpha\bullet\beta=\hskip-0.15cm\sum_{s=1}^n(-1)^{|\alpha|+1+|\partial_{i_s}|(|\partial_{i_1}|+\dots+|\partial_{i_{s-1}}|)+(|\tilde{\beta}^J|+1)(|\partial_I|-|\partial_{i_s}|)}
\tilde{\alpha}^I(x,\xi)\partial_{i_s}(\tilde{\beta}^J(x,\xi))\partial_{i_1}\wedge\dots\widehat{s}\dots\wedge\partial_{i_n}\wedge\partial_J.
$$
Finalement, on a une relation de Leibniz de la forme:
$$
[\alpha,\beta\wedge\gamma]_S=[\alpha,\beta]_S\wedge\gamma+(-1)^{(deg(\beta)-1)deg(\alpha)}\beta\wedge[\alpha,\gamma]_S.
$$
On en d\'eduit que $\Big(T_{poly}(\mathbb{R}^{p| q}),\wedge, [~,~]_S, d=0\Big)$ est une $(0,1)$-alg\`ebre diff\'erentielle gradu\'ee puisque
$\Big(T_{poly}(\mathbb{R}^{p| q}),\wedge\big)$ est une alg\`ebre commutative gradu\'ee et $\Big(T_{poly}(\mathbb{R}^{p| q})[-1], [~,~]_S\Big)$
est une alg\`ebre de Lie gradu\'ee telle que $\wedge$ et $[~,~]_S$ v\'erifie l'identit\'e de Leibniz.

\

\subsection{Exemple 4}(Le crochet de Poisson sur le super-espace $\mathbb{R}^{p| q}$)

\

Reprenons le super-espace $\mathbb{R}^{p| q}$ avec les variables $x_i$ et $\xi_j$ de degr\'es:
$$
|x_i|=2 \ \hbox{et} \ |\xi_j|=1.
$$
L'alg\`ebre $S(\mathbb{R}^{p| q})=S(\mathbb{R}^p)\bigotimes\bigwedge(\mathbb{R}^q)$ devient une alg\`ebre commutative gradu\'ee pour ce degr\'e:
$$
f\pt g=(-1)^{|f||g|}g\pt f.
$$
Cette alg\`ebre admet des d\'erivations gradu\'ees $\partial_{x_i}$ et $\partial_{\xi_j}$ avec $|\partial_{x_i}|=-2$ et $|\partial_{\xi_j}|=-1$. Suivant \cite{[AIPP-B]}, on d\'efinit un crochet de Poisson homog\`ene de degr\'e $m\in\mathbb{Z}$ sur $S(\mathbb{R}^{p| q})$ en posant:
$$
\{f,g\}=(-1)^{m|f|}\sum_{i,j}(-1)^{|\partial_j|(|f|+|\partial_i|)}\omega^{ij}(x,\xi)\partial_i(f)\partial_j(g),
$$
o\`u $\omega^{ij}$ est une fonction polyn\^ome. Dire que $\{~,~\}$ est un crochet de Poisson, c'est \`a dire que $\Big(S(\mathbb{R}^{p| q})[-m],\{~,~\}\Big)$ est une alg\`ebre de Lie gradu\'ee. Ceci se traduit par trois conditions sur les $\omega^{ij}$:
\begin{align*}&{\bf (i):} \ \ \ |\omega^{ij}|+|\partial_i|+|\partial_j|=m, \ \forall i,j,\\&
{\bf (ii):} \ \ \ \omega^{ij}(x,\xi)=(-1)^{|\partial_i||\partial_j|+m+1}\omega^{ji}(x,\xi), \ \forall i,j,\\&
{\bf (iii):} \ \ \sum_k(-1)^{|\partial_\ell|(m+|\partial_i|)}\omega^{\ell,k}(x,\xi)\partial_k\big(\omega^{ji}(x,\xi)\big)
+(-1)^{|\partial_j|(m+|\partial_\ell|)}\omega^{j,k}(x,\xi)\partial_k\big(\omega^{i\ell}(x,\xi)\big)\\&\hskip1.6cm
+(-1)^{|\partial_i|(m+|\partial_j|)}\omega^{i,k}(x,\xi)\partial_k\big(\omega^{\ell j}(x,\xi)\big)=0, \ \forall i,j,\ell.
\end{align*}
Enfin, la relation de Leibniz est bien v\'erifi\'ee pour $\pt$ et $\{~,~\}$. Elle s'\'ecrit:
$$
\{f,g\pt h\}=\{f,g\}\pt h+(-1)^{|g|(|f|+m)}g\pt\{f,h\}.
$$
On obtient, donc, que $\Big(S(\mathbb{R}^{p| q}),\pt,\{~,~\},d=0\Big)$ est une $(0,m)$-alg\`ebre diff\'erentielle gradu\'ee.

\vskip9cm


\end{document}